\documentclass[12pt]{article}

\usepackage{amsmath,amsthm,amsfonts,amssymb,epsf,epsfig}
\evensidemargin
.25in \textheight 8.5in \topmargin 0in 
\title{Contact Structures, $\sigma$-Confoliations, and Contaminations in 3-Manifolds}
\author{ Ulrich Oertel \thanks{Research supported by  Max Planck Institute, Bonn; the National Science Foundation, grant DMS 9803293}
\hbox{\ \ }\&\hbox{\ }  
Jacek \'Swi\c atkowski \thanks{Research supported by  Max Planck Institute, Bonn; the National Science Foundation; and
the Polish Scientific Research Committee (KBN), grant 5 P02A 035 20}}
\date{Revised April, 2003}

\newtheorem{thm}{Theorem}[section] 
\newtheorem{lemma}[thm]{Lemma}
\newtheorem{fact}[thm]{Fact}
\newtheorem{corollary}[thm]{Corollary}
\newtheorem{proposition}[thm]{Proposition}

\theoremstyle{definition}
\newtheorem{defn}[thm]{Definition}

\newtheorem{ex}[thm]{Example}

\newtheorem{remark}[thm]{Remark}
\theoremstyle{remark}
\begin{document}
\maketitle
%\tableofcontents

\def\reals{\mathbb R}
\def\rationals{\mathbb Q}
\def\complex{\mathbb C}
\def\naturals{\mathbb N}
\def\integers{\mathbb Z}

\def\proj{P}
\def\hyp {\hbox {\rm {H \kern -2.8ex I}\kern 1.15ex}}
\def\intr{\text{int}}
\def\inter{\ \raise4pt\hbox{$^\circ$}\kern -1.6ex}
\def\Cal{\cal}
\def\from{:}
\def\inverse{^{-1}}
\def\Max{\text{Max}}
\def\Min{\text{Min}}
\def\embed{\hookrightarrow}
\def\Genus{\text{Genus}}
\def\Z{Z}
\def\X{X}
\def\roster{\begin{enumerate}}
\def\endroster{\end{enumerate}}
\def\definition{\begin{defn}}
\def\enddefinition{\end{defn}}
\def\subhead{\subsection\{}
\def\theorem{thm}
\def\endsubhead{\}}
\def\head{\section\{}
\def\endhead{\}}
\def\example{\begin{ex}}
\def\endexample{\end{ex}}
\def\ves{\vs}
\def\mZ{{\mathbb Z}}
\def\M{M(\Phi)}
\def\bdry{\partial}
\def\hop{\vskip 0.15in}
\def\trip{\vskip 0.09in}

\centerline{Rutgers University, Newark;  Wroc\l aw University, Wroc\l aw }

\abstract 
We propose in this paper a method for studying contact structures in 3-manifolds
by means of branched surfaces.  We explain what it means
for a contact structure to be carried by a branched surface embedded in
a 3-manifold.  To make the transition from contact structures to branched surfaces, we first
define  auxiliary objects called $\sigma$-confoliations and pure contaminations, both
generalizing contact structures. We study various deformations 
of these objects and show that the $\sigma$-confoliations and pure contaminations obtained by suitably modifying  
a contact structure remember the contact structure up to isotopy.

After defining tightness for all pure contaminations in a natural way, generalizing the definition of tightness for contact
structures, we obtain some conditions on (the embedding of) a branched surface in a 3-manifold sufficient to
guarantee that any pure contamination carried by the branched surface is
tight.  We also find conditions sufficient to prove that a branched surface carries only overtwisted (non-tight)
contact structures.

Our long-term goal in developing these methods is twofold:  Not only do we want to study tight contact structures
and pure contaminations, but we also wish to use them as tools for studying 3-manifold topology.  These
structures can exist in manifolds having either infinite or finite fundamental group, which suggests that they might be especially
difficult to use effectively, but which also makes them especially attractive.
\endabstract

\section{Introduction }\label{Introduction}

We will usually assume that $M$ is a closed orientable or oriented 3-manifold.

\definition  Suppose a field $\xi$ of planes on an oriented
3-manifold
$M$ is locally defined by a form
$\alpha$.  The field is integrable, and defines a {\it foliation} if $\alpha\wedge
d\alpha=0$.  It defines a {\it positive contact structure} ({\it positive confoliation}) if 
$\alpha\wedge d\alpha>0$ ($\alpha\wedge d\alpha\ge 0$).  It defines a {\it negative contact structure} ({\it negative
confoliation}) if 
$\alpha\wedge d\alpha<0$ ($\alpha\wedge d\alpha\le 0$).  A positive (negative) confoliation is {\it contact at a point} if 
$\alpha\wedge d\alpha>0$ ($\alpha\wedge d\alpha<0$) at that point.
\enddefinition

A (positive) contact structure is a
field which is nowhere locally integrable due to the fact that the planes of the
field locally look something like the tangent planes to the blades of a 
propeller (with the appropriate sense).  A positive confoliation is a similar plane
field, except that the limiting case of integrable fields is allowed, as is a
mixture of integrable and positive non-integrable behavior.   In this paper our working definition of a confoliation is often in
terms of local charts called ${\cal B}$-charts and ${\cal C}$-charts, see Section \ref{Reduction}.

The monograph \cite{ET:Confoliations} of Yakov Eliashberg and 
William Thurston made clear the close relationship between 
contact structures and foliations.  If further evidence is needed, 
the recent papers \cite{HKM2:Convex,HKM:TightContact}
confirm this close relationship. The work of Eliashberg and Thurston 
led us to consider the possibility of using branched
surfaces to study contact structures.  To make this approach possible, 
one must show that a contact structure can be
modified such that it is ``carried" by a branched surface, 
and one must further show that the modification can be done in
such a way that the contact structure is not lost, or that 
it is recoverable from the modified object. As a first step in this
direction  we will show that if a contact structure is deformed slightly, 
to a special kind of confoliation called $\sigma$-confoliation, 
then the contact structure is indeed remembered.

\definition  A {\it $\sigma$-confoliation} is a
confoliation $\xi$ in a closed
3-manifold $M$ for which there exists
a smooth compact surface $F$ in $M$ such that:
\item{(1)} $F$ is everywhere tangent to $\xi$, 
\item{(2)} $\xi$ is contact in the complement $M- F$,
\item{(3)} each connected component of $F$ has nonempty boundary
(equivalently, $F$ admits a handle decomposition consisting of $0$-handles
and $1$-handles only).

\trip

\noindent
$F$ is then called the {\it maximal integral surface}
of the $\sigma$-confoliation $\xi$.
\enddefinition

In the remainder of the paper, we will usually assume that all 
contact structures, confoliations, and contaminations (yet
to be defined) are positive.

In studying modifications of contact structures into
confoliations we distinguish the following classes of deformations.

\definition\label{SigmaDeformationDef} 
A {\it deformation of contact structures} or {\it pure deformation} is
a smooth deformation $\xi_t:t\in[a,b]$ of plane fields in $M$,
such that for each $t\in[a,b]$ the plane field
$\xi_t$ is a contact structure.  Similarly we define a {\it deformation of
confoliations} as a smooth 1-parameter family of confoliations.

A  deformation of confoliations (contact structures) $\xi_t$ is an {\it isotopy of confoliations (contact structures)} if
there is a smooth 1-parameter family $\phi_t$ of diffeomorphisms of $M$ such that $\phi_0$ is the identity
and $(\phi_t)_*(\xi_0)=\xi_t$.

A {\it  $\sigma$-deformation} is 
a smooth deformation $\xi_t:t\in[a,b]$ of plane fields in $M$,
such that for each $t\in[a,b]$ the plane field
$\xi_t$ is a $\sigma$-confoliation, and moreover the maximal integral surfaces 
$F_t$ of $\xi_t$ behave well
with respect to the parameter $t$ in the following sense.
There exists a 1-parameter family $h_t:t\in[a,b]$ of diffeomorphisms
of $M$, a partition of the interval
$[a,b]$ into a finite number of nontrivial disjoint consecutive
subintervals $I_1,I_2,\dots,I_k$,
and a sequence $\Sigma_1,\Sigma_2,\dots,\Sigma_k$
of compact smooth surfaces in $M$, such that:
\item{(1)} for each $t\in \intr(I_j)$ we have $h_t(F_t)=\Sigma_j$;
\item{(2)} for each $1\le j\le k-1$ the surfaces $\Sigma_j$ and 
$\Sigma_{j+1}$ are distinct, and one of them
is contained in the other;
\item{(3)} if $\Sigma_j\subset\Sigma_{j+1}$ (or conversely
$\Sigma_{j+1}\subset\Sigma_j$), then
the closure $cl(\Sigma_{j+1}-\Sigma_j)$
(cl($\Sigma_j-\Sigma_{j+1}$) respectively) is
a compact surface with piecewise smooth boundary;
\item{(4)} 
each component of the larger of the surfaces
$\Sigma_j$ and $\Sigma_{j+1}$ is obtained from the smaller by 
attaching 1-handles and adding 0-handles.
\enddefinition

In terms of these definitions, Gray's Theorem can be stated as follows:
If $\xi_0$ and $\xi_1$ are connected by a deformation of contact structures, then
they are isotopic.

In Sections \ref{Reduction}, \ref{Splitting} and \ref{SplitToCarry} we show how to construct various $\sigma$-deformations,
both the ones that decrease and the ones that increase  maximal
integral surfaces. Our main result concerning $\sigma$-deformations
is the following generalization of Gray's Theorem. 

\begin{thm}\label{SigmaDeformationThm} {\bf $\sigma$-Deformation Theorem.}  
Any two contact structures on a closed 3-manifold
$M$ which are connected  by a $\sigma$-deformation are isotopic.  
\end{thm}

Theorem \ref{SigmaDeformationThm} implies that a contact structure is
``remembered" or determined by a $\sigma$-confoliation connected to the contact
structure by a $\sigma$-deformation.  To see this, suppose that a $\sigma$-confoliation $\xi$ 
is connected to contact structures $\xi_0$ and $\xi_1$ by $\sigma$-deformations. Combining the
$\sigma$-deformations, we obtain a $\sigma$-deformation from $\xi_0$ to $\xi_1$, and therefore by
the theorem $\xi_0$ and $\xi_1$ are isotopic.

The statement of Theorem \ref{SigmaDeformationThm} should be compared with Proposition 2.3.2 of \cite{ET:Confoliations},
which implies that smooth deformations through arbitrary confoliations from one contact structure to another do not preserve the
isomorphism class of the contact structures.

 One of our most ambitious goals is to understand phenomenon of tightness 
of contact structures using
$\sigma$-confoliations and branched surfaces. Towards the end of understanding tightness via 
$\sigma$-confoliations and branched surfaces, we
give the following definition which extends the notion of tightness 
to all $\sigma$-confoliations. 

\begin{definition}\label{SigmaConfoliationTightDef}
 A $\sigma$-confoliation $\xi$ is $\it tight$, if there is
no
2-dimensional disc $D$ smoothly embedded in $M$, such that the boundary
$\partial D$ is tangent to $\xi$ and disjoint from the maximal integral
surface $F$ of $\xi$, and $D$ itself is transverse to $\xi$ at $\partial D$.
\end{definition}

Since we use Theorem \ref{SigmaDeformationThm} to identify a 
contact structure with a $\sigma$-confoliation if the two
structures are related by a
$\sigma$-deformation, the following theorem is needed to show that 
the above definition is consistent with the definition
for contact structures.

\begin{thm} \label{TightInvarianceTheorem}  
Tightness of $\sigma$-confoliations is an invariant of
$\sigma$-deformation.
\end{thm}

\hop

The statement of Theorem \ref{TightInvarianceTheorem} should be compared with Proposition 3.6.2 of \cite{ET:Confoliations},
which shows that smooth deformations through arbitrary confoliations do not preserve the property of tightness.

Now we introduce branched surfaces, and we explain what it means for a 
$\sigma$-confoliation to be carried by a branched
surface.  The following ideas, involving branched surfaces and 
contaminations will be explained more fully in Section
\ref{Splitting}.

\begin{figure}[ht]
\centering
\scalebox{1.0}{\includegraphics{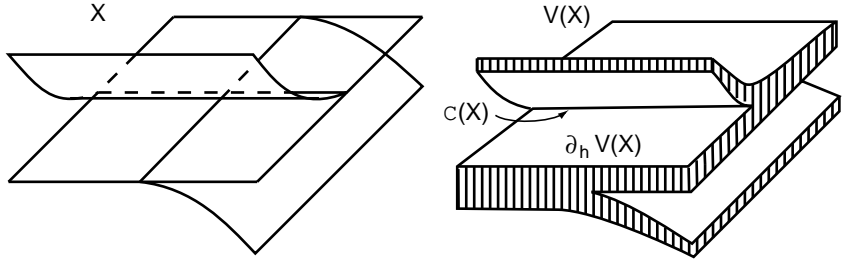}}
\caption{\small Local model for branched surface and fibered neighborhood.}
\label{ContactBranchedModel}
\end{figure}

A closed {\it branched surface} $B$  with generic branch
locus is a space with smooth structure locally modelled on the space $X$ 
shown in Figure 1.  (A neighborhood of each point of
$B$ is diffeomorphic to a neighborhood of a point in $X$.)   
The {\it branch locus} $K$ of
$B$ is the union of points of $B$ none of whose neighborhoods is a manifold.  
The branch locus is (homeomorphic to) a cell complex of dimension
$1$, but it can also be regarded as the image of the {\it branch curves} 
${\cal C}(B)$
immersed in $B$.  The {\it sectors} of $B$ are the completions (in
a path metric coming from a Riemannian metric on $B$) of the components of
$B- K$.  The sense of branching at the boundary of a sector is described using the following convention.  If sectors
$W,X,Y$ are adjacent along an arc $\gamma$ of branch locus, and if
$W\cup Y$ and $X\cup Y$ are smooth, we say that branching along the arc
$\gamma\subset \bdry Y$ is {\it inward} for $Y$ and {\it outward} for $X$.  We indicate the sense of branching using a
transverse arrow pointing into
$Y$.   If $B$ is embedded in a
3-manifold,  then a {\it fibered neighborhood} $V(B)$ of $B$ in $M$ is a
closed regular neighborhood of $B$ foliated by interval fibers, as shown in
Figure 1, with the frontier $\bdry V(B)$ cut by {\it cusp curves} 
${\cal C}(B)$ into smooth surfaces called the {\it components of horizontal
boundary} and denoted $\bdry_h V(B)$.  The reader should not confuse $\bdry_hV(B)$ with $\bdry V(B)$, which is not smooth. 
We note that $\bdry_hV(B)$ is not embedded in $M$ or in $\bdry V(B)$, and when we say that a subspace is contained in $\bdry_hV(B)$
we mean that the embedding of the subspace in $M$ factors through an embedding in $\bdry_hV(B)$.   There is a projection map
$\pi\from V(B)\to V(B)/\sim$ which maps $V(B)$ to a quotient space 
in which fibers of $V(B)$ are collapsed to points.  Notice that $V(B)/\sim$ can be identified with $B$,  
and that we identify the branch curves
${\cal C}(B)$ immersed in $B$ with the cusp curves, 
also denoted ${\cal C}(B)$, which are embedded in $\bdry V(B)$.

More detailed definitions of branched surfaces can be found in other papers, 
e.g. \cite{DGUO:EssentialLaminations, LMUO:Nonnegative, 
UO:MeasuredLaminations}, or in the original source for
branched manifolds, R. Williams, \cite {RW:BranchedSurfaces}.

\begin{definition} \label{ContaminationDefB}
Suppose $B$ is a branched surface embedded in a 3-manifold $M$.  
Suppose a smooth plane field $\xi$ is defined on $V(B)$ such that:

\item{(1)} $\xi$ is transverse to the interval fibers of $V(B)$;

\item{(2)} $\xi$ is tangent to $\bdry_hV(B)$;

\item{(3)} $\xi$ is a confoliation in $V(B)$.
\trip

\noindent Then $\xi$ is called a {\it contamination carried by $B$}. 
\end{definition}

\begin{figure}[ht]
\centering
\scalebox{1.0}{\includegraphics{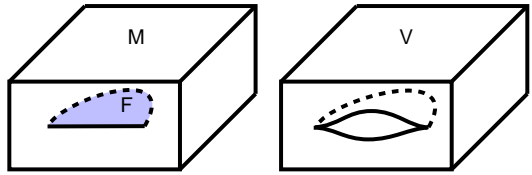}}
\caption{\small Cutting a $\sigma$-confoliation.}
\label{ContactSplit}
\end{figure}

Suppose now that $\xi$ is a $\sigma$-confoliation in $M$.  
Then we can {\it cut} $M$ on the maximal integral surface $F$ 
to obtain $V$, see Figure 2, a {\it manifold with inward cusps}. 
This is a manifold with all the characteristics of a
neighborhood $V(B)$ of a branched surface, except there is no foliation 
by interval fibers.  The manifold has a smooth
structure in a sense which will be made precise in Section \ref{Splitting}.  The boundary also has a smooth structure, 
but has inward cusps at a 1-submanifold ${\cal C}(V)$ of the boundary.  We
use $\bdry_h V$, the {\it horizontal boundary}, to denote the components of
$\bdry V$ cut on cusp curves.  As before, we observe that $\bdry_h V$ is not embedded in $M$, and that
when we say that a subspace of $M$ is contained in $\bdry_hV(B)$
we mean that the embedding of the subspace in $M$ factors through an embedding in $\bdry_hV(B)$.

If the manifold with inward cusps $V$ was obtained by cutting on $F$, 
then $\bdry_hV$ consists of two copies of the integral surface
$F$ embedded in $\bdry V$.   In this case the manifold
$V$ can be embedded smoothly in
$M$, such that the complement is topologically a product
$\hbox{int}(F)\times(0,1)$.  The $\sigma$-confoliation
$\xi$ induces on $V$ a smooth plane field, which we denote
$\xi'$.  This is the pure contamination in $M$ induced by
the $\sigma$-confoliation $\xi$.  Regarding $V$ as a subspace of $M$, 
the contamination $\xi'$
is a partial plane field defined only on $V\subset M$.  The operation converting $\xi$ to $\xi'$ is called {\it
cutting} $\xi$ on the maximal integral surface.  
 
More generally:

\begin{definition} \label{ContaminationDef}
A {\it (pure) contamination} in $M$ is a smooth plane
field $\xi$ defined on a manifold $V$ with inward cusps 
embedded smoothly in $M$, with the property that $\xi$ is
tangent to $\bdry_hV$ and $\xi$ is a (contact structure) confoliation in $\intr (V)$.  
The manifold $V$ is called the {\it support} of the contamination $\xi$. 
\end{definition}

 We note that just as a confoliation is a
hybrid object obtained from foliations and contact structures, so a contamination is a hybrid object obtained from smooth
laminations and contact structures.

We do not use arbitrary contaminations in this paper, only pure contaminations and certain other contaminations called
$\sigma$-contaminations. In Section \ref{Splitting} we will give precise
definitions of $\sigma$-contaminations, $\sigma$-deformations of
$\sigma$-contaminations, and pure deformations of pure contaminations.  Roughly speaking, a $\sigma$-contamination with
support $V$ is a contamination such that the integral surface consists of $\bdry_h V$ together with a compact surface
whose interior lies in the interior of $V$, and a $\sigma$-deformation is a deformation through $\sigma$-contaminations. 
Similarly, a pure deformation is a deformation through pure contaminations.  It should be clear that
$\sigma$-confoliations are
$\sigma$-contaminations.

The operation of {\it splitting} a contact structure $\xi$ converts the contact structure to a pure contamination $\xi''$
in two steps.  First a $\sigma$-deformation converts $\xi$ to a $\sigma$-confoliation $\xi'$, then cutting on the maximal
integral surface converts
$\xi'$ to a pure contamination $\xi''$.  

We would like to represent an arbitrary contact structure $\xi$ in $M$ 
by a pure contamination obtained from $\xi$ by means of splitting, with additional property that this contamination is
carried by some branched surface
$B\embed M$. The following theorem shows this is possible. 

\begin{thm}\label{SplitToCarryThm} Let $\xi$ be a contact structure 
defined on a closed 3-manifold $M$ and let $\Omega$
be a smooth 1-foliation transverse to $\xi$.  There is a 
$\sigma$-deformation $\xi_t$, $0\le t\le 1$, from $\xi=\xi_0$ to a $\sigma$-confoliation $\xi_1$ 
such that $\xi_t$ is transverse to $\Omega$ at all times $t$ and the
maximal integral surface
$F$ for $\xi_1$ cuts all leaves of $\Omega$ into intervals.  
It follows that cutting $\xi_1$ on $F$ yields a manifold $V$ with inward cusps with
a foliation by intervals induced by $\Omega$, and thus $V$ is
isomorphic to the fibered neighborhood of a branched surface $B\embed M$.
Moreover, the resulting pure contamination supported by $V$ is clearly
carried by this branched surface and is a splitting of $\xi_0$.
\end{thm}

\hop
We would like to be able to recover a contact structure from a pure contamination
obtained from the contact structure by a splitting. This can be done using
a (roughly) opposite operation which we call pinching. Pinching
can be applied to any pure contamination $\xi$ in $M$ supported by a manifold $V$ 
satisfying the following condition: 
for any connected component $U$ of $M- V$ the frontier $\bdry U$ 
consists of two homeomorphic components
$F$ of the horizontal boundary $\bdry_hV$ (where $\partial F\ne\emptyset$)
and $U$ itself is homeomorphic
to the product $\hbox{int}(F)\times(0,1)$ (intuitively, $U$ is homeomorphic
to a thickened $F$ with outward cusps at $\partial F$). 
If this condition is satisfied, we will say that
contamination $\xi$ has {\it product complementary pieces}.
Note that any pure contamination obtained by a splitting of a contact structure
clearly has product complementary pieces.

{\it Pinching} a pure contamination $\xi$ with product complementary pieces consists of
the following steps
(which will be described more precisely in Section \ref{Pinching}). First, we identify
the two components $F$ in the frontier of any complementary
component $U$ via a diffeomorphism equal to the identity on $\partial F$
(intuitively, 
we collapse $U$ into a single integral surface $F$ in $M$). 
We will call this operation a {\it glueing} and note that it depends not only
on a glueing map, but also on the choice of smooth structure at the glueing locus
(this will be discussed in detail in Section \ref{Pinching}).
Glueing gives us $\xi'$, a sort of $\sigma$-confoliation, which is in
general not smooth but only continuous (as a plane field in $M$)
at the glueing locus, i.e.
at the integral surfaces $F$ obtained by collapsing the components $U$.
We call such an object a {\it semi-smooth $\sigma$-confoliation}.
We then modify $\xi'$ to a contact structure, by simultaneously smoothing it and
eliminating its maximal integral surface. This can be done by means of a
{\it semi-smooth $\sigma$-deformation}.  Semi-smooth $\sigma$-deformations will be defined in Section \ref{Pinching} and
their existence will be proved there. 

The operation inverse to a splitting is an example of a pinching, since we first
reglue and then apply the reverse $\sigma$-deformation (which is of course a special case 
of a semi-smooth $\sigma$-deformation). 

The following theorem shows that a contamination obtained from a contact structure by splitting remembers the contact
structure.

\begin{thm}\label{ContactizationThm} 
(1) Let $\xi$ be a pure contamination in
a 3-manifold $M$ with product complementary pieces. Then any two contact
structures on $M$ obtained from $\xi$ by pinching are isotopic.

\noindent
(2) If $\xi$ is a pure contamination obtained from
a contact structure $\eta$ on $M$ by splitting (e.g. as in Theorem \ref{SplitToCarryThm}), 
then any contact structure obtained from $\xi$ by pinching is isotopic to $\eta$.

\end{thm}

The above theorem is a special case of a more general result, Lemma \ref{TwoSmoothingsLem}. In Section \ref{Pinching}, we
define splitting and pinching in the context of contaminations.  Each splitting or pinching of a pure contamination yields
another  pure contamination.  The splitting and pinching operations on pure contaminations generalize the splitting and pinching
operations defined above for contact structures.  We then prove that any two
contact structures connected by any sequence of splittings and pinchings
(of pure contaminations) are isotopic.

\hop
\begin{definition} Let $\xi$ be a contact structure on a 3-manifold $M$,
and let $B$ be a branched surface embedded in $M$. We say that
$B$ {\it carries} $\xi$ if $B$ carries a pure contamination obtained
from $\xi$ by a splitting.
\end{definition}

It is not clear which branched surfaces $B$ embedded
in 3-manifolds carry contact structures. An obvious necessary condition
is that $B$ has product complementary pieces (which means that its
fibered neighborhood $V(B)$ has product complementary pieces).
In the separate paper 
\cite {UOJS:ContaminationCarrying}, we give a sufficient condition for a branched
surface to carry a positive pure contamination. Together with the requirement
of product complementary pieces, this gives a sufficient condition to carry
a positive contact structure. On the other hand, a branched
surface with product complementary pieces may carry two or more
nonisotopic contact structures. Therefore in general one cannot recover
a contact structure from a branched surface which carries it.

It turns out however, that one can impose conditions on a branched surface $B$ so that any contact structure carried is
overtwisted.  We are able to distinguish such branched surfaces using
``pseudo-transversal knots."

\definition 
Let $B$ be a branched surface embedded in $M$ with
product complementary pieces. A null-homologous knot $\kappa$ in $M$
is a {\it pseudo-transversal} for $B$ if it satisfies the following
conditions: (1) $\kappa$ intersects $B$ transversely away from the branch locus, (2) the intersection
of $\kappa$ with each complementary component $U$ of $M- B$
is a collection of unknotted segments monotonically passing
from one to the other component of the frontier $\partial U$
(so that these segments form a braid in $U$).
\enddefinition

It is possible to define a {\it self linking number} $l(\kappa)$
of a knot $\kappa$ pseudo-transversal to $B$, 
in a way analogous to the self-linking (or the Thurston-Bennequin invariant) of
a null-homologous knot transverse to a contact structure (see Section
\ref{Pseudotransversal} for more details).
If a branched surface $B$ carries a contact structure $\xi$ then
$\kappa$ determines (uniquely up to
transverse isotopy) a knot $\kappa'$ transverse to $\xi$,
with the same self-linking. Transversal knots in tight contact
structures are known to satisfy the so called Bennequin inequality 
$$
l(\kappa')\le -\chi(S),
$$
where $\chi(S)$ is the Euler characteristic of any Seifert surface $S$
for $\kappa'$, see \cite{YE:LegendrianKnots}. By what was said above, the same inequality holds
for pseudo-transversal knots $\kappa$, for any branched surface $B$
that carries a tight contact structure. Thus we get the following.

\begin{thm}\label{PseudoTransversalThm}
Let $B$ be a branched surface with product complementary
pieces embedded in a 3-manifold $M$. Suppose that there is a knot in
$M$ which is pseudo-transversal to $B$ and which violates the Bennequin inequality. 
Then any contact structure in $M$ carried by $B$ is overtwisted.
\end{thm}

It is clear that a branched surface $B$ uniquely determines the homotopy
type (as a plane field) of any contact structure it carries.
In view of the fact that an overtwisted contact structure
is uniquely determined (up to isotopy) by its homotopy type (see \cite {YE:OvertwistedClassification}),
a branched surface $B$ as in \ref{PseudoTransversalThm} carries at most one
contact structure (up to isotopy) and this structure is overtwisted.

One of the most important goals of our research program is to prove results about 3-manifold topology using contact
structures and contaminations.  It is well-known that among contact structures on a 3-manifold, only tight ones can
contain information about the topology of the 3-manifold.  It is therefore necessary to find a suitable definition of
tightness for pure contaminations.  The definition proposed by Eliashberg and Thurston for tightness of confoliations
gives a clear indication for defining tightness of pure contaminations:

\begin{definition} Suppose $\xi$ supported by $V\embed M$ is a pure contamination in $M$.  Then a disc
$D\embed M$ is an {\it overtwisting disc} for $\xi$ if the boundary
$\partial D$ is tangent to $\xi$, $D$ itself is transverse to $\xi$ at $\partial D$, and in case $\bdry D\subset\bdry_hV$,
$\bdry D$ does not bound a disc in $\bdry_hV(B)$.  The pure contamination $\xi$ is {\it tight} if there are no overtwisting discs
for $\xi$.
\end{definition}

We will show that in case the pure contamination $\xi$ represents a contact structure, this definition implies tightness according to
previous  definitions.  There are fairly obvious candidate definitions for tightness of arbitrary contaminations, but in this paper we
restrict our attention to pure contaminations.  

If we want to use tight contaminations to prove results about 3-manifolds, it will be necessary to understand modifications which
preserve the property of tightness.  Our knowledge in this domain is currently limited, but what we know is described in Section
\ref{ContaminationTightnessInvariance}.

In Sections \ref{Tightness} and \ref{ContaminationsTightness}
we prove the ``internal tightness" of some special classes of pure contaminations carried by branched surfaces.   If $\xi$ is a pure
contamination carried by $B$, then 
$\xi$ is {\it internally (universally) tight} if there are no overtwisting discs in  $\intr(V(B))$ ($\intr(\tilde V(B))$), where
$\tilde V(B)$ is the universal cover of $V(B)$.

\begin{thm} \label{TightnessTheoremOne}  Suppose a branched surface $B\embed M$ carries a positive $\sigma$-contamination
$\xi$.  Suppose that $\bdry_hV(B)$ has no sphere components and suppose that the branch locus of 
$B$ consists of a simple system of essential closed curves, i.e. there are no double points of the branch locus, and the smooth closed
curves of the branch locus are $\pi_1$-injective.   Then $\xi$ is internally universally tight.
\end{thm}

A somewhat stronger version of this theorem (Theorem \ref{TightnessTheorem}) is proved in Section \ref{Tightness}, which has the
following corollary.  The definition of a ``immersed twisted disc of contact" is given in Section \ref{Tightness}, but the following
corollary applies, in particular, when there are no double points of the branch locus.
 
\begin{corollary}\label{ContaminationTightnessCorollary}
Suppose that $B\embed M$ is a branched surface satisfying the following conditions:
\item{(1)} the branch locus of $B$ consists of a system of essential closed curves;
\item{(2)} there are no immersed twisted discs of contact in $B$;
\item{(3)} the boundary $\partial V(B)$ is incompressible in $M$;
\item{(4)} $\bdry_hV(B)$ does not contain sphere components or disc components.
\trip

\noindent
Then any pure contamination carried by $B$ is tight.
\end{corollary}

The above corollary can be used to construct examples of tight contaminations, see Section \ref{ContaminationsTightness}.
There are conditions one can impose on a branched surface $B$ sufficient to ensure than any contamination carried by $B$ represents
a tight contact structure.

It remains now to say a few more words about our research program.  One obvious goal is to shed light on the subject of contact
structures in 3-manifolds using our new methods.  Establishing the foundations of our theory was difficult.  If we improve and refine our
theory enough, we may be able to offer a different approach to the subject of contact structures, sufficiently powerful to give
alternative proofs of many known results, and also to uncover new results. Our other goal is to use tight contact structures and tight
pure contaminations to prove results about 3-manifolds.  To do this it would be best to be able to work entirely with branched surfaces,
since branched surfaces are finite combinatorial objects closer to the traditional tools in 3-manifold topology.  For the neatest
possible theory, it would therefore be best to prove theorems of the following form:  (1) If $B\embed M$ has properties $X$, then any
pure contamination
$\xi$ carried by $B$ is tight;  (2) If $\xi$ is a tight pure contamination in $M$, it is carried by a branched surface $B$ with
properties $X$.  Results of this kind were proved for incompressible surfaces and essential laminations, see
\cite{WFUO:IncompressibleViaBranched} and \cite{DGUO:EssentialLaminations}.  For tight pure contaminations this promises to be more
difficult.

We thank Yakov Eliashberg, Ko Honda, and William Kazez for their help, and gratefully acknowledge the support of the
National Science Foundation, the Polish State Committee for Scientific Research (KBN), and the Max Planck Institute, Bonn.

\section{Eliminating maximal integral surfaces} \label{Reduction}

As a first step towards proving Theorem \ref{SigmaDeformationThm}, we will show in this section, given a
$\sigma$-confoliation
$\zeta_0$ on a closed 3-manifold $M$, how to construct a $\sigma$-deformation $\zeta_s$ from $\zeta_0$ to a contact
structure
$\zeta_1$. 

The elimination of the maximal integral surface will be
achieved using a sequence of deformations, each of which either eliminates a disc component from
the maximal integral surface $F$, or removes a 1-handle in $F$.  

We must establish some further notation and definitions before describing precisely the deformations we need.
Let $\xi_t:t\in[a,b]$ be a $\sigma$-deformation, and
$h_t:t\in[a,b]$ a 1-parameter family of diffeomorphisms
of the underlying 3-manifold $M$ as in Definition \ref{SigmaDeformationDef}.
We will say that the family $h_t$ puts
$\xi_t$ into a {\it canonical form} (with respect to maximal
integral surfaces).

Let $I_1,I_2,\dots,I_k$ be the sequence of subintervals
as in Definition \ref{SigmaDeformationDef}, corresponding
to families $\xi_t$ and $h_t$ as above. For $1\le j\le k-1$, denote by $t_j$ the
common endpoint of the subintervals $I_j$ and $I_{j+1}$.
Points $t_1,t_2,\dots,t_{k-1}$ will be called {\it discontinuity
points} of a $\sigma$-deformation $\xi_t$.

A $\sigma$-deformation with a single discontinuity point $t_1$
is a {\it reduction}, if the (strict)
inclusion $\Sigma_2\subset\Sigma_1$ holds; it is an {\it expansion}
in the opposite case. Clearly, each $\sigma$-deformation can be
decomposed into a combination of reductions and expansions.
In this section we describe examples of $\sigma$-deformations
that we call {\it controlled reductions}. We use combinations
of controlled reductions to show that each $\sigma$-confoliation
can be connected by a $\sigma$-deformation to a contact structure.
We postpone the discussion of expansions until later in this section.

Controlled reductions will be defined in terms of certain charts, which we now 
define.

\definition \label{CChartDef}
 Let $(r,\theta,z)$ be the cylindrical coordinates
in $\reals^3$, and let ${\cal C}=\{(r,\theta,z)\in \reals^3:r<1,-1<z<1\}$. A {\it $\cal C$-chart or cylinder chart} for a
confoliation
$\xi$ in a 3-manifold $M$ is a smooth embedding $\psi:{\cal C}\to M$, such that:
\item{(1)} the images in $M$ of all curves $\{r=\text{const},\theta=\text{const}\}$
in $\cal C$ are transverse to $\xi$;
\item{(2)} the images in $M$ of all curves $\{\theta=\text{const},z=\text{const}\}$
in $\cal C$ are everywhere tangent to $\xi$.
\enddefinition

Despite the restrictions imposed by conditions (1) and (2) above,
there is quite some freedom in constructing $\cal C$-charts for a confoliation $\xi$.  We state without proof the
following easy fact.

\begin{fact}\label{CylinderChartFact}
For any
disc $D^2\subset M$ radially foliated by curves tangent to $\xi$,
and for any vector field $X$ transverse to $\xi$ in an open
neighborhood $U$ of $D^2$ and transverse to $D^2$,
there exists a $\cal C$-chart $\psi:{\cal C}\to M$ such that
\item{(1)} $\psi(\{(r,\theta,z)\in{\cal C}:z=0\})=\hbox{int}(D^2)$;
\item{(2)} $\psi({\cal C})\subset U$;
\item{(3)} the image in $M$ of any curve $\{ r=\text{const},\theta=\text{const}\}$ in $\cal C$ is an integral curve of the
vector field $X$.
 \end{fact}

\begin{definition}\label{BChartDef} Let $(x,y,z)$ be the cartesian coordinates
in $\reals^3$, and let ${\cal B}=\{(x,y,z)\in \reals^3:x,y,z\in(-1,1)\}$.
A {\it $\cal B$-chart or box chart} for a confoliation $\xi$ in a 3-manifold $M$
is a smooth embedding $\psi:{\cal B}\to M$, such that
\item{(1)} the images in $M$ of all curves $\{x=\text{const},y=\text{const}\}$
in $\cal B$ are transverse to $\xi$;
\item{(2)} the images in $M$ of all curves $\{x=\text{const},z=\text{const}\}$ in $\cal B$ are everywhere tangent to $\xi$.
\end{definition}

The freedom for constructing $\cal B$-charts for a confoliation $\xi$
is similar to that for $\cal C$-charts, see Fact \ref{CylinderChartFact} . We omit the details.

The advantage of working with $\cal C$-charts and $\cal B$-charts
for a confoliation $\xi$ is that in the coordinates provided
by these charts $\xi$ can be expressed in a unique way
as the kernel of a 1-form $dz+f(r,\theta,z)d\theta$ ($dz-f(x,y,z)dx$ respectively),
for some smooth function $f$. We will call the function $f$ as above a {\it slope function}. Deformations of the slope
function $f$ then yield deformations of the corresponding
confoliation $\xi$ inside a chart neighborhood in $M$. To perform such deformations
we need the following two lemmas.

\begin{lemma}\label{BoxChartLemma}  Let $\omega=dz-f(x,y,z)dx$ be a 1-form
and $\xi=\ker(\omega)$ be the induced plane field in $\cal B$.
Then $\xi$ is a positive confoliation if and only if $\partial f/\partial y(x,y,z)\ge0$
for each $(x,y,z)\in{\cal B}$. Moreover, $\xi$ is contact at a point $(x,y,z)\in{\cal B}$ if and only if $\partial
f/\partial y(x,y,z)>0$.
\end{lemma}

\begin{lemma}\label{CylinderChartLemma}
 Let $\omega=dz+f(r,\theta,z)d\theta$ be a 1-form
and $\xi=\ker(\omega)$ be the induced plane field in $\cal C$.
\item{(a)} The form $\omega$ is well-defined and smooth in $\cal C$ if and only if
$f(r,\theta,z)=r^2\cdot h(r,\theta,z)$ for some smooth function
$h:{\cal C}\to \reals$.
\item{(b)} The plane field $\xi$ is a positive confoliation if and only if
${{\partial f}/{\partial r}}(r,\theta,z)\ge0$ for all points in $\cal C$ with $r>0$.
\item{(c)} The plane field $\xi$ is a positive contact structure
at a point $(r,\theta,z)\in{\cal C}$ with $r>0$ if and only if
 ${{\partial f}/{\partial r}}(r,\theta,z)>0$.
\item{(d)} The plane field $\xi$ is a positive contact structure
at a point $(0,\theta,z)\in{\cal C}$ if and only if
the function $h$
as in (a) satisfies the condition $h(0,\theta,z)>0$.
\end{lemma}

The proof of Lemma \ref{BoxChartLemma} is straightforward, and can be found in [ET].
We will include the slightly more difficult proof of Lemma \ref{CylinderChartLemma}.

\begin{proof} To prove (a), consider coordinates $(x,y,z)$ in $\cal C$
with $x=r\cos\theta$, $y=r\sin\theta$. In these coordinates we have
$$
d\theta={x\,dy-y\,dx\over x^2+y^2} \hbox{ and hence }
\omega=dz+{f(x,y,z)\over x^2+y^2}(x\,dy-y\,dx).
$$
For a function $h=f/r^2$ we then have $\omega=dz+h(x,y,z)(x\,dy-y\,dx)$,
and therefore $\omega$ is smooth if and only if the function $h:{\cal B}\to \reals$
is well-defined and smooth.

 Having proved (a), the other parts of the lemma follow
by direct calculation.
\end{proof}

The following lemma
establishes the existence of a $\sigma$-deformation which eliminates a disc from the maximal integral surface of a
$\sigma$-confoliation.  Any deformation of the kind described in the lemma is called a {\it controlled reduction} of a disc
component.

\begin{lemma}\label{DiscReductionLemma}  Suppose $\xi$ is  a $\sigma$-confoliation with maximal integral surface $F$, and
suppose $D$ is a disc component of $F$.  There is a $\sigma$-deformation of
$\xi$, with one discontinuity point, which reduces the maximal integral surface from $F$ to $F- D$.
The deformation can be chosen to be supported in an arbitrarily small
neighborhood $U$ of the component $D$, i.e. it can be chosen to be constant
outside of $U$.
\end{lemma}

\begin{proof}  
Given a neighborhood $U$ of $D$ and a vector field $X$
in $U$ transverse to $\xi$, whose role is purely technical, it is not difficult to find,
essentially using Fact \ref{CylinderChartFact},
 a $\cal C$-chart $\psi:{\cal C}\to M$ for $\xi$,
with the following properties.

\hop
\noindent{\it Properties C:}
\item{(1)} $\psi({\cal C})$ is contained in $U$ and disjoint from
$F- D$;
\item{(2)} $D$ is contained in $\psi({\cal C})$ and $\psi^{-1}(D)=\{(r,\theta,z)\in{\cal C}:z=0, r\le1-\epsilon\}$, where
$\epsilon$ is a small positive number;
\item{(3)} the images in $M$ of the curves
$\{(r,\theta,z)\in{\cal C}:r=\text{const},\theta=\text{const}\}$ are all contained in the
integral curves of the vector field $X$.
\hop

Let $f:{\cal C}\to \reals$ be the slope function (as defined just before
Lemma \ref{BoxChartLemma}) for the $\sigma$-confoliation $\xi$, in some $\cal C$-chart
coordinates as above. Then there exists a smooth family $f_t:t\in[0,1]$
of functions $f_t:{\cal C}\to \reals$ having the following properties.

\hop
\noindent{\it Properties PC:}

\item{(1)} $f_t=f$ for all $t\le1/2$;
\item{(2)} for each $t\in[0,1]$ we have $f_t=h_t\cdot r^2$,
with $h_t:{\cal C}\to \reals$ smooth and smoothly depending on $t$;
\item{(3)} $h_t(0,\theta,z)>0$ for all points $(0,\theta,z)$ in $\cal C$
 and for all $t>1/2$;
\item{(4)} $\partial f_t/\partial r(r,\theta,z)>0$ for all
$(r,\theta,z)\in{\cal C}$ with $r>0$ and for all $t>1/2$;
\item{(5)} for each $t\in[0,1]$ we have $f_t=f$ in a neighborhood of
the boundary of $\cal C$ in $\reals^3$.
\trip

\noindent
Note that, by Lemma \ref{CylinderChartLemma}, conditions (2), (3) and (4) are already
satisfied by the slope function $f$, except that (3) and (4) do not hold
for $z=0$ and $r\le1-\epsilon$ (i.e. at points of
the preimage $\psi^{-1}(D)$), where the corresponding
inequalities are replaced by equalities.
 The existence of a family $f_t$ as above is then fairly obvious.

By the properties of functions $f_t$, it follows from Lemma \ref{CylinderChartLemma}
that a 1-parameter family $\xi_t:t\in[0,1]$ of plane fields in $M$ defined by

$$
\xi_t=\left\{\begin{array}{cc}\ker(dz+f_t d\theta)& \hbox{in $\cal C$-chart coordinates
inside $\psi({\cal C})$}\\
   \xi& \hbox{ outside $\psi({\cal C})$}\end{array}\right.
$$

\noindent is a $\sigma$-deformation as required, with the discontinuity point
$t_1=1/2$.
\end{proof}

The next lemma, analogous to the previous lemma, 
establishes the existence of a $\sigma$-deformation which eliminates a 1-handle or ``bridge" from the maximal integral surface $F$
 of a
$\sigma$-confoliation.  Thus, a bridge is a rectangular subsurface of $F$ with two opposite sides contained in $\bdry F$
and the remaining sides contained in $\intr (F)$, except their endpoints.  Any deformation of the kind described in the
lemma is called a {\it controlled reduction} of a bridge.

\begin{lemma}\label{BridgeReductionLemma}  Suppose $\xi$ is  a $\sigma$-confoliation with maximal integral surface $F$, and
suppose $B$ is a bridge in $F$.  There is a $\sigma$-deformation of  $\ \xi$, with one discontinuity point, which reduces
the maximal integral surface from $F$ to $\hbox{cl}({F- B})$. The deformation can be chosen to be supported in an
arbitrarily small neighborhood $U$ of the bridge $B$, i.e. it can be chosen to be constant
outside of $U$.
\end{lemma}

\begin{proof}

Since the proof of this lemma is in many respects similar to that
of Lemma \ref{DiscReductionLemma}, our arguments here are more sketchy.

Let $U\subset M$ be a neighborhood of the bridge $B$ and $X$ a
vector field in $U$ transverse to $\xi$, whose role is again purely technical.
Then there exists a $\cal B$-chart $\psi:{\cal B}\to M$
for $\xi$ with the following properties.

\hop
\noindent{\it Properties B:}
\item{(1)} $\psi({\cal B})$ is contained in $U$;
\item{(2)} $B$ is contained in $\psi({\cal B})$ and
$\psi^{-1}(B)=\{(x,y,z)\in{\cal B}:z=0,
|x|\le1-\epsilon,\ \phi_1(x)\le y\le\phi_2(x)\}$,
where $\epsilon>0$ is a small number and
$\phi_1,\phi_2:[-1+\epsilon,\ 1-\epsilon]\to(-1,1)$ are two functions
with $\phi_1<\phi_2$;
 ,\item{(3)} $\psi^{-1}(F)=\psi^{-1}(B)\cup             \{(x,y,z)\in{\cal B}:z=0,\ |x|\ge1-\epsilon\}$;
\item{(4)} the images in $M$ of the curves
$\{(x,y,z)\in{\cal B}:x=\text{const},y=\text{const}\}$ are all contained in the integral curves
of the vector field $X$.
\trip

See Figure \ref{ContactBChart}.

\begin{figure}[ht]
\centering
\scalebox{1.0}{\includegraphics{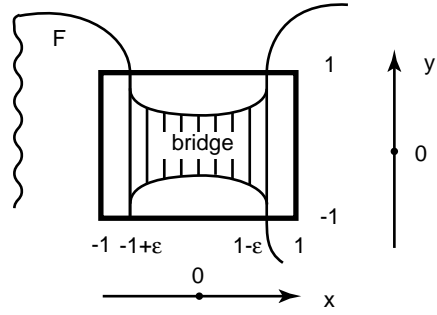}}
\caption{\small $\cal B$-chart for a bridge.}
\label{ContactBChart}
\end{figure}

Let $f:{\cal B}\to \reals$ be the slope function for the $\sigma$-confoliation
$\xi$, in some $\cal B$-chart coordinates as above.
Then there exists a smooth family $f_t:t\in[0,1]$
of functions $f_t:{\cal B}\to \reals$ having the following properties.

\hop
\noindent{\it Properties PB:}
\item{(1)}  $f_t=f$ for all $t\le1/2$;\item{(2)} if $t>1/2$, then
$f_t(x,y,z)=0$ for all points $(x,y,z)\in{\cal B}$ with
$z=0$ and $|x|\ge1-\epsilon$ and
$\partial f_t/\partial y(x,y,z)>0$ for all other points in $\cal B$;
\item{(3)} for each $t\in[0,1]$ we have $f_t=f$ in a neighborhood
of the boundary of $\cal B$ in $\reals^3$.
\trip

In view of above properties of functions $f_t$, it follows from
Lemma \ref{BoxChartLemma} that a 1-parameter family $\xi_t:t\in[0,1]$
 of plane fields in $M$ defined by

$$
\xi_t=\left\{\begin{array}{cc}\ker(dz-f_t dx)& \hbox{in $\cal B$-chart coordinates inside $\psi({\cal B})$}\\
   \xi& \hbox{ outside $\psi({\cal B})$}\end{array}\right.
$$

\noindent is a $\sigma$-deformation as required, with the discontinuity
point $t_1=1/2$.
\end{proof}

Combining simple $\sigma$-deformations, it is possible to prove the following: 

\begin{proposition} \label{IntegralSurfaceEliminationProp} For any $\sigma$-confoliation $\xi$ there is a $\sigma$-deformation
$\xi_t$, $0\le t\le 1$, a composition of finitely many controlled reductions, changing
$\xi_0=\xi$ to a contact structure $\xi_1$.
\end{proposition}

\section{Proof of $\sigma$-Deformation Theorem}\label{DeformationTheorem}

We shall prove Theorem
\ref{SigmaDeformationThm} by constructing a 2-parameter family of $\sigma$-confolia\-tions.  Let
$\xi_t:t\in[0,1]$ be a
$\sigma$-deformation between contact structures $\xi_0$ and $\xi_1$, as in the theorem.
We will use ideas developed in Section \ref{Reduction} to construct a smooth
2-parameter family $\xi_{s,t}:s\in[0,1],t\in[0,1]$ of confoliations with the following properties:

\item{(1)} $\xi_{0,t}=\xi_t$ for each $t\in[0,1]$;
\item{(2)} $\xi_{s,0}$ and $\xi_{s,1}$ are contact structures
for all $s\in[0,1]$;
\item{(3)} $\xi_{1,t}$ is a contact structure for each $t\in[0,1]$.
\trip

The following proposition yields a 2-parameter family in a slightly more general setting.

\begin{proposition}\label{TwoParameterProp}
 Let $\xi_t:t\in[0,1]$ be a $\sigma$-deformation.
Then there exists a smooth 2-parameter family $\xi_{s,t}$ of
positive confoliations satisfying the following conditions:
\item{(1)} $\xi_{0,t}=\xi_t$ for each $t\in[0,1]$;
\item{(2)} both $\xi_{s,0}:s\in[0,1]$ and $\xi_{s,1}:s\in[0,1]$ are combinations of controlled reductions, deforming $\xi_0$
and $\xi_1$ respectively into contact structures;
\item{(3)} $\xi_{1,t}$ is a contact structure for each $t\in[0,1]$.
\end{proposition}

\begin{proof} [Proof of Theorem \ref{SigmaDeformationThm} using Proposition \ref{TwoParameterProp}]
When $\xi_t$ is a $\sigma$-deformation from a contact structure to another contact structure, the 2-parameter
smooth family 
$\xi_{s,t}$ provides a smooth 1-parameter family of contact structures from
$\xi_0=\xi_{00}$ to $\xi_1=\xi_{10}$ by combining the following 1-parameter families in order, following three
consecutive sides of the parameter space $I\times I$ for $(s,t)$ :
$\xi_{s0},\ 0\le s\le 1$ followed by
$\xi_{1t},\ 0\le t\le 1$ followed in turn by $\xi_{(1-s)1},\ 0\le s\le 1$.  Now using Gray's Theorem we can replace the smooth
1-parameter family by an isotopy, i.e., a 1-parameter family $\phi_t:M\to M$ of diffeomorphisms such that for all $t$,
$\phi_{t*}(\xi_0)=\xi_t$.
\end{proof}

From the theorem and proposition, it is now clear that for every $\sigma$-confolia\-tion $\xi_0$, if one gets 
 $\xi_1$ from $\xi_0$ via {\it any} $\sigma$-deformation $\xi_t,\ 0\le t\le 1$, then one always obtains the same contact structure up to
isotopy, and we call this the {\it contact structure associated to $\xi_0$}.  Furthermore, this contact structure is actually associated
to an equivalence class of $\sigma$-confoliations, the equivalence relation being $\sigma$-deformation.  The existence of the
$\sigma$-deformation $\xi_t$ from $\xi_0$ to some contact structure $\xi_1$ is guaranteed by Proposition  
\ref{IntegralSurfaceEliminationProp}.

We will prove the proposition by proving a sequence of lemmas, some of which are special cases of the proposition.
The first lemma will deal with deformations $\xi_t$ having
no discontinuity point. The proof in this case will consist of
constructing a 2-parameter family of
combinations of controlled reductions. Two more lemmas
prove the proposition for deformations $\xi_t$ with one discontinuity point.
This will require a more skilled use of special charts,
to construct smooth parametric families of controlled reductions
extended by deformations with no discontinuity points.
The final lemma is used as an inductive step for combining 2-parameter families.  This will require a
delicate argument for showing how to glue to each other 2-parameter families of confoliations corresponding to
deformations with single discontinuity points.

\begin{lemma}\label{NoDiscontinuityLemma} Proposition \ref{TwoParameterProp} is true if the $\sigma$-deformation $\xi_t$
has no discontinuities.
\end{lemma}
\begin{proof}
By putting such a deformation into a canonical form (as explained
at the beginning of Section \ref{Reduction}), we will assume
that all the confoliations $\xi_t $ have the same maximal integral surface, which we denote by $\Sigma$.

Let $D$ be a connected component of $\Sigma$ which is a disc,
and $X$ be a vector field in $M$ transverse to $D$. Choose
a neighborhood $U$ of $D$ such that for each $t\in[0,1]$
the field $X$ is transverse to $\xi_t$ in $U$.
Then there exists a smooth family $\psi_t:t\in[0,1]$ of maps
$\psi_t:{\cal C}\to M$ such that for each $t\in[0,1]$

\hop
\item{(1)} $\psi_t$ is a $\cal C$-chart for $\xi_t$;
\item{(2)} $\psi_t$ has Properties C, as required for a $\cal C$-chart $\psi$ in the proof of Lemma
\ref{DiscReductionLemma}.

\trip

For each $t\in[0,1]$,
let $f^t:{\cal C}\to \reals$ be the slope function for $\xi_t$ with respect
to a $\cal C$-chart $\psi_t$ as above.
Clearly, the family $f^t:t\in[0,1]$ of functions depends smoothly
on $t$.
Consider then a smooth 2-parameter family $f^t_s:s\in[0,1],t\in[0,1]$
of functions $f^t_s:{\cal C}\to \reals$ with the following properties:

\hop
 \item{(1)} $f^t_0=f_t$ for each $t\in[0,1]$;
\item{(2)} for any fixed $t\in[0,1]$ the family $f^t_s :s\in[0,1] $ has Properties PC, as required in the proof of Lemma
\ref{DiscReductionLemma} (where the role of parameter $t$ in
$f_t$ is taken by $s$ in the family $f^t_s$).

\trip

\noindent The existence of a 2-parameter family $f^t_s$ as above is as obvious as
the existence of a family $f_t$ in the proof of Lemma \ref{DiscReductionLemma}.

Consider a 2-parameter family $\xi^t_s$ of $\sigma$-confoliations
defined by

$$
\xi^t_s=\left\{\begin{array}{cc}\ker(dz+f^t_s d\theta)& \hbox{in cylindrical coordinates induced by $\psi_t$ inside
$\psi_t({\cal C})$}\\
   \xi& \hbox{ outside $\psi_t({\cal C})$}\end{array}\right.
$$

\noindent We shall call any family $\xi^t_s$ obtained by the construction as above
a {\it parametric controlled reduction} of a disc component.

Clearly one can, in a completely analogous way, define the notion of
a parametric controlled reduction of a bridge in the maximal integral surface $\Sigma$. The proof is then
completed by combining consecutively in $s$-time a sequence of parametric controlled reductions of discs and/or bridges to
obtain the  family $\xi_{s,t}$.
\end{proof}

The following lemma will be used to show that Proposition \ref{TwoParameterProp} is true for a $\sigma$-deformation
$\xi_t$ with one discontinuity.

\begin{lemma} \label{OneDiscontinuityLemma} Let $\xi_t:t\in[0,1]$ be a $\sigma$-deformation with
exactly one discontinuity point.  Suppose the maximal
integral surface $F_t$ of $\xi_t$ reduces from $\Sigma_1$ to $\Sigma_2$, i.e.
$\Sigma_2\subset\Sigma_1$ and
$F_t=\Sigma_1$ for $t\le t_1$ while $F_t=\Sigma_2$ for $t> t_1$(for some $t_1\in(0,1)$). Then there exists
a smooth 2-parameter family $\xi_{s,t}$ of $\sigma$-confoliations such that
\item{(1)} $\xi_{0,t}=\xi_t$ for each $t\in[0,1]$;
\item{(2)} $\Sigma_2$  is the maximal integral surface for all confoliations $\xi_{1,t}$ with $t\in[0,1]$ and $\xi_{s,1}$
with $s\in[0,1]$;
\item{(3)} $\xi_{s,0}:s\in[0,1]$ is a sequence of controlled reductions.
\end{lemma}

\begin{proof}
There is a sequence of operations of deleting a bridge
or a disc component, which converts $\Sigma_1$ into $\Sigma_2$.  To prove the lemma, it is then sufficient to consider the case when
$\Sigma_2$ is obtained from $\Sigma_1$ using sequence consisting of a single operation, deleting a bridge or disc component, since the
general case will then follow by induction. Therefore, we assume that, to get $\Sigma_2$, it is sufficient to delete from $\Sigma_1$ a
disc component
$D$. The reasoning in the case when $\Sigma_2$ is obtained from $\Sigma_1$ by deleting a bridge
is analogous, and we omit it.

Since $D$ is a connected component of the maximal integral surfaces in $\xi_t$
for all $t\le t_1$, we may consider a smooth family
$\psi_t:t\in[0,t_1]$ of maps $\psi_t:{\cal C}\to M$,
with the same properties as a family $\psi_t$ of maps
considered in the proof of Lemma \ref{NoDiscontinuityLemma}.
This family may then be extended smoothly to a family $\psi_t:t\in[0,t_1+\epsilon]$ of maps, for some small positive
$\epsilon$, which all satisfy the following two properties:
\hop
\item{(1)} $\psi_t$ is a $\cal C$-chart for $\xi_t$;
\item{(2)} $\psi_t({\cal C})$ is disjoint from the surface $\Sigma_2=
\Sigma_1- D$ for each $t\in[0,t_1+\epsilon]$.
\trip
\noindent
Since for $t>t_1$ the disc $D$ is no longer an integral surface in $\xi_t$,
we do not require that, for $t>t_1$, $D$ be the core of the each $C$-charts $\psi_t$. In fact, for these
parameters we do not even require that $D$ be contained in $\psi_t({\cal C})$.

For each $t\in[0,t_1+\epsilon]$, let $f_t:{\cal C}\to \reals$ be the slope
function for $\xi_t$ with respect to a $\cal C$-chart $\psi_t$ as above.
Consider then a smooth 2-parameter family $f_{s,t}:t\in[0,t_1+\epsilon],
 s\in[0,1]$ of functions $f_{s,t}:{\cal C}\to \reals$ with the following properties:
\hop
\item{(1)} $f_{0,t}=f_t$ for each $t\in[0,t_1 +\epsilon]$;\item{(2)} for any fixed $t\in[0,t_1+\epsilon]$ the family
$f_{s,t}:s\in[0,1]$ has Properties PC,
as required for a family $f_t$ of functions in the proof of Lemma \ref{DiscReductionLemma}
(where the role of parameter $t$ in $f_t$ is taken by $s$ in the family
$f_{s,t}$);
\item{(3)} for each $t\in[t_1+\epsilon/2,t_1+\epsilon]$ and for each
$s\in[0,1]$ we have $f_{s,t}=f_t$.
\trip

\noindent
As before,
we regard the existence of
such a family $f_{s,t}$ as fairly obvious.
In view of Lemma \ref{CylinderChartLemma}, it follows from the above properties
(1)-(3) of a family $f_{s,t}$ that the 2-parameter family
$\xi_{s,t}$ of plane fields in $M$ defined by

$$
\xi_{s,t}=\left\{\begin{array}{cc}\ker(dz+f_{s,t}d\theta)& \hbox{in cylindrical coordinates induced
           }\\
&\hbox{by $\psi_t$ inside $\psi_t({\cal C})$, if $t\le t_1+\epsilon$    }\\
   \xi_t& \hbox{otherwise,}\end{array}\right.
$$

\noindent is well-defined, consists of $\sigma$-confoliations, and satisfies
the properties (1)-(3) asserted in the lemma.
\end{proof}

\begin{corollary}  Proposition \ref{TwoParameterProp} is true for a $\sigma$-deformation $\xi_t$ with exactly one
discontinuity.
\end{corollary}

\begin{proof}
By putting the deformation $\xi_t$ into a canonical form,
and by reversing parameter $t$ if necessary, we assume that the maximal
integral surface $F_t$ of $\xi_t$ reduces from $\Sigma_1$ to $\Sigma_2$, i.e.
$\Sigma_2\subset\Sigma_1$ and
$F_t=\Sigma_1$ for $t\le t_1$ while $F_t=\Sigma_2$ for $t> t_1$(for some $t_1\in(0,1)$).
We apply the previous lemma to the deformation $\xi_t$ which yields a 2-parameter family $\xi_{s,t}$ with the property
that $\xi_{1,t}$ has no discontinuities. We follow the 2-parameter family $\xi_{s,t}$ in
$s$-time by a 2-parameter family as in Lemma \ref{NoDiscontinuityLemma} to eliminate $\Sigma_2$.
\end{proof}

Now we deal with the case of two discontinuities, which contains all the ideas to complete the proof of Proposition
\ref{TwoParameterProp}.

\begin{lemma}
Proposition 3.1 is true for a $\sigma$-deformation $\xi_t$
with exactly two discontinuities.
\end{lemma}

\begin{proof}  Suppose the discontinuities $t_1,t_2$ are contained in the intervals
$(0,1/2)$ and $(1/2,1)$ respectively, which we can always assume by
reparametrizing in $t$ if necessary. Let $\xi_t^1:t\in[0,1/2]$ be defined
by $\xi_t^1:=\xi_t$, and similarly $\xi_t^2:t\in[1/2,1]$ by $\xi_t^2:=\xi_t$.
Then $\xi_t^1$ and $\xi_t^2$ are $\sigma$-deformations with one discontinuity
each. Let $\xi_{t,s}^1:t\in[0,1/2],s\in[0,1]$ and 
$\xi_{t,s}^2:t\in[1/2,1],s\in[0,1]$ be smooth 2-parameter families
constructed as in Lemma 3.3. 

The obvious problem with simply glueing
families $\xi^1_{s,t}$ and $\xi^2_{s,t}$ into a single family $\xi_{s,t}$
is that deformations $\xi^1_{s,1/2}:s\in[0,1]$ and $\xi^2_{s,1/2}: s\in[0,1]$ need not to coincide. Both of these
deformations are sequences of controlled reductions, but the sequences
of bridges and disc components deleted from the common initial
maximal integral surface may be very different in the two deformations. To overcome this
difficulty, we need to prove the following claim, which provides
a link between deformations $\xi^1_{s,1/2}$ and $\xi^2_{s,1/2}$.

\hop\noindent
{\bf Claim.} {\it If all $\sigma$-confoliations in a family
$\xi^2_{s,1/2}:s\in[0,1]$ are sufficiently  close to the $\sigma$-confoliation $\xi^2_{0,1/2}$, then there exists a
smooth 2-parameter family
$\eta_{u,v}:u\in[0,1],v\in[0,1]$ of positive confoliations, such that:
 \item{(1)} $\eta_{u,0}=\xi^1_{u,1/2}$ for each $u\in[0,1]$;
\item{(2)} $\eta_{0,v}=\xi^2_{v,1/2}$ for each $v\in[0,1]$;
\item{(3)} all confoliations $\eta_{u,1}:u\in[0,1]$ and
$\eta_{1,v}:v\in[0,1]$ are contact structures.}

\trip
Before proving the claim, we will finish the proof of the lemma,
assuming that the claim is true. To do this, note that a 2-parameter family $\xi_{s,t}$ obtained from Lemma
\ref{OneDiscontinuityLemma} can be always constructed
in such a way that the restricted deformation $\xi_{s,0}:s\in[0,1]$
consists of confoliations as close to $\xi_{0,0}$ as necessary.
This is due to the fact that a controlled reduction
of a given bridge or of a given disc component can be made to consist
of confoliations arbitrarily close to the initial one. In particular,
it is always possible to choose a family $\xi^2_{s,t}:s\in[0,1],\ t\in [1/2,1]$ so that the restricted deformation
$\xi^2_{s,1/2}:s\in[0,1]$ satisfies the assumption of the Claim. The lemma then follows by glueing
families $\xi^1_{s,t}$, $\xi^2_{s,t}$ and $\eta_{u,v}$,
and performing some cosmetic adjustment in order to make the
resulting family $\xi_{s,t}$ smooth, see Figure 3. We omit further details.

\begin{figure}[ht]
\centering
\scalebox{1.0}{\includegraphics{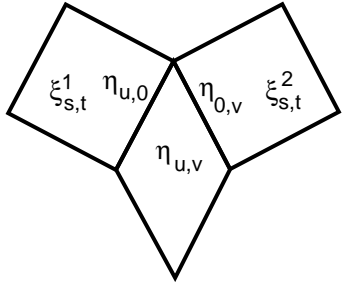}}
\caption{\small Identifications for domain of 2-parameter family.}
\label{ContactTwoParameter}
\end{figure}

 \vskip5pt\noindent
{\it Proof of Claim:} Put both $\sigma$-deformations
$\xi^1_{u,1/2}:u\in[0,1]$ and $\xi^2_{v,1/2}:v\in[0,1]$
into canonical forms, with sequences
$\Sigma_1,\Sigma_2 ,\dots,\Sigma_k=\emptyset$ and $\Sigma'_1,\Sigma'_2,\dots,\Sigma'_l=\emptyset$ of maximal integral
surfaces respectively, so that still $\xi^1_{0,1/2}=\xi^2_{0,1/2}$
and hence also $\Sigma_1=\Sigma'_1$. Moreover, let $t_1,t_2,\dots,t_{k-1}$
be the sequence of discontinuity points for $\xi^1_{u,1/2}:u\in[0,1]$.
The proof will essentially go by induction with respect to $k$.

As a first step in this induction we shall construct a smooth
2-parameter family $\eta^1_{u,v}:u\in[0,t_1+\delta],v\in[0,1]$
of confoliations, where $\delta>0$ is such that $t_1+\delta<t_2$, satisfying conditions:

\hop
\item{(1)} $\eta^1_{u,0}=\xi^1_{u,1/2}$ for each $u\in[0,t_1+\delta]$;
\item{(2)} $\eta^1_{0,v}=\xi^2_{v,1/2}$ for each $v\in[0,1]$;
\item{(3)} all confoliations $\eta^1_{u,1}:u\in[0,t_1+\delta]$
are contact structures;
\item{(4)} all confoliations $\eta^1_{t_1+\delta,v}:v\in[0,1]$
are contact at all points of the set $M- \Sigma_2$.
\trip

 Since, by the assumption, $\xi^1_{u,1/2}:u\in[0,1]$ is a combination of controlled reductions, it follows that $\Sigma_2$ is
obtained from $\Sigma_1$ by deleting a bridge or a disc component. Assume the first of those two cases, so that
$\Sigma_2=cl(\Sigma_1- B)$, where $B$ is a bridge in $\Sigma_1$.
The argument for the second case is analogous, and we will omit it.

Let $\psi:{\cal B}\to M$ be the $\cal B$-chart used for determining
the controlled reduction $\xi^1_{u,1/2}:u\in[0,t_1+\delta]$, as in
Lemma \ref{BridgeReductionLemma}, and let $f_u:u\in[0,t_1+\delta]$ be the corresponding family of
slope
 functions.  Further, let $X$ be the vector field in a neighborhood of the bridge $B$, corresponding to $\psi$ in the
same lemma, Lemma \ref{BridgeReductionLemma}. If confoliations $\xi^2_{v,1/2}:v\in[0,1]$ are close enough to the
confoliation
$\xi^2_{0,1/2}$, then
$X$ is transverse to all of them (in the same neighborhood of $B$), and there exists a smooth
1-parameter family $\psi_v:v\in[0,1]$ of maps $\psi_v:{\cal B}\to M$,
with $\psi_0=\psi$, such that for each $v\in[0,1]$:

\hop
\item{(1)} $B$ is contained in $\psi_v({\cal B})$ and
$$
\psi_v^{-1}(B)=\{(x,y,z)\in{\cal B}:z=g_v(x,y), |x|\le1-\epsilon,
\phi_1(x)\le y\le\phi_2(x)\},
$$
\noindent where $\epsilon$, $\phi_1$ and $\phi_2$ are the same as for $\psi$,
and $g_v:v\in[0,1]$ is a family of smooth functions $g_v:(-1,1)^2\to(-1,1)$,
smoothly depending on $v$;
\item{(2)}
$\psi_v^{-1}(\Sigma_2)=\{(x,y,z)\in{\cal B}:z=g_v(x,y),|x|\ge1-\epsilon\}$,
where $\epsilon$ and $g_v$ are the same as in (2);
\item{(3)} the  images by $\psi_v$ of the curves$\{(x,y,z)\in{\cal B}:x=\text{const},y=\text{const}\}$ are all contained
in the integral curves of the vector field $X$.
\trip

The main difficulty in constructing a family $\psi_v:v\in[0,1]$
as above, is to make all the $\cal B$-charts in the family
contain $B$ in the image (condition (2)).
The difficulty comes from the fact that,
in general, $B$ is not integral in all
of the confoliations $\xi^2_{v,1/2}:v\in[0,1]$. This is where we must make use of the assumption of the claim, which
clearly allows the construction of a desired family.
 
The construction of a family $\eta^1_{u,v}$ of confoliations is now based on
the same idea used in previous lemmas dealing with special cases of Proposition \ref{TwoParameterProp}. For each
$v\in[0,1]$, let $f^v:{\cal B}\to \reals$ be the slope function for the confoliation
$\xi^2_{v,1/2}$, in the $\cal B$-chart $\psi_v$. Then there exists
a smooth 2-parameter family $f_u^v:u\in[0,t_1+\delta],v\in[0,1]$
of functions $f_u^v:{\cal B}\to \reals$ such that:

\hop
\item{(1)} $f_u^0=f_u$ for each $u\in[0,t_1+\delta]$, where $f_u$
is  the family of slope functions for the controlled reduction $\xi^1_{u,1/2}:u\in[0,t_1+\delta]$, as mentioned above;
\item{(2)} $f_0^v=f^v$ for each $v\in[0,1]$;
\item{(3)} $\partial f_u^v/\partial y\ge0$ for each $u\in[0,t_1+\delta]$
and for each $v\in[0,1]$;
\item{(4)} $\partial f_{t_1+\delta}^v/\partial y>0$ outside the set
$\psi_v^{-1}(\Sigma_2)$;
\item{(5)} $f_u^v=f^v$ at a neighborhood of the boundary of
$\cal B$ in $\reals^3$, for all $u\in[0,t_1+\delta]$ and all $v\in[0,1]$.
\trip

\noindent
The existence of such a family $f_u^v$ of functions is fairly clear.

In view of the above properties of functions $f_u^v$,
it follows from Lemma \ref{BoxChartLemma} that a 2-parameter family
$\eta^1_{u,v}:u\in[0,t_1+\delta],v\in[0,1]$ of confoliations, defined by
$$
\eta^1_{u,v}=\left\{\begin{array}{cc}\ker(dz-f_u^v\cdot dx)& \hbox{in cylindrical coordinates induced}\\
                      &\hbox{ by $\psi_v$ inside $\psi_v({\cal B})\ \ \ \ \ \ \ \ \ \ \ \ \ \ \ \ \ \ \  $ }\\
   \xi^2_{v,1/2}& \hbox{outside $\psi_v({\cal B})$}\end{array}\right.
$$
has the required properties.

To finish the proof of the Claim, note that it is always possible first to choose a deformation
$\xi^2_{v,1/2}:v\in[0,1]$, and then a related family $\eta^1_{u,v}$ as above, so that the confoliations
$\eta^1_{t_1+\delta,v}:v\in[0,1]$ are all arbitrarily close to the confoliation
$\eta^1_{t_1+\delta,0}=\xi^1_{t_1+\delta,1/2}$.
This observation allows an inductive construction of a family
$\eta_{u,v}$ required in the claim, and so the proof is finished.
\end{proof}

\begin{proof} [Proof of Proposition \ref{TwoParameterProp}] The proposition can be proved by induction as follows.
Suppose it is true for a deformation with $k$ discontinuities.  Replacing $\xi^1_{t,s}$ in the previous lemma
by the $k$-discontinuity deformation, the lemma remains true. Therefore, the proposition is true
for a deformation with $k+1$ discontinuities.
\end{proof}

\section{Invariance of tightness}\label{Invariance}

This section is devoted to the proof of Theorem \ref{TightInvarianceTheorem}, showing that tightness of a
$\sigma$-confoliation is invariant under $\sigma$-deformation.

In view of Theorem \ref{SigmaDeformationThm} and Proposition \ref{IntegralSurfaceEliminationProp},
Theorem \label{ConfoliationTightInvarianceTheorem} reduces to the following.

\begin{thm} \label{TightReductionTheorem} Let $\xi'$ be a contact structure obtained
from a confoliation $\xi$ by a sequence of controlled reductions.
Then $\xi'$ is tight iff $\xi$ is tight.
\end{thm}

Indeed, if $\xi_t:t\in[0,1]$ is a $\sigma$-deformation, and $\xi_0'$
and $\xi_1'$ are contact structures obtained from $\xi_0$ and $\xi_1$
respectively, by means of sequences of controlled reductions
(which exist due to Proposition \ref{IntegralSurfaceEliminationProp}), then by Theorem
\ref{SigmaDeformationThm} those contact structures are isotopic. Hence, Theorem \ref{TightReductionTheorem} and the
invariance of tightness under isotopy of contact structures imply that $\xi_0$
is tight iff $\xi_1$ is tight, which is exactly the assertion
of Theorem \ref{TightInvarianceTheorem}. We therefore now devote the remainder of the section to proving
Theorem \ref{TightReductionTheorem}. To do this, we need several preparatory results.

\begin{lemma} \label{DiskReductionDiffeoLemma} Let $D$ be a disc which is a connected component 
of the maximal integral surface for a $\sigma$-confoliation 
$\xi$, and consider a controlled reduction of $D$
which turns $\xi$ into a new $\sigma$-confoliation $\xi'$.
Denote by $p$ the point of the underlying 3-manifold $M$ with
cylindrical coordinates $(0,0,0)$ in the $\cal C$-chart $\psi$
used to describe the above controlled reduction. Then the restricted
$\sigma$-confoliations $\xi|_{M- D}$ and $\xi'|_{M-\{p\}}$
are isotopic by a diffeomorphism $g:M-\{p\}\to M- D$ 
with the following properties:
\item{(1)} $g$ coincides with the identity outside $\psi(\cal C)$;
\item{(2)} each curve in $M$ corresponding by $\psi$
to a curve $\{z=\text{const},\theta=\text{const}\}$ in $\cal C$ is mapped
by $g$ to itself.
\end{lemma}

\begin{proof} Recall from Section \ref{Reduction} that in the $\cal C$-chart $\psi$
used for the description of the controlled reduction of $D$ we have
$$
\xi=\ker(dz+f\cdot d\theta) \hbox{ \quad and \quad} \xi'=\ker(dz+f'\cdot d\theta),
$$
where the functions $f,f':{\cal C}\to \reals$ satisfy the following conditions:

\item{(1)} $f(r,\theta,z)=r^2\cdot h(r,\theta,z)$ and 
$f'(r,\theta,z)=r^2\cdot h'(r,\theta,z)$ for some smooth functions
$h,h':{\cal C}\to \reals$;

\item{(2)} $h'>0$ everywhere in $\cal C$, and $h>0$ except at the set
$\{z=0,r\le1-\epsilon\}=\psi^{-1}(D)\subset{\cal C}$ at which $h=0$;

\item{(3)} $\partial f'/\partial r(r,\theta,z)>0$ for all points
$(r,\theta,z)$ in $\cal C$ with $r>0$;

\item{(4)} $\partial f/\partial r(r,\theta,z)>0$ for all points
$(r,\theta,z)$ in ${\cal C}-\psi^{-1}(D)$ with $r>0$;

\item{(5)} $f=f'$ at a neighborhood of the boundary of $\cal C$
in $\reals^3$.
\trip

Define a map $u:{\cal C}-\{z=0,r\le1-\epsilon\}\to{\cal C}-
\{(0,0,0)\}$ of the form $u(r,\theta,z)=(w(r,\theta,z),\theta,z)$, by putting
$$
w(r,\theta,z)=(f_{\theta,z}')^{-1}\circ f_{\theta,z}(r),
$$
where $f_{\theta,z}(r)=f(r,\theta,z)$ and $f_{\theta,z}'(r)=f'(r,\theta,z)$.
Due to the properties (1)-(5) above $u$ is a well-defined homeomorphism
equal to the identity at a neighborhood of the boundary of ${\cal C}$
in $\reals^3$. Moreover, $u$ is easily seen to be regular at all points $(r,\theta,z)$ 
in its domain with $r>0$. Therefore, to prove that $u$ is a diffeomorphism,
it is sufficient to show that it is regular at points with $r=0$. This
follows however from the fact that $u$ can be expressed as a composition
of the map $(r,\theta,z)\to (r\cdot h^{1/2},\theta,z)$ with the inverse
of the map $(r,\theta,z)\to(r\cdot(h')^{1/2},\theta,z)$, and both those
maps are regular at points with $r=0$ and $z\ne 0$ due to property (2).

A straightforward argument shows that, by the definition, the diffeomorphism
$u$ maps the structure $\xi$ to the structure $\xi'$. Since moreover,
by property (5), $u$ is equal to the identity close to the boundary
of $\cal C$, it can be extended to a map $g$ as required, by putting
$g=id$ outside the chart neighborhood $\psi(\cal C)$. This finishes
the proof. 
\end{proof}

\begin{lemma} \label{BridgeReductionDiffeoLemma} Let $\xi_t:t\in[0,1]$ be a controlled reduction of a
bridge $B$, as in Section \ref{Reduction}. Denote by $F_0$ and $F_1$
the maximal integral surfaces for $\xi_0$ and $\xi_1$ respectively,
and by $\gamma$ an arc in $F_0$ defined by
$\gamma=\psi(\{(x,y,0)\in{\cal B}:|x|<1-\epsilon, f_1(x,y,0)=0\})$,
where $\psi:{\cal B}\to M$ is the $\cal B$-chart used to
describe the controlled reduction $\xi_t$, and the rest of
notation is as in Section \ref{Reduction}.
Then the contact structures $\xi_0|_{M- F_0}$ and
$\xi_1|_{(M- (F_1\cup\gamma)}$ are isotopic
by a contactomorphism which coincides with the identity outside
the set $U$ supporting the deformation $\xi_t$.
\end{lemma}

The argument for proving this lemma is analogous
as that for Lemma \ref{DiskReductionDiffeoLemma}, and we omit it.

\hop

Let $F$ be a compact surface, each component of which has nonempty boundary.
A {\it core} in $F$ is a graph $\Gamma$ smoothly embedded in the interior
of $F$, whose image is a deformation retract of $F$. In what follows
we shall identify $\Gamma$ with its image in $F$.

Crucial for proving Theorem \ref{TightReductionTheorem} is the following

\begin{proposition} \label{ReductionDiffeoProposition} (1) Let $\xi$ be a $\sigma$-confoliation in $M$, $F$ its
maximal integral surface, and $\Gamma$ a core in $F$.
Then, given a neighborhood $U$ of $F$, there exists 
a sequence of controlled reductions
which changes $\xi$ into a contact structure $\xi'$ in such a way
that the restricted contact structures $\xi|_{M- F}$
and $\xi'|_{M-\Gamma}$ are isomorphic by a contactomorphism $h$
equal to the identity outside $U$.

(2) A contactomorphism $h$ as above can be chosen to satisfy the following
additional property. For arbitrarily small regular neighborhood $W$ of the
core $\Gamma$ in $M$ the restriction of $h$ to the preimage $h^{-1}(M- W)$
extends to a diffeomorphism of $M$. Regularity of a neighborhood $W$ means
that it is the interior of a compact submanifold with boundary in $M$.

\end{proposition}

\begin{proof} (1) Note first that, by the definition of a core, the pair
$(F,\Gamma)$ is diffeomorphic to the pair $(N\Gamma,\Gamma)$,
where $N\Gamma$ is a small tubular neighborhood in $F$ of $\Gamma$.
It follows that there exists a 
decomposition of the surface $F$ into a union of a disjoint family
$D_v:v\in V\Gamma$ of discs corresponding to the vertices of $\Gamma$,
and a disjoint family $B_e:e\in E\Gamma$ of bridges corresponding
to the edges of $\Gamma$.
Such a decomposition may be easily chosen to have an additional property
that for each edge $e$ of $\Gamma$ the intersection $B_e\cap\Gamma$ 
is a closed subarc contained in the interior of $e$, and it meets the
boundary $\partial B_e$ transversely.

We will construct a required system of controlled reductions
using a system of bridges and discs as above. Two steps of this construction
will consist of first reducing all the bridges and then all
the resulting discs. We describe these two steps separately.

\vskip5pt\noindent
{\bf Step 1. Reduction of the bridges $B_e$.}

Choose a system of $\cal B$-charts $\psi_e:{\cal B}\to M$
around the bridges $B_e$ which satisfy Properties B (as in the proof \ref{BridgeReductionLemma}) and:

\item{(1)} the chart neighborhoods $\psi_e(\cal B)$ in $M$ are
pairwise disjoint;
\item{(2)} for each edge $e$ we have $\psi_e^{-1}(\Gamma)=
\{(x,y,z)\in{\cal B}: y=0,z=0\}$.
\trip

\noindent
For each edge $e$ of $\Gamma$ consider a controlled reduction
of the bridge $B_e$, determined as in Section \ref{Reduction} with respect
to the $\cal B$-chart $\psi_e$. Denote by $\xi_t:t\in[0,1/2]$
a composition of the controlled reductions as above.  We may assume for
each $e$ the following additional property: if $\xi_{1/2}=\ker(dz-f^edx)$
in the chart $\psi_e$, then $f^e(x,0,0)=0$ for each $x\in(-1,1)$.
It follows then from Lemma \ref{BridgeReductionDiffeoLemma} applied inductively that
the contact structures $\xi|_{M- F}$ and 
$\xi_{1/2}|_{M-(\cup D_v\cup\Gamma)}$ are isotopic as required.

\vskip5pt\noindent
{\bf Step 2. Reduction of the discs $D_v$.}

Note that the maximal integral surface of the $\sigma$-confoliation
$\xi_{1/2}$, as obtained in Step 1, is a union of the family
$D_v:v\in V\Gamma$ of pairwise disjoint discs. Choose a system
of $\cal C$-charts $\psi_v:{\cal C}\to M$ around those discs
which satisfy Properties C (as in the proof of Lemma \ref{DiscReductionLemma}) and:

\item{(1)} the chart neighborhoods $\psi_v(\cal C)$ in $M$ are
pairwise disjoint;
\item{(2)} for each vertex $v$ there exists a finite set $\Theta_v$
such that $\psi_v^{-1}(\Gamma)=\{(r,\theta,z)\in{\cal C}:z=0,
\theta\in\Theta_v\}$.
\trip

\noindent
For each vertex $v$ of $\Gamma$ consider a controlled reduction
of the disc $D_v$ determined as in Section \ref{Reduction} with respect
to the $\cal C$-chart $\psi_v$. Denote by $\xi_t:t\in[1/2,1]$
a composition of the controlled reductions above.
By Lemma \ref{DiskReductionDiffeoLemma}, for each $v$ there is a diffeomorphism 
$g_v:\psi_v({\cal C})- D_v\to\psi_v({\cal C})-\psi_v(0,0,0)$
which maps the restricted structure $\xi_{1/2}$ to the restricted
structure $\xi_1$. Since, by the same lemma, the map $g_v$ maps each curve
$\psi_v(\{\theta=\text{const},z=\text{const}\})$ to itself, we may view it
as a contactomorphism between the structures 
$\xi_{1/2}|_{\psi_v({\cal C})- (D_v\cup\Gamma)}$ and
$\xi_1|_{\psi_v(\cal C)-\Gamma}$. Since moreover each of
the maps $g_v$ coincides with the identity close to the boundary
of $\psi_v(\cal C)$, the system of these maps can be extended
by the identity to a contactomorphism $g$ between the structures
$\xi_{1/2}|_{M-(\cup D_v\cup\Gamma)}$ and $\xi_1|_{M-\Gamma}$.
Together with Step 1, this gives an isotopy between the contact structures
$\xi|_{M- F}$ and $\xi_1|_{M-\Gamma}$, and therefore
the composition $\xi_t:t\in[0,1]$ of the deformations from Steps 1 and 2
satisfies the assertions of the proposition.

(2) The second assertion in \ref{ReductionDiffeoProposition} follows from the fact
that a diffeomorphism $h:M- F\to M-\Gamma$ is in fact obtained
as the final map in a smooth 1-parameter family of embeddings of the complement
$M- F$ into $M$, starting at the inclusion map. 
This 1-parameter family can be associated to controlled
reductions described in the proof of part (1). Restriction of this 1-parameter
family to the preimage $h^{-1}(M- W)$ then yields an isotopy of embeddings
of a submanifold $h^{-1}(M- W)$ into $M$, and this can be extended
to an isotopy of diffeomorphisms of $M$ by the standard results concerning isotopy
extensions. 

\end{proof}

\begin{proof} [Proof of Theorem \ref{TightReductionTheorem}] Call any disc $\Delta$ as in 
Definition \ref{SigmaConfoliationTightDef} an {\it overtwisting disc} of a $\sigma$-confoliation
(compare \cite{ET:Confoliations}). Then tightness boils down to nonexistence
of an overtwisting disc. To prove Theorem \ref{TightReductionTheorem} we shall show
that such a disc exists for $\xi$ iff it exists for $\xi'$.  This follows fairly directly from
Proposition \ref{ReductionDiffeoProposition} as follows.

Assume first that there is an overtwisting disc for $\xi$, and call it $\Delta$.
Since the boundary $\partial\Delta$ is disjoint with the maximal integral surface $F$
for $\xi$, we can choose a neighborhood $U$ of $F$ disjoint from $\partial\Delta$.
By Proposition \ref{ReductionDiffeoProposition}, it is then possible
to deform $\xi$ into a contact structure $\xi''$
by a sequence of controlled reductions with supports contained in $U$, so that
$\xi$ and $\xi''$ coincide outside $U$. Then $\Delta$ is clearly an overtwisting
disc for $\xi''$. However, since due to Theorem 1.3 the contact structures $\xi''$
and $\xi'$ are isotopic, it follows that there exists an overtwisting disc for $\xi'$.

Conversely, assume that $\Delta$ is an overtwisting disc for $\xi'$.
By Theorem 1.3, we may assume that $\xi'$ is obtained from $\xi$ by a sequence
of controlled reductions as in Proposition \ref{ReductionDiffeoProposition},
with a core $\Gamma$ and a contactactomorphism $h$ as in the statement of this
proposition. Now, note that if the boundary curve $\partial\Delta$ is deformed
within the class of smooth closed embedded curves tangent to $\xi'$, then it
still bounds an overtwisting disc for $\xi'$. Moreover, it is always possible to deform
$\partial\Delta$ as above so that after deformation it is disjoint from the
graph $\Gamma$. Therefore we may assume that the boundary $\partial\Delta$
of an overtwisting disc $\Delta$ for $\xi'$ is disjoint from $\Gamma$.

Let $h:M- F\to M-\Gamma$ be a contactomorphism as prescribed
by Proposition \ref{ReductionDiffeoProposition}. Then the curve $h^{-1}(\partial\Delta)$
is clearly disjoint with $F$ and tangent to $\xi$. 
Choose a regular neighborhood $W$ of the graph $\Gamma$ disjoint with $\partial\Delta$.
Let $H:M\to M$ be a diffeomorphism whose restriction to the preimage 
$h^{-1}(M- W)$ coincides with $h$, as prescribed by part (2)
of Proposition \ref{ReductionDiffeoProposition}. Then the disc
$H^{-1}(\Delta)$ is easily seen to be transverse to $\xi$ at its boundary
and hence it is an overtwisting disc for $\xi$.  \end{proof}

\section{Contaminations and splitting}\label{Splitting}

In this section we give a precise definition of a contamination as an
object with a smooth structure. We also discuss various operations
on pure contaminations, especially an operation of splitting. A special
case of this operation, applied to contact structures, was described
in the introduction.

\begin{figure}[ht]
\centering
\scalebox{1.0}{\includegraphics{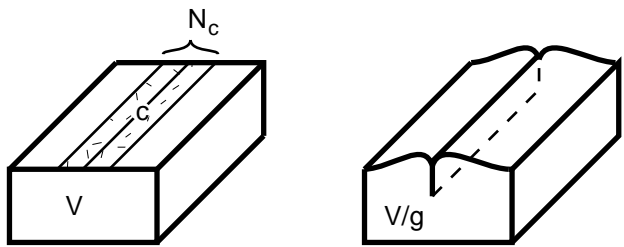}}
\caption{\small The smooth germ at the cusp locus.}
\label{ContactSmoothGerm}
\end{figure}

A {\it manifold with inward cusps} is a topological oriented 3-manifold $V$
with boundary $\partial V$ and with distinguished collection ${\cal C}(V)$ of pairwise
disjoint simple closed curves in $\partial V$ called {\it cusp curves}.
Each cusp curve $c\in {\cal C}(V)$ is furthermore equipped with a germ $g_c$
of a homeomorphic identification of the two sides of a thin neighborhood $N_c$ of $c$
in $\partial V$. We assume also that for each $c\in {\cal C}(V)$
the germ $g_c$ extends the identity map on $c$ (intuitively, $g_c$ is a germ
of a glueing map for $V$ which identifies the two sides of the strip
$N_c$ after folding it along $c$), see Figure \ref{ContactSmoothGerm}.

A smooth structure on a manifold with inward cusps $V$ is given by the following data:
\item{(1)} ordinary charts (including boundary charts) at open neighborhoods in $V$
disjoint with ${\cal C}(V)$;
\item{(2)} germs of charts at points of ${\cal C}(V)$ described more precisely as follows.
\trip

\noindent
Let $c\in {\cal C}(V)$ be a cusp curve of $V$ and let $p\in c$. For a homeomorphism $g$
representing germ $g_c$ let $V/g$ be the quotient of $V$ with respect to $g$ viewed
as a glueing map. A {\it germ of a chart at $p$} is a germ of a homeomorphism
$U_p\to U\subset \reals^3$, where $U_p$ is a sufficiently small neighborhood of $p$
in $V/g$ (so that $U_p$ is homeomorphic to an open set in $\reals^3$).   

We of course assume that transition maps for all charts and germs of charts are smooth.
We notice that a smooth structure on $V$ induces smooth structures on the components
of the horizontal boundary $\partial_hV$ and on the cusp curves in ${\cal C}(V)$. 
It then turns out that the germs $g_c$ are smooth. On the
other hand, the boundary $\partial V$ as a whole does not inherit a canonical
smooth structure.  Charts for $\partial V$ at points of cusp curves are not canonically
determined.

With the above definition of a smooth structure on a manifold with inward cusps $V$,
one can speak of various smooth objects on $V$, like functions, 
vector fields, etc.
This justifies the definition of a (smooth)
contamination in $V$ given in the introduction (Definition 
\ref{ContaminationDef}). 
It also becomes clear what
it means for $V$ to be smoothly embedded in a 3-manifold $M$. Strictly speaking, 
an embedding $V\embed M$ is not injective near cusp curves; 
it respects identifications
at thin neighborhoods of cusp curves induced by the germs $g_c$.
Finally, {\it diffeomorphisms} of manifolds with inward cusps (as well as
general smooth maps between them) and
{\it isomorphisms of contaminations} are now well-defined.

Given a branched surface $B$, a fibered neighborhood $V(B)$ of $B$, as described
in the introduction, is a well-defined (up to diffeomorphism) smooth manifold
with inward cusps. Furthermore, if $B$ is embedded in a 3-manifold $M$, its fibered
neighborhood $V(B)$ is also embedded in $M$, uniquely up to isotopy of embeddings.
Similarly, a contamination $\xi$ carried by a branched surface $B$ embedded in $M$
is well-defined only up to isotopy of embeddings $V(B)\embed M$.

\medskip

A {\it (pure) deformation of contaminations} on a fixed manifold with inward cusps $V$ is
a 1-parameter family $\xi_t$ of (pure) contaminations smoothly depending on the parameter $t$.
It is worth noting that pure contaminations connected by a pure deformation (unlike contact structures
on closed manifolds connected by a pure deformation) need not be isomorphic. One of the reasons is that such
contaminations may behave in a significantly different way near the boundary (the rate
of convergence to tangency when approaching $\partial_hV$ is a sensitive invariant
of isomorphism). On the other hand, even the restricted contact structures
in $\hbox{int}(V)$
for contaminations connected by a pure deformation probably need not be isomorphic. 
For our purposes in this paper it is sufficient to view pure contaminations up to an
even weaker relation than pure deformation. To describe this relation,
we need some further definitions.

A {\it $\sigma$-contamination} in a manifold with inward cusps $V$ is a smooth
plane field $\xi$ in $V$ tangent to $\partial_hV$, for which there exists a compact
surface $F$ (with piecewise smooth boundary) smoothly embedded in $V$, 
such that
\item{(1)} $\hbox{int}(F)$ is embedded in $\hbox{int}(V)$, the intersection 
$F\cap\partial V$ is contained in cusp curves ${\cal C}(V)$
and $F$ is tangent to $\partial_hV$ at the intersection;
\item{(2)} the intersection $\partial F\cap {\cal C}(V)$ is a finite number of arcs and
the components of  
$[\partial F\cup {\cal C}(V)]-\hbox{int}[\partial F\cap {\cal C}(V)]$
are smooth closed curves (i.e. the segments of $\partial F\cap\hbox{int}(V)$ 
smoothly attach to the segments of ${\cal C}(V)- F$);
\item{(3)} for each connected component $F_i$ of $F$, 
$\partial F_i- {\cal C}(V)\ne\emptyset$;
\item{(4)} $\xi$ is everywhere tangent to $F$ and it is contact on
$\hbox{int}(V)- F$.
\item{(5)} At points where $F$ meets ${\cal C}(V)$, $\bdry_hV\cup F$ is smooth.
\trip

\noindent
The surface $F$ is called the {\it maximal integral surface} of a $\sigma$-contamination.

Observe that both pure contaminations and $\sigma$-confoliations are special cases of
$\sigma$-contami\-nations. Condition (3) in the above definition corresponds, 
in the definition
of a $\sigma$-confoliation, to the condition that the maximal integral surface
of a $\sigma$-confoliation cannot have a closed component. 
Condition (2) allows cutting $V$ along $F$ to get another manifold
with inward cusps, as will be described later.
Components $F_i$
intersecting ${\cal C}(V)$ smoothly attach to components of $\partial_hV$, making
the entire locus of integrability of $\xi$ a smooth branched surface, improperly embedded in $M$.

\medskip
A {\it $\sigma$-deformation of $\sigma$-contaminations} 
is a smooth 1-parameter family $\xi_t:t\in[a,b]$
of plane fields on a manifold with inward cusps $V$, satisfying the following
conditions:
\item{(1)} for each $t$, $\xi_t$ is a $\sigma$-contamination on $V$;
\item{(2)} there exists a 1-parameter family $h_t:t\in[a,b]$ of diffeomorphisms
of $V$, a partition of the interval
$[a,b]$ into a finite number of nontrivial disjoint consecutive
subintervals $I_1,I_2,\dots,I_k$,
and a sequence $\Sigma_1,\Sigma_2,\dots,\Sigma_k$
of compact surfaces in $V$ such that:

\item{(a)} for each $t\in \intr(I_j)$ we have $h_t(F_t)=\Sigma_j$ 
(where $F_t$ is the maximal integral surface of $\xi_t$);

\item{(b)} for each $1\le j\le k-1$ the surfaces $\Sigma_j$ and $\Sigma_{j+1}$ are distinct, and one of them
is contained in the other;

\item{(c)} if $\Sigma_j\subset\Sigma_{j+1}$ (or conversely
$\Sigma_{j+1}\subset\Sigma_j$), then
the closure cl($(\Sigma_{j+1}-\Sigma_j)$)
(cl($\Sigma_j-\Sigma_{j+1}$) respectively) is
a compact surface with piecewise smooth boundary;

\item{(d)} each component of the larger of the unions
$\Sigma_j\cup {\cal C}(V)$ and $\Sigma_{j+1}\cup {\cal C}(V)$ is obtained from the smaller one by
attaching 1-handles and adding 0-handles, where 1-handles may be attached to ${\cal C}(V)$.
\trip

Notice that a $\sigma$-deformation of $\sigma$-confoliations
(as defined in the introduction) is a special case of a $\sigma$-deformation
for $\sigma$-contaminations.

In analogy with the results of Section \ref{Reduction} we get the following
lemma which shows that it is possible to eliminate the maximal
integral surface in any $\sigma$-contamination, by means of
a $\sigma$-deformation.

\begin{lemma}\label{PurifySigmaContaminationLem}
For each $\sigma$-contamination $\xi$ there exists a $\sigma$-defor\-ma\-tion
which modifies $\xi$ to yield a pure contamination.
\end{lemma}

To prove this lemma one can work with special charts as in Section \ref{Reduction},
performing controlled reductions associated to those charts. This requires
an obvious adaptation of types of special charts to the situations
when the integral surface meets the cusp locus ${\cal C}(V)$. We omit the
details.

By the arguments as in Section \ref{DeformationTheorem}, generalized as indicated above,
one gets also the following extension of the $\sigma$-Deformation
Theorem (Theorem \ref{SigmaDeformationThm}).

\begin{thm}\label{ContaminationSigmaDeformationThm}
Any two pure contaminations which are connec\-ted by a $\sigma$-deformation
of $\sigma$-contaminations are also connected by a pure deformation.
\end{thm}

We now turn to describing an operation of splitting for pure contaminations.
This operation consists of a sequence of two operations:  first, a $\sigma$-deformation which introduces a nonempty maximal integral
surface followed, second,  by cutting on this maximal integral surface. 
More precisely, the cutting is applied to a
given $\sigma$-contamination $\xi$ on a manifold with inward cusps $V$,
with maximal integral surface $F$, and cutting $V$ on $F$
yields a new smooth manifold with inward cusps $V'$. The cusp locus
${\cal C}(V')$ of $V'$ is then equal to 
$[\partial F\cup {\cal C}(V)]-\hbox{int}[\partial F\cap {\cal C}(V)]$
and the germs of identification close to new parts of ${\cal C}(V')$ are induced
by the identification of the two copies of $F$ which reverse the cutting operation. There is also an obvious smooth structure on $V'$
induced from $V$, and the canonical smooth map $s:V'\to V$.
We obtain a plane field $\xi'$
on $V'$ induced by $s$, which is clearly a pure contamination. We will call
$\xi'$ the {\it pure contamination obtained by cutting $\xi$} on
the maximal integral surface.

\begin{definition}\label{SplittinDef} A
{\it splitting} of a pure contamination $\xi$ in $V$ is an operation consisting of a
$\sigma$-deformation of $\xi$ to a $\sigma$-contamination $\xi'$ with 
maximal integral surface $F$, followed by cutting
$\xi'$ (and its support $V'$, which is the same as the support $V$ of $\xi$) to yield a new pure contamination $\xi''$ with
support $V''$.
\end{definition}

If $V''$ is obtained from $V$ by cutting as above,
then an embedding of $V$ in a 3-manifold $M$ induces an embedding of $V''$ in $M$, well-defined up to isotopy of embeddings. 
Thus,
any splitting of a pure contamination embedded in $M$ yields a pure contamination embedded
in $M$.  One case of interest to us occurs when splitting
a contact structure in $M$ yields a pure contamination in $M$ with
product complementary pieces, as defined in the introduction.

The next section contains descriptions of certain splittings used to prove Theorem \ref{SplitToCarryThm}.

\section{Contact structures are carried by branched surfaces}\label{SplitToCarry}

In this section we prove Theorem \ref{SplitToCarryThm}, showing that any contact structure $\xi$ defined on a 3-manifold
$M$ is $\sigma$-deformable to a $\sigma$-confoliation, which, when cut on its maximal integral surface, is a
pure contamination carried by a branched surface.  Recall that the statement of our theorem is actually stronger:  Given any
1-foliation $\Omega$ transverse to $\xi$, we can find a $\sigma$-deformation $\xi_t$ transverse to $\Omega$ at all times
$t$, $0\le t\le 1$, such that the maximal integral surface of $\xi$ cuts leaves of $\Omega$ into intervals.

\begin{proof} {\it Theorem \ref{SplitToCarryThm}.}  We begin by choosing any surface $F$ everywhere transverse to
$\Omega$ which cuts the leaves of $\Omega$ into intervals.  Finding such a surface is not difficult.  The idea of the proof is to
$\sigma$-deform $\xi$ to $\xi'$, introducing integral surfaces near $F$ such that for every intersection of a leaf of
$\Omega$ with $F$ there is an intersection nearby on the leaf of $\Omega$ with the integral surface of $\xi'$.   

If the surface 
$F$ is in a suitable general position with respect to $\xi$, then the plane field $\xi$ induces on $F$ a singular line
field whose integral curves form the {\it characteristic foliation} or {\it induced foliation}.  In fact, the surface $F$ can be
isotoped slightly such that the induced foliation on $F$ is a singular foliation with only elliptic and hyperbolic singularities, all
of these being in the interior of $F$.  The induced foliation on $F$ may have tangencies on $\bdry F$, i.e. points where the induced
foliation is tangent ot $\bdry F$, and these tangencies are of two types as shown in Figure \ref{ContactRectangles}, again assuming a
suitable general position.    

At the singularities
$F$ is tangent to
$\xi$.  For every elliptic singularity of
$F\embed M$ we choose a $\cal C$-chart centered at the singularity.  Recalling ${\cal C}=\{(r,\theta, z)\in \reals ^{3}
: r<1,\ -1<z<1\}$, we choose a $\cal C$-chart such that $\{z=0\}$ is a disc in $F$, the curves $r=\text{const},\
\theta=\text{const}$ are segments in leaves of $\Omega$, and radii $\theta=\text{const},\ z=0$ are contained in leaves of
the induced foliation on
$F$.  Here, abusing notation, we identify $\cal C$ with its image in $M$.  This choice of a $\cal C$-chart is possible by
Fact \ref{CylinderChartFact}.

For every $\cal C$-chart ${\cal C}_i$ associated to an elliptic point on $F$, we may assume after possibly reparametrizing
$r$, that there is a $\sigma$-deformation which introduces an integral disc at $r\le 1/2$, $z=1/2$, i.e. to one side of
$F$.  Also, we may assume that the support of the $\sigma$-deformation is contained in the interior of the $\cal
C$-chart, and is disjoint from $\{z<1/8\}$.

At each hyperbolic singularity of the induced foliation on $F$, we choose a disc $D\embed M$ transverse to $\Omega$ and
tangent to $F$ at the singularity such that the induced foliation on $D$ has exactly one elliptic point at the same
singularity. After possibly replacing $D$ by a smaller concentric disc, we can then choose a $\cal C$-chart such that
$\{z=0\}=D$, the curves $r=\text{const},\
\theta=\text{const}$ are segments in leaves of $\Omega$, and radii $\theta=\text{const},\ z=0$ are contained in leaves of
the induced foliation on
$D$.  Again this is possible by Fact \ref{CylinderChartFact}.  We can further require, again by possibly replacing $D$
by a smaller disc, that each segment $r=\text{const},\
\theta=\text{const}$ intersects each of $F$ and $D$ exactly once, and that $F\cap {\cal C}\subset \{|z|<1/4\}$.  By
reparametrizing $r$, we can suppose that there is a $\sigma$-deformation with support in $z>3/8$ introducing an integral
surface  $r\le 1/2$, $z=1/2$.

For every $\cal C$-chart ${\cal C}_i$ associated to an elliptic or hyperbolic point in $F$, we project the disc $\{z=1/2,
r\le 1/4\}$ vertically along the 1-foliation to a disc in $F$, and we denote by $\Delta$ the union of these discs,
assumed disjoint, in $F$.  In fact, in the above construction we may arrange that there are exactly four tangencies of the induced
foliation with every boundary component of $\bdry (F-\inter\Delta)$ corresponding to a hyperbolic singularity
and no tangencies on any boundary component corresponding to an elliptic singularity.  Further, the tangencies are all of the
same type; locally near a tangency point, a leaf of the induced foliation tangent to the boundary lies inside
$F-\inter\Delta$.   We want the induced foliation on
$F-\inter\Delta$ to be a foliation by intervals.   That is not the case, since the induced foliation might contain closed leaves or
non-compact leaves.  To achieve our goal we will enlarge
$\Delta$ by adding more
$\cal C$-charts. 

For any $p\in (F-\Delta)\embed M$, we can choose a $\cal C$-chart with $\{p\}=\{r=0,\ z=0\}$, such that the disc
$D=\{z=0\}$ becomes a disc transverse to $\Omega$, such that the curves $r=\text{const},\
\theta=\text{const}$ are segments in leaves of $\Omega$, and such that radii $\theta=\text{const},\ z=0$ are contained in
leaves of the induced foliation on
$D$.  As before, we also require, again by possibly replacing $D$
by a smaller disc, that each segment $r=\text{const},\
\theta=\text{const}$ intersects each of $F$ and $D$ exactly once, and that $F\cap {\cal C}\subset \{|z|<1/4\}$.  By
reparametrizing $r$, we can suppose that there is a $\sigma$-deformation with support in $z>3/8$ introducing an integral
surface  $r\le 1/2$, $z=1/2$.  

For a suitable finite collection of $\cal C$-charts as above,  mutually disjoint and also disjoint from $\cal C$-charts
associated to elliptic and hyperbolic singularities on $F$, if we add to $\Delta$ the  discs $\{z=1/2,
r\le 1/4\}$ projected to $F$, we ensure that the induced foliation of $F-\inter\Delta$ has only interval leaves, with
finitely many tangencies at the boundary. 
We can assume that there are exactly two tangencies (of the same type as for 
hyperbolic singularities) at boundary components corresponding to discs of
$\Delta$ for a non-singular point, see Figure \ref{ContactRectangles}.

Now, let $F_0$ be a slightly larger surface than $F$, containing $F$
in its interior. We may assume that the induced foliation on $F_0-\inter{\Delta}$
has the same properties as the one on $F-\inter{\Delta}$, as described above. Then $F-\inter\Delta$ can be covered by finitely many
(rectangular) discs of the form $R_i=\{(x,y)\in \reals^2: x,y\in (-1,1)\}$, with the property that the segments $y=\text{const}$ lie in
leaves of the induced foliation on $F_0$, and contain entire leaves of the induced foliation on $F-\inter\Delta$, see Figure
\ref{ContactRectangles}.  We also require that for each
$p\in F_0$, $p$ belongs to a most two $R_i$'s.

\begin{figure}[ht]
\centering
\scalebox{1.0}{\includegraphics{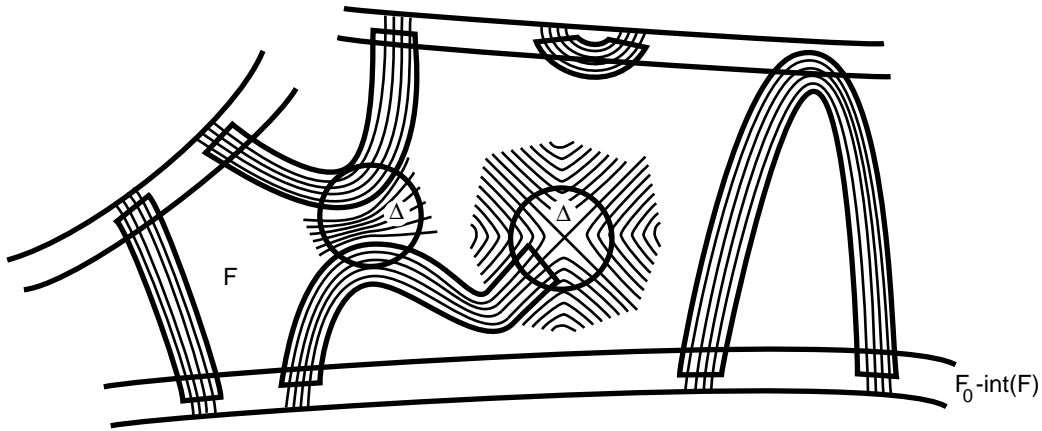}}
\caption{\small Some rectangles $R_i$.}
\label{ContactRectangles}
\end{figure}

For each $R_i$ we now find a $\cal B$-chart ${\cal B}_i=\{(x,y,z)\in \reals^2: x,y,z\in (-1,1)$ 
such that (the images of) the curves $\{x=const,y=const\}$ in ${\cal B}_i$ are contained in the leaves
of $\Omega$; such that projection of ${\cal B}_i$ to $F_0$ along the leaves of $\Omega$ coincides with $R_i$; such that
the $(x,y)$-coordinates in ${\cal B}_i$ are induced from the $(x,y)$-coordinates in $R_i$ via this
projection; and such that the curve $\{y=0\}$ in $R_i$ coincides with (the image of) the curve $\{z=0,y=0\}$
in ${\cal B}_i$.

If the rectangles $R_i$ are sufficiently thin in their $y$-coordinates 
(which can always be assumed),
we can choose the ${\cal B}_i$'s so that for each $i$ the rectangle $R_i$ is contained in (the image of)
${\cal B}_i$, and
for any cylindrical chart ${\cal C}_j$ described above we have  ${\cal B}_i\cap{\cal C}_j\subset \{|z|<1/8\}$, where the
$z$-coordinate is in ${\cal C}_j$.  This means the box charts are away from the support of the $\sigma$-deformations
in the cylinder charts.  

For every ${\cal B}_i$ consider a {\it central rectangle} $P_i$ in ${\cal B}_i$ defined by
$P_i=\{z=0\}$. Observe that due to our assumptions, any rectangle $P_i$ (viewed as a subset of $M$)
intersects $F_0$ exactly at the curve $\{z=0,y=0\}$, and that the two parts of $P_i$ corresponding to
$y>0$ and $y<0$ lie on opposite sides of the surface $F_0$, see Figure \ref{ContactCentralRectangles}.
From this we easily deduce that the central rectangles $P_i$ are pairwise disjoint in $M$,
since the parts that overlap when projected to $F_0$ lie on opposite sides of $F_0$.

\begin{figure}[ht]
\centering
\scalebox{1.0}{\includegraphics{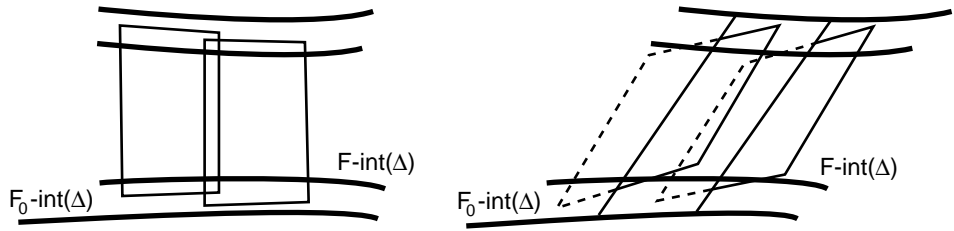}}
\caption{\small Overlapping rectangles $R_i$ and corresponding central rectangles $P_i$.}
\label{ContactCentralRectangles}
\end{figure}

Observe now that, by restricting to a smaller range of the $z$-coordinate around $z=0$, we may choose
smaller ${\cal B}$-charts ${\cal B}_i'$ which contain the $P_i$'s as central rectangles,
and which are pairwise disjoint.
We will show that it is possible to perform a $\sigma$-deformation supported in ${\cal B}_i'$
which produces a disc component of the maximal integral surface contained in the central rectangle $P_i$
and as close to this rectangle as we wish. Note that this is sufficient for completing the proof,
since projections of such discs to $F_0$ still cover $F-\inter\Delta$.

Construction of a $\sigma$-deformation as required is very similar to the construction of
a controlled reduction that eliminates a bridge, as in Section \ref{Reduction}. It can be performed by deforming
the slope function for the contact structure $\xi$ expressed in the corresponding ${\cal B}$-chart
${\cal B}_i'$. The differences are that now we introduce a piece of integral surface
instead of eliminating it, and that the whole process concerns an isolated disc component
of the maximal integral surface instead of a bridge that connects two components of a surface where $\xi$ is integrable.
Nevertheles, the details are completely analogous to those described in Section \ref{Reduction}, and we omit them.
\end{proof}

\section{Three examples in $S^3$}\label{S3Examples}

Given a contact structure $\xi$ in a 3-manifold, if we find a branched surface carrying $\xi$ by the procedure described in
Section \ref{SplitToCarry}, we will almost certainly produce an unnecessarily complicated branched surface.  In practice, one
should look for the simplest and most natural branched surface reflecting the geometry of the contact structure.  However, we
shall see that if the branched surface is not ``fine" enough, it may not encode the tightness or overtwistedness of the contact
structure(s) carried by it.  Clearly, examples of contact structures in $S^3$ are among the most fundamental in the theory, so we
examine, in this section, branched surfaces for some contact structures in $S^3$.  Our explanations will be somewhat informal.

\begin{figure}[ht]
\centering
\scalebox{1.0}{\includegraphics{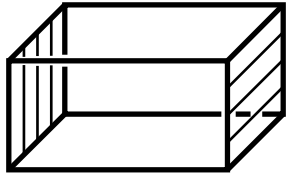}}
\caption{\small The 3-sphere.}
\label{ContactBlockS3}
\end{figure}

We begin with a model of $S^3$ suited to our purposes.  Namely, we regard $S^3$ as a quotient space obtained from
$T^2\times[-1,1]$, where $T^2$ is a torus, by identifying slope $\infty$ circles in $T^2\times \{-1\}$ to points and similarly
collapsing slope
$0$ circles to points in $T^2\times \{1\}$.  In our figures we will represent the sphere as a cube $[-1,1]\times[-1,1]\times[-1,1]$ in
$\reals^3$, see Figure \ref{ContactBlockS3}, where the squares $y=\text{const}$ represent the tori.

\begin{figure}[ht]
\centering
\scalebox{1.0}{\includegraphics{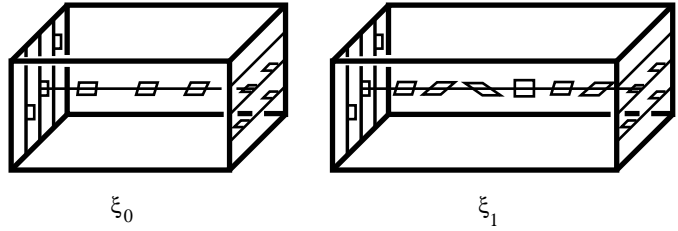}}
\caption{\small Contact structures $\xi_0$ and $\xi_1$.}
\label{ContactXi0Xi1}
\end{figure}

\begin{figure}[ht]
\centering
\scalebox{1.0}{\includegraphics{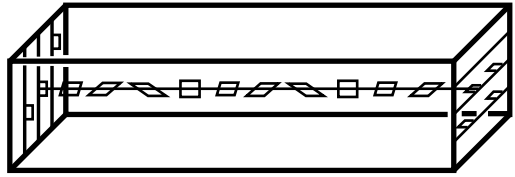}}
\caption{\small Contact structure $\xi_2$.}
\label{ContactXi2}
\end{figure}

Now we describe three positive contact structures on $S^3$.   All three have horizontal lines ${p}\times [-1,1]$ as Legendrians,
$p\in T^2$.  In our model these appear as lines $z=\text{const},\ x=\text{const}$. The first contact structure, which we denote
$\xi_0$ is defined by the property that planes rotate $\pi/2$ going from $y=-1$ to $y=+1$ on each horizontal Legendrian.  The
contact structures $\xi_1$ and $\xi_2$ are similar, but with a rotation of $3\pi/2$ and $5\pi/2$, respectively, on horizontal
Legendrians.  See Figures \ref{ContactXi0Xi1}, \ref{ContactXi2}.  The reader may recognize these structures as the standard
structure $\xi_0$, the once overtwisted structure $\xi_1$, and the twice overtwisted structure $\xi_2$.  We note that adding
additional twists does not give new overtwisted structures, though this may not be obvious to the reader.

We also remark that the Hopf fibration can easily be recognized in our model of $S^3$.  The circle leaves are given by lines of
slope 1 in each plane $y=\text{const}$, namely lines of the form $y=a, z=x+c$, where $c$ is any constant and $-1<a<1$. In
addition, there are two circle leaves obtained by performing the identifications at $y=\pm1$.

\begin{figure}[ht]
\centering
\scalebox{1.0}{\includegraphics{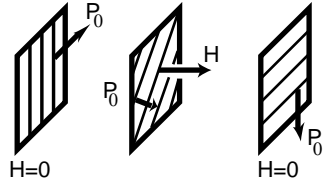}}
\caption{\small The vector fields $P_0$ and $H$.}
\label{ContactField}
\end{figure}

Now that we have our contact structures described explicitly, we will choose a transverse 1-foliation $\Omega_j$ for
each $\xi_j$ to begin our construction of branched surfaces.  The 1-foliations are all constructed in the same way.  First notice
that the induced foliation of
$\xi_j$ on the section
$T^2_y$, $-1<y<1$, is a foliation by constant slope curves (lines in the model), and the slope rotates from $\infty$ to $0$ through
an angle of
$\pi/2+j\pi$.  We define a vector field $P_j$ on our cube model tangent to planes $y=\text{const}$ and an orthogonal unit vector
to the foliation induced on these planes by $\xi_j$.  The vector field induces a vector field on $S^3$ which we also denote
$P_j$.  Let $H=(0,f(y),0)$ be a horizontal vector field, where $f$ is zero at $y=\pm1$ and positive elsewhere, chosen to give a
smooth vector field on the quotient space.  We now let $\Omega_j$ be the 1-foliation obtained by integrating $X_j=P_j+H$, see
Figure \ref{ContactField}.  Notice that the figure applies only for $\xi_0$ and $\xi_2$, not $\xi_1$, because $P$ is opposite for
$\xi_1$ at the right end of the model.

The reader can verify that $\Omega_j$ has two circle leaves, the circles obtained by performing identifications at $y=\pm1$, and
all other (line) leaves limit on both of these closed curves.  Regarding the leaves as orbits of the flow determined by the vector
field, as $t\to\infty$
the line leaves spiral towards the right circle, as $t\to -\infty$ they spiral to the left circle leaf.

\begin{figure}[ht]
\centering
\scalebox{1.0}{\includegraphics{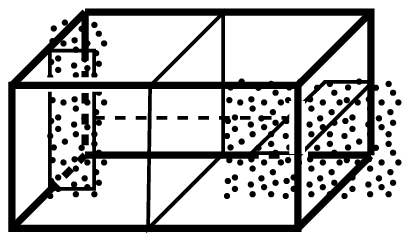}}
\caption{\small Splitting Discs.}
\label{ContactSplittingDiscs}
\end{figure}

We are now ready to begin constructing branched surfaces.  In all cases, we use two discs to split the contact structure.  We
perform a small $\sigma$-deformation of $\xi_j$ to make these integral surfaces.  The discs are shown in Figure
\ref{ContactSplittingDiscs}.  (The dotted Legendrian joining the two discs in the figure should be ignored for now.)  Cutting on the
splitting discs yields a manifold $V_j$ with inward cusps.  The foliation $\Omega_j$ foliates $V_j$ by intervals, giving $V_j$ the
structure $V(\hat B_j)$ of a fibered neighborhood of a branched surface $\hat B_j$.  

\begin{figure}[ht]
\centering
\scalebox{1.0}{\includegraphics{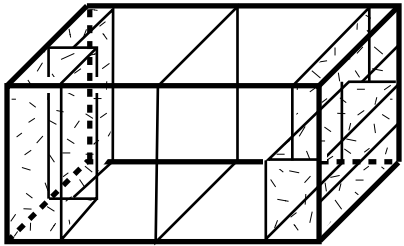}}
\caption{\small Regions corresponding to sectors of branched surface for $\xi_0$ and $\xi_2$.}
\label{ContactS3Sectors}
\end{figure}

After cutting the leaves of $\Omega_j$ on the discs of splitting we get three kinds of segments.  Segments of the first kind have
both ends on the left cutting disc; these segments are contained in the left region shown in Figure \ref{ContactS3Sectors}.  (The
figure is valid for $\xi_0$ and $\xi_2$, not $\xi_1$.) Segments of the second kind have both ends on the right cutting disc; these
segments are contained in the right region shown in Figure \ref{ContactS3Sectors}.  Segments of the third kind have one end on the
left cutting disc and one end on the right cutting disc.  The segments of the third kind all cross the middle torus $\{y=0\}$ and
are contained in the middle region shown in Figure \ref{ContactS3Sectors}.   For $\xi_1$ the right end of Figure
\ref{ContactS3Sectors} must be modified.

\begin{figure}[ht]
\centering
\scalebox{1.0}{\includegraphics{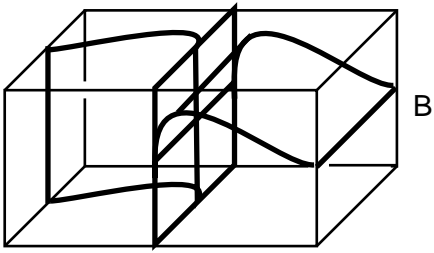}}
\caption{\small Branched surface carrying $\xi_0$ and $\xi_2$.}
\label{ContactS3BranchedInBox}
\end{figure}

For each type of segment in the cut $\Omega_j$, we obtain an embedded subsurface of our branched surface $\hat B_j$.  For segments of the
first (or second) kind, we obtain a disc sector in the branched surface, because the segments can be parametrized by their
endpoints in the corresponding splitting disc.  Every segment of the third kind intersects the
torus
$y=0$ exactly once, so these segments yield a torus.  To understand how the discs are attached to the torus to form a branched
surface, we examine the leaves of $\Omega_j$ passing through the boundaries of the regions in Figure \ref{ContactS3Sectors}.  Due to the
invariance of $\Omega_j$ under the action of the group $T^2$, these are circles, so the discs are attached to the torus on simple
closed curves in the torus.  To understand the sense of branching in the cases $j=0,2$, note that as we move a segment corresponding to a
point in the torus up in the model, the corresponding leaf of $\Omega_j$ suddenly intersects the splitting disc on the right an
additional time.  This means the sense of branching where the right disc sector is attached to the central torus is as shown in Figure
\ref{ContactS3BranchedInBox}.  We denote by $B=\hat B_0=\hat B_2$ the branched surface shown in the figure, which is nothing but a
double Reeb branched surface in $S^3$.  There is another branched surface
$\hat B_1=B_1$  which carries $\xi_1$, and which has the opposite sense of branching on the right side.  The double Reeb branched surface
$B_1$ is shown on the left side of Figure \ref{ContactPseudoTransversalS3}.

\begin{figure}[ht]
\centering
\scalebox{1.0}{\includegraphics{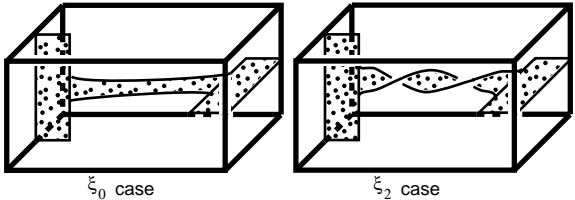}}
\caption{\small Splitting discs joined to get branched surfaces for $\xi_0$ and $\xi_2$.}
\label{ContactSplittingDiscsJoined}
\end{figure}

We are not done because we now have a branched surface $B$ carrying both $\xi_0$ and $\xi_2$, one of which is tight and the other
overtwisted.  In general, our goal is to distinguish tight and overtwisted contact structures using branched surfaces, so this
branched surface is not satisfactory.  By splitting further, we will obtain two branched surfaces $B_0$ and $B_2$, each a
splitting of $B$, such that $B_0$ carries $\xi_0$ and $B_2$ carries $\xi_2$.  
The splitting is obtained by joining 
   the two splitting discs with a thin bridge. This bridge can be produced 
   by a $\sigma$-deformation performed inside a properly chosen 
   $\cal B$-chart around the legendrian segment shown
   in Figure \ref{ContactSplittingDiscs} as a dotted segment.  A $\sigma$-deformation that one can use is,
   roughly speaking, the opposite of the controlled reduction which eliminates a bridge, 
   as described in Section \ref{Reduction}.
The resulting splitting discs are shown in Figure \ref{ContactSplittingDiscsJoined}.

\begin{figure}[ht]
\centering
\scalebox{1.0}{\includegraphics{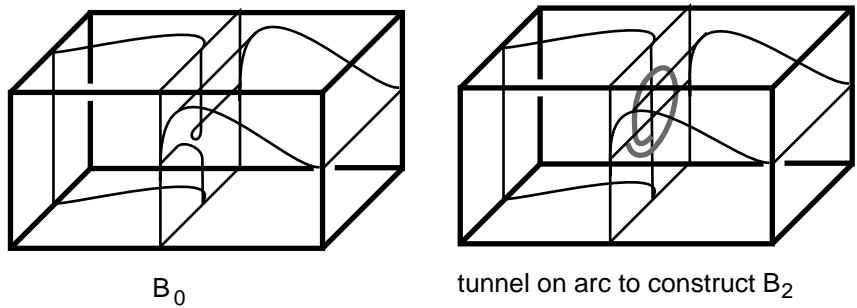}}
\caption{\small  Branched surfaces for $\xi_0$ and $\xi_2$.}
\label{ContactS3B0B2InBlock}
\end{figure}

The effect of extending the splitting surface is to introduce a tunnel joining two parts of the branch locus of $B$.  In the
two cases there are two different ways of tunnelling, one to yield $B_0$ and the other to yield $B_2$, see figure
\ref{ContactS3B0B2InBlock}.  For $\xi_2$, the figure shows the arc along which one tunnels to obtain $B_2$, joining the two
complementary components of $B$.  Another picture of
$B_2$ is shown in Figure \ref{ContactPseudoTransversalS3}, together with a picture of $B_1$, the branched surface carrying
$\xi_1$.

We will show in Section \ref{Pseudotransversal} that $B_1$ and $B_2$ each carry a unique contact structure, which is overtwisted.

\section{Pinching and the proof of Theorem \ref{ContactizationThm}}\label{Pinching}

In this section we deal with the operation of ``pinching" for pure contaminations,
which is something like the reverse of splitting. We show
that using this operation one can recover a contact structure carried 
by a branched 
surface $B$ from any pure contamination carried by $B$ representing this structure
(Theorem \ref{ContactizationThm}).

\medskip

Roughly speaking, a pinching of a pure contamination $\xi$ is an operation which is the composition
of the following operations:
\item{(1)} identification (or glueing) in a manifold $V$ supporting $\xi$
some diffeomorphic parts of two adjacent components 
of the horizontal boundary $\partial_hV$,
thus getting a new manifold with inward cusps $V'$;
\item{(2)} 
pushing forward $\xi$ to $V'$, which gives the plane field  
$\xi'$, for which the glueing locus is the maximal integral surface,
but which is not necessarily smooth at this locus; 
\item{(3)} eliminating the integral surface from  and simultaneously
smoothing $\xi'$. 
\trip

\noindent
These steps require more careful description which we give below.
\medskip

\begin{figure}[ht]
\centering
\scalebox{1.0}{\includegraphics{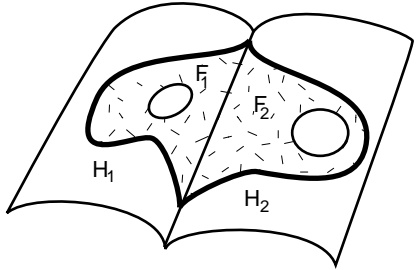}}
\caption{\small Pinching: identify $F_1$ and $F_2$ using glueing diffeomorphism $h$.}
\label{ContactPinch}
\end{figure}

Let $H_1,H_2$ be two adjacent components of the 
horizontal boundary $\partial_hV$.
Suppose that for $i=1,2$, $F_i\subset H_i$ is a surface 
(with piecewise smooth
boundary) satisfying the following conditions, see Figure \ref{ContactPinch}:
\item{(1)} $\partial F_i\subset \hbox{int}(H_i)\cup(H_1\cap H_2)$ and 
the closure of $H_i- F_i$ is a surface (possibly empty) with 
smooth boundary;
\item{(2)} there is a diffeomorphism $h:F_1\to F_2$ which coincides 
with each of the germs $g_c$
close to $c$, for any $c\subset H_1\cap H_2$ (so that in particular, 
$F_1\cap c=F_2\cap c$).
\trip

\noindent
A diffeomorphism $h$ as above will be called a 
{\it glueing diffeomorphism} for $V$.  Strictly speaking, we will always assume that a glueing diffeomorphism $h$ is
equipped with a germ of its extension to a neighborhood of $F_1$ in $H_1$, and that this extending germ coincides with each of the
germs $g_c$ on the overlapping parts.

The quotient space $V/h$ with respect to a glueing diffeomorphism $h$
can be equipped with a structure of a smooth manifold with inward cusps, in such a way
that the quotient map $V\to V/h$ is smooth. Observe however that there is some freedom
of choice for a smooth structure at points of the glueing locus (we will face this problem
later in Lemma \ref{TwoSmoothingsLem}). 
Pushing forward (through the quotient map $V\to V/h$) a pure contamination $\xi$
supported on $V$, we get a plane field $\xi'$ on $V/h$, 
smooth everywhere except possibly at the glueing locus. 
The behaviour of $\xi'$ at the glueing locus motivates the following.

\begin{definition}\label{HalfSmoothContaminationDef}
A {\it semi-smooth $\sigma$-contamination} is a plane field $\xi$ on a manifold with inward
cusps $V$ which satisfies all the properties of a $\sigma$-contamination, except that
at the maximal integral surface $F$ it may be semi-smooth (and not just smooth) in
the following sense: for a sufficiently small neighborhood $U$ of a point of $F$,
$\xi$ is smooth when restricted to any of the two (closed) sides of $U$ with respect
to $F$. In other words, $\xi$ is semi-smooth  
if the induced plane field
in the manifold $V'$ obtained from $V$ by cutting along the maximal integral surface
of $\xi$ is a smooth pure contamination.
\end{definition}

Clearly, the plane field $\xi'$ obtained from a pure contamination $\xi$ by
a glueing $h$ is a semi-smooth $\sigma$-contamination on $V/h$.

A {\it semi-smooth $\sigma$-deformation} is a 1-parameter family $\xi_t$
of semi-smooth $\sigma$-contamina\-tions, smoothly depending on $t$, which
satisfies all other properties of a $\sigma$-deforma\-tion.
It turns out that semi-smooth $\sigma$-contaminations and semi-smooth
$\sigma$-deformations share basic properties with their smooth counterparts.
We formulate them in the two results below. Proofs of these results, 
in principle similar to the smooth case, require a complicated 
technique of moving charts (extending the technique of special charts
as presented in Sections \ref{Reduction} and \ref{DeformationTheorem}). We will present this technique and
sketch the proofs in the next section.

\begin{lemma}\label{HalfSmoothPurifyLem}
For any  semi-smooth $\sigma$-contamination $\xi$ there exists a
semi-smooth $\sigma$-deformation which turns $\xi$ into a pure contamination.
Moreover, for any semi-smooth $\sigma$-deformation $\xi_t$ there exists
a 1-parameter family $\xi_{t,s}:s\in[0,1]$ of semi-smooth $\sigma$-deformations
(so that $\xi_{t,s}$ depends smoothly on both $t$ and $s$) which turns $\xi_t=\xi_{t,0}$
parametrically into pure contaminations $\xi_{t,1}$, for all $t$.
\end{lemma}

\begin{proposition}\label{HalfSmoothSigmaDeformationProp}
Any two contact structures in a closed 3-manifold $M$ which are connected by
a semi-smooth $\sigma$-deformation, are isotopic.
\end{proposition}

We are now ready to give a precise description of an operation of pinching.

\begin{definition}\label{PinchingDef}
{\it Pinching} is an operation which converts a pure contamination $\xi_1$ into another
pure contamination $\xi_2$ by means of the following two steps:
\item{(1)} modifying $\xi_1$ into a semi-smooth $\sigma$-contamination $\xi_1'$
by means of a glueing of the manifold with inward cusps supporting $\xi_1$
and by pushing forward $\xi_1$;
\item{(2)} modifying $\xi_1'$ to a smooth pure contamination $\xi_2$ by a semi-smooth
$\sigma$-deformation (as in Lemma \ref{HalfSmoothPurifyLem}).
\end{definition}

We will use the term ``glueing" not only for manifolds $V$ with inward cusps, but also for contaminations $\xi$ with
support $V$.

Observe that the inverse operation of any splitting is an example of pinching,
but the inverse of a pinching is not easily seen to be a splitting
(this is because each $\sigma$-deformation is a semi-smooth $\sigma$-deformation,
but not vice versa). Nevertheless, we have the following.

\begin{thm}\label{SplitPinchThm}
Two contact structures on a closed 3-manifold $M$ which are connected by a sequence
of splittings and pinchings are isotopic.
\end{thm}

Notice that Theorem \ref{ContactizationThm} of the introduction is a special case of 
the above theorem. The proof, which is given at the end of the section, requires many preparatory results. We start with a result
explaining how the choice of the smooth structure in the supporting manifold
influences the semi-smooth $\sigma$-contamination obtained by glueing.

\begin{lemma}\label{TwoSmoothingsLem}
Let $\xi$ be a pure contamination on a manifold with inward cusps $V$ and let
$h$ be a glueing diffeomorphism for $V$. Let $(V/h)_1$ and $(V/h)_2$ be
two smooth structures of a manifold with inward cusps for the quotient $V/h$.
Denote by $\xi_1$ and $\xi_2$ the semi-smooth $\sigma$-contaminations induced by $\xi$
in $(V/h)_1$ and $(V/h)_2$ respectively. Then there is a semi-smooth
$\sigma$-deformation which turns $\xi_1$ into $\xi_2$.
\end{lemma}

\begin{proof}
Denote by $i:(V/h)_1\to(V/h)_2$ the identity map on $V/h$ and notice that 
it is a semi-smooth isomorphism of $\xi_1$ and $\xi_2$, i.e. $i^*(\xi_1)=\xi_2$.
On the other hand, $(V/h)_1$ and $(V/h)_2$ are diffeomorphic (as smooth manifolds
with inward cusps) through a diffeomorphism $f:(V/h)_1\to(V/h)_2$ which is the identity
on the glueing locus and which approximates $i$. Then there exists an isotopy 
$f_t:(V/h)_1\to(V/h)_2$, $t\in[0,1]$, of semi-smooth diffeomorphisms $f_t$
(where the semi-smoothness reflects the bad behaviour at the glueing locus of $h$),
with $f_0=f$ and $f_1=i$. Put $\eta_t=f_t^*\xi_1$ and note that $\eta_t:t\in[0,1]$
is a semi-smooth $\sigma$-deformation connecting $\eta_0=f^*\xi_1$ and
$\eta_1=i^*\xi_1=\xi_2$. Since $f^*\xi_1$ is isomorphic to $\xi_1$, the lemma follows.
\end{proof}

The next lemma deals with glueings (of pure contaminations with product complementary pieces)
along the entire horizontal boundary, which yield semi-smooth 
$\sigma$-confolia\-tions.

\begin{lemma}\label{TwoGlueingsLem} 
Let $\xi$ be a pure contamination in a closed 3-manifold $M$, and suppose that
$\xi$ has product complementary pieces. Let $h_1,h_2$ be two glueing diffeomorphisms
that collapse all components of the complement, so that the quotients
$V/h_1$ and $V/h_2$ are closed topological manifolds homeomorphic to $M$.
Let $(V/h_1)_0$ and $(V/h_2)_0$ be smooth structures 
in the corresponding quotients, and denote by $\xi_1,\xi_2$ the   
semi-smooth $\sigma$-confoliations induced from $\xi$ in $(V/h_1)_0$ and $(V/h_2)_0$
respectively. Then $\xi_1$ and $\xi_2$ are connected by a semi-smooth
$\sigma$-deformation.
\end{lemma}

\begin{proof}
Since both glueing diffeomorphisms $h_1,h_2$ collapse all product complementary
pieces, there is an isotopy $h_t:t\in[1,2]$ between them, through glueing 
diffeomorphisms $h_t$. It is possible to choose smooth structures $(V/h_t)_0$
in the quotients $V/h_t$, for $t\in[1,2]$, so that the corresponding family $\xi_t'$ of semi-smooth 
$\sigma$-confoliations induced by $\xi$ represents a
semi-smooth
$\sigma$-deformation. Since, by Lemma \ref{TwoSmoothingsLem}, the semi-smooth
$\sigma$-confoliations $\xi_i$ and $\xi_i'$ are connected by a semi-smooth
$\sigma$-deformation, for $i=1,2$, the proof is finished.
\end{proof}

The next lemma directly prepares the ground for proving Theorem \ref{SplitPinchThm}.

\begin{lemma}\label{GluePinchGlueLem}
Let $\xi_2$ be a pure contamination with product complementary pieces in a closed
3-manifold $M$, and suppose that it is obtained by pinching from another such
pure contamination $\xi_1$. Let $\xi_1'$ and $\xi_2'$ be any two semi-smooth
$\sigma$-confoliations in $M$ obtained by glueing from $\xi_1$ and $\xi_2$
respectively. Then $\xi_1'$ and $\xi_2'$ are connected by a semi-smooth
$\sigma$-deformation.
\end{lemma}

\begin{proof}
Denote by $\eta$ the semi-smooth $\sigma$-contamination obtained after the glueing
part of the pinching from $\xi_1$ to $\xi_2$ (before eliminating the maximal integral
surface and smoothing). Then $\eta$ and $\xi_2$ are connected by
the semi-smooth $\sigma$-deformation occurring in the remaining part of this pinching.
We denote this semi-smooth $\sigma$-deformation by $\eta_t:t\in[0,1]$, where
$\eta_0=\eta$ and $\eta_1=\xi_2$. We apply an operation of glueing parametrically 
to $\eta_t$, thus getting $\eta_t'$, a family of semi-smooth $\sigma$-confoliations
being a semi-smooth $\sigma$-deformation. By Lemma \ref{TwoGlueingsLem},
there is a semi-smooth $\sigma$-deformation between $\eta_1'$ and $\xi_2'$. On the 
other hand, since $\eta_0'$ is obtained from $\xi_1$ by the two succesive glueings
$\xi_1\to\eta=\eta_0$ and $\eta_0\to\eta_0'$, we may view it as obtained from $\xi_1$
by a single glueing. Applying Lemma \ref{TwoGlueingsLem} again, we get that
there is a semi-smooth $\sigma$-deformation between $\xi_1'$ and $\eta_0'$.
Composing the above semi-smooth $\sigma$-deformations we connect $\xi_1'$ with $\xi_2'$,
and the lemma follows.  \end{proof}

\begin {proof} {\it Theorem \ref{SplitPinchThm}.}  Suppose that $\xi_0,\xi_n$ are contact structures in a closed 3-manifold $M$
connected by a sequence $\xi_0,\xi_1,\dots,\xi_{n-1},\xi_n$ of pure contaminations,
so that $\xi_{i+1}$ is obtained from $\xi_i$ either by splitting or by pinching,
for $i=0,1,\dots,n-1$. For each $0\le i\le n$, let $\xi_i'$ be a semi-smooth
$\sigma$-confoliation in $M$ obtained from $\xi_i$ by glueing (we also assume that
$\xi_0'=\xi_0$ and $\xi_n'=\xi_n$). If $\xi_{i+1}$ is obtained from $\xi_i$
by pinching, Lemma \ref{GluePinchGlueLem} implies that there is a semi-smooth
$\sigma$-deformation between $\xi_i'$ and $\xi_{i+1}'$. The same is true if $\xi_{i+1}$ is obtained from $\xi_i$ by splitting,
since then we may view $\xi_i$ as obtained from $\xi_{i+1}$ by pinching.
Thus $\xi_0$ and $\xi_n$ are connected by a semi-smooth
$\sigma$-deformation. Applying Proposition 
\ref{HalfSmoothSigmaDeformationProp} we conclude that the contact structures
$\xi_0$ and $\xi_n$ are isotopic.
\end{proof}

\section{Moving charts and properties of semi-\break smooth $\sigma$-deformations}\label{MovingCharts}

In this section we sketch the proofs of Lemma
\ref{HalfSmoothPurifyLem} and Proposition \ref{HalfSmoothSigmaDeformationProp}.
These proofs go basically along the same lines as proofs of the related results
in Sections \ref{Reduction} and \ref{DeformationTheorem}, where we
use special charts corresponding to bridges (or 1-handles) and discs
(0-handles) of some handle decomposition of the maximal integral surface.
The main difference is that here we not only need to eliminate 
the maximal integral surface, but also we need to get smoothness from the semi-smoothness
in the process. This forces us to work with ``moving" charts,
and not just with the fixed special charts as in Sections \ref{Reduction} and \ref{DeformationTheorem}.
 
To avoid inessential technical complications we will restrict the exposition
to the case of semi-smooth $\sigma$-confoliations,
and skip the details that have to be added in the general case.
We start by introducing semi-smooth analogues of special charts.

\medskip
Let $B$ be a bridge in the maximal integral surface $F$ of a semi-smooth
$\sigma$-confoliation $\xi$ in a 3-manifold $M$ and let $X$ be a vector field in $M$
transverse to $\xi$. Recall that ${\cal B}=\{(x,y,z):|x|<1,\,|y|<1,\,|z|<1\}$.
It is not difficult to get a smooth chart $\varphi:{\cal B}\to M$ satisfying the following
conditions:
\item{(1)} images by $\varphi$ of the lines $\{x=\text{const},\,y=\text{const}\}$ in $\cal B$
are contained in the integral curves of the vector field $X$;
\item{(2)} 
$\varphi^{-1}(B)=\{(x,y,z)\in{\cal B}:z=0,\,|x|\le1-\epsilon,\,
\phi_1(x)\le y\le\phi_2(x)\}$ and
$\varphi^{-1}(F)=\varphi^{-1}(B)\cup\{(x,y,z)\in{\cal B}:z=0,\,|x|\ge1-\epsilon\}$.
\trip

\noindent
We will call a chart $\varphi$ as above an {\it $X$-chart for $B$}.

\begin{definition}\label{HalfSmoothChartDef}
With notation as above, a {\it semi-smooth $\cal B$-chart for $\xi$ consistent with
an $X$-chart $\varphi$ for $B$}
is an injective continuous map $\psi:{\cal B}\to M$ 
of the form $\psi(x,y,z)=\varphi(x,y,g(x,y,z))$, with $g(x,y,0)=0$, such that

\item{(1)} denoting ${\cal B}^+=\{(x,y,z)\in{\cal B}:z\ge0\}$ 
and ${\cal B}^-=\{(x,y,z)\in{\cal B}:z\le0\}$, 
the restrictions $\psi^{\pm}=\psi|_{{\cal B}^{\pm}}$ of $\psi$ 
to the sets ${\cal B}^\pm$ respectively, are smooth embeddings;

\item{(2)} the images by $\psi$ of the lines $\{x=\text{const},\,z=\text{const}\}$
in $\cal B$ are tangent to $\xi$.
\end{definition}

Notice that a semi-smooth $\cal B$-chart as above satisfies all properties
of an ordinary $\cal B$-chart and all properties B of Section \ref{Reduction}, except
smoothness at the set $\{(x,y,z)\in{\cal B}:z=0\}$.  
Observe also that a real function $g:{\cal B}\to(-1,1)$ occurring in the definition
clearly has to satisfy the following properties:
\item{(1)} $g$ depends smoothly on variables $x$ and $y$; 
\item{(2)} for all $x,y\in(-1,1)$, 
the one variable function $g(x,y,\cdot):(-1,1)\to(-1,1)$ is a continuous injective
function whose restrictions to intervals $\{-1< z\le0\}$
and $\{0\le z<1\}$ are smooth embeddings.
\trip

\noindent
The existence of a semi-smooth $\cal B$-chart for any $X$-chart $\varphi$ for a bridge $B$ 
follows from a construction that
requires solving a parametric system of ordinary differencial
equations (whose solutions correspond to the curves tangent to $\xi$ mentioned 
in condition (2) in Definition \ref{HalfSmoothChartDef}). We omit the details.

\medskip
Semi-smooth $\cal B$-charts make it possible to express semi-smooth $\sigma$-confoliations
in terms of slope functions. More precisely, let $\xi$ be a semi-smooth $\sigma$-confoliation with maximal integral surface $F$, and let $B$ be a bridge in $F$.
Let $\psi$ be a semi-smooth $\cal B$-chart for $\xi$, consistent with an $X$-chart
$\varphi$ for $B$.  As in the case of ordinary 
$\cal B$-charts, $\xi$ can be expressed in both half-charts $\psi^\pm$
(i.e. restrictions of $\psi$ to the subsets ${\cal B}^\pm$)
as the kernel of a (unique) 1-form $dz-f^\pm(x,y,z)dx$. We will call the pair 
$f^\pm:{\cal B}^\pm\to \reals$ of real functions as above the {\it slope function}
for $\xi$ in the chart $\psi$. Note that, because of semi-smoothness of $\psi$,
in general $\xi$ cannot be expressed in terms of a single smooth slope function.
Notice also that any slope function $(f^\pm)$ for a semi-smooth $\cal B$-chart $\psi$
clearly satisfies the following conditions:
\item{(1)} $f^\pm(x,y,z)=0$ for $(x,y,z)\in(\psi^\pm)^{-1}(F)$;
\item{(2)} $\partial f^\pm/\partial y(x,y,z)>0$ for $(x,y,z)\notin(\psi^\pm)^{-1}(F)$.
\trip

In a similar way, one also defines semi-smooth 
$\cal C$-charts for disc components of $F$,
and slope functions for these charts.

Our goal is to construct a semi-smooth $\sigma$-deformation which eliminates
a bridge from the integral surface $F$ of a semi-smooth
$\sigma$-confoliation $\xi$. As before, we will do this in terms of a deformation
of a slope function. Deformation of a slope function in a $\cal B$-chart 
corresponds to rotating a plane field around the integral curves $\{x=\text{const},\,z=\text{const}\}$.
Since in our semi-smooth case the foliation (of a part of $M$ corresponding to the chart) 
by these curves is not smooth at $F$, we cannot (in general) get a smooth plane field
by modifying the slope function only. Thus, we have to allow
the chart to move, in order to modify the family $\{x=\text{const},\,z=\text{const}\}$ of curves into
a smooth foliation near the bridge. Simultaneously, modifying the slope function
in the moving chart, we will eliminate the bridge from the integral surface.
To carry out this plan we start with the following.

\begin{definition}\label{SmoothXDeformDef}
Let $\psi(x,y,z)=\varphi(x,y,g(x,y,z))$ be a semi-smooth $\cal B$-chart for $\xi$ 
consistent with an $X$-chart $\varphi$. A {\it smoothing $X$-deformation} for $\psi$ is
a 1-parameter family $\psi_t:t\in[0,1]$, depending smoothly on $t$, 
of maps $\psi_t:{\cal B}\to M$ of the form 
$\psi_t(x,y,z)=\varphi(x,y,g_t(x,y,z))$, such that:

\item{(1)} $\psi_0=\psi$ and for each $t\in[0,1]$ $\psi_t$ is a continuous embedding; 

\item{(2)} for each $t\in[0,1]$ and for all $x,y\in(-1,1)$ we have
$\psi_t(x,y,0)=\psi(x,y,0)$, and the restrictions $\psi_t^\pm$ of $\psi_t$
to the subsets ${\cal B}^\pm$ respectively, are smooth embeddings;

\item{(3)} for each $t\ge1/2$ the map $\psi_t$ 
is smooth at points $(x,y,z)\in{\cal B}$
with $z=0$ and $|x|<1-\epsilon+\rho$, for some $0<\rho<\epsilon$;

\item{(4)} there is $\delta>0$ such that for all $t\in[0,1]$ the images
by $\psi_t$ of the curves $\{x=\text{const},\,z=\text{const}\}$ in $\cal B$ are
tangent to $\xi$ at all points $(x,y,z)\in{\cal B}$ with $|x|\ge 1-\delta$ or 
$|y|\ge 1-\delta$ or $|z|\ge 1-\delta$.
\end{definition}

The smoothing $X$-deformation will be our tool for eliminating a bridge from the maximal
integral surface of a semi-smooth $\sigma$-confoliation.
Thus, in condition (3) above we require smoothness for $\psi_t$ with $t\ge1/2$
in order to get smoothness near the eliminated bridge in modified
semi-smooth $\sigma$-confoliation.
The role of the small constant $\rho$ in this condition is related
to the fact that, according to the definition, a semi-smooth 
$\sigma$-confoliation is smooth in some thin neighborhood $N$ of $\partial F$
in the maximal integral surface $F$. 
The set $\{(x,y,z)\in{\cal B}:z=0,\,1-\epsilon<|x|<1-\epsilon+\rho\}$
will correspond to such a neighborhood $N$ after eliminating the bridge $B$
from the maximal integral surface.

Condition (4) means roughly that the foliation by the images
of the curves $\{x=\text{const},\,z=\text{const}\}$ does not change near the frontier
of the chart neighborhoods in $M$, when the parameter $t$ runs through $[0,1]$.
This will make it possible to modify $\xi$ inside the chart neighborhoods only.
The constant $\delta$ occurring in this condition has to be so small that
the preimages $\psi_t^{-1}(B)=\psi^{-1}(B)$ of the bridge $B$ are contained in the set
$(x,y,z)\in{\cal B}:|y|\le 1-\delta$; it also has to be smaller than $\epsilon$.

\medskip
Observe that a smoothing $X$-deformation $\psi_t$ for $\psi$
is fully determined by a 1-parameter family $g_t:t\in[0,1]$
of functions $g_t:{\cal B}\to(-1,1)$, smoothly depending on $t$, satisfying 
the following conditions:

\item{(1)} $g_0=g$ and for all $t\in[0,1]$ and all $x,y\in(-1,1)$ the maps
$g_t(x,y,\cdot)$ are continuous and injective;  

\item{(2)} for all $t\in[0,1]$ and all $x,y\in(-1,1)$ we have 
$g_t(x,y,0)=g(x,y,0)$
and $g_t(x,y,\cdot)$ is a smooth embedding when restricted to intervals $\{-1< z\le0\}$
or $\{0\le z<1\}$;

\item{(3)} for each $t\ge 1/2$ the function $g_t$ is smooth at 
$\{(x,y,z)\in{\cal B}:z=0,|x|<1-\epsilon+\rho\}$, for some $0<\rho<\epsilon$;

\item{(4)} there is $\delta>0$ such that for all $t\in[0,1]$ 
we have that locally $g_t(x,y,z)=g(x,y,u_t(x,z))$ for some functions $u_t$
and for $|x|\ge 1-\delta$ or 
$|y|\ge 1-\delta$ or $|z|\ge 1-\delta$.
\trip

\noindent
Note that these conditions correspond to the conditions (1)-(4) in Definition
\ref{SmoothXDeformDef}.

To eliminate a bridge by means of a semi-smooth $\sigma$-deformation we need
to enrich a smoothing $X$-deformation corresponding to this bridge with a consistent
deformation of the slope function $f^\pm$. To do this, we need further restrictions
on our smoothing $X$-deformation.
The next definition and lemma explain the details of this idea.

\begin{definition}\label{DelicateDef}
With notation as in Definition \ref{SmoothXDeformDef} for all related objects, let $\psi_t$ be 
a smoothing $X$-deformation. Due to property (4) in
this definition, the semi-smooth $\sigma$-confoliation $\xi$ can be expressed
in charts $\psi_t$, near their frontiers, by means of slope functions. More precisely,
put ${\cal B}_\delta:=\{(x,y,z)\in{\cal B}:|x|\ge1-\delta,|y|\ge1-\delta,|z|\ge1-\delta\}$,
${\cal B}_\delta^+:=\{(x,y,z)\in{\cal B}_\delta:z\ge0\}$ and
${\cal B}_\delta^-:=\{(x,y,z)\in{\cal B}_\delta:z\le0\}$. Then there is the unique
family $h_t^\pm:{\cal B}_\delta^\pm\to \reals$ of functions such that $\xi$ is equal
to the kernel $\ker[dz-h_t^\pm(x,y,z)dx]$ in the chart $\psi_t$.

We say that $\psi_t$ is {\it slope-extendable}, if for each $t\in[0,1]$ and for each $(x,z)$
not in $\{ z=0, |x|\le 1-\epsilon \}$ we have 
$h_t^\pm(x,-1+\delta,z)<h_t^\pm(x,1-\delta,z)$.
\end{definition}

Observe that the inequality as above is automatically satisfied by the slope function
$f^\pm$ for $\xi$, due to monotonicity with respect to the variable $y$. 
Thus, this inequality is a necessary condition for partial slope functions $h_t^\pm$ 
to be extendable to appropriate slope functions (for the constructed deformation) 
in the whole $\cal B$.

As we shall see below, slope-extendability makes it possible to extend functions $h_t^\pm$
to a deformation of slope functions associated to $\psi_t$,
which describes a semi-smooth $\sigma$-deformation eliminating a bridge from 
the maximal integral surface. It is thus an important question whether slope-extendable
smoothing $X$-deformations exist. The answer is given in the following.

\begin{lemma}\label{DelicateExistLem}
Suppose we are given a semi-smooth $\sigma$-confoliation $\xi$, a bridge $B$
in the maximal integral surface for $\xi$, a vector field $X$ transverse to $\xi$,
and an $X$-chart $\varphi$ for $B$. Then we can find a semi-smooth $\cal B$-chart $\psi$ 
for $B$, consistent with $\varphi$, for which there exists a slope-extendable smoothing
$X$-deformation $\psi_t$.
\end{lemma}

We delay the proof of above lemma, and first show how to apply it to eliminate
a bridge from a semi-smooth $\sigma$-confoliation $\xi$.

\begin{proposition}\label{SemiSmoothBridgeEliminationProp}
Let $\xi$ be a semi-smooth $\sigma$-confoliation with the maximal integral surface $F$,
and let $B$ be a bridge in $F$. Then there exists a semi-smooth $\sigma$-deformation
$\xi_t:t\in[0,1]$, with $\xi_0=\xi$, which eliminates $B$ from $F$, i.e. for which
the maximal integral surface of $\xi_1$ is equal to the closure of $F- B$ 
and there is 
only one discontinuity point for the maximal integral surface of $\xi_t$ 
(where the surfaces changes
from $F$ to the closure of $F- B$).
\end{proposition}

\begin{proof} Choose a vector field $X$ transverse to $\xi$ and an $X$-chart $\varphi$
for $B$. By Lemma \ref{DelicateExistLem}, there exists a semi-smooth $\cal B$-chart $\psi$
for $\xi$ consistent with $\varphi$, which admits a slope-extendable smoothing $X$-deformation
$\psi_t$. Let $h_t^{\pm}:{\cal B}_\delta^\pm \to \reals$ be the partial slope
functions for $\psi_t$. By the slope-extendability of $\psi_t$, we can extend $h_t^\pm$
to a family $f_t^\pm:{\cal B}^\pm \to \reals$ of smooth (pairs of) functions, smoothly 
depending on $t\in [0,1]$, satisfying the following conditions:

\item{(1)} for each $t\in[0,1/2]$ we have $\partial f_t^\pm/\partial y(x,y,z)>0$ for  $(x,y,z)\in \psi_t^{-1}(M- F)$ and $\partial f_t^\pm/\partial y(x,y,z)=0$ for
$(x,y,z)\in \psi_t^{-1}(F)$;

\item{(2)} for each $t\in(1/2,1]$ we have $\partial f_t^\pm/\partial y(x,y,z)>0$ for  $(x,y,z)\in \psi_t^{-1}(M- cl(F- B))$ and $\partial f_t^\pm/\partial y(x,y,z)=0$ for $(x,y,z)\in \psi_t^{-1}(cl(F- B))$;

\item{(3)} for each $t\in[0,1]$ the 1-forms $dz-f_t^\pm(x,y,z)dx$ in ${\cal B}^\pm$,
after pushing them forward through $\psi_t$, determine a well-defined 1-form $\omega_t$
in $\psi_t({\cal B})$
(i.e. the two pushed-forward 1-forms on the sets $\psi_t({\cal B}^\pm)$ respectively coincide on the intersection of these sets, thus
giving a sort of semi-smooth 1-form);

\item{(4)} for each $t\in[0,1/2]$ the 1-form $\omega_t$ is smooth at all points
$\psi_t(x,y,0)$ outside $F$ or in a small neighborhood $N$ of $\partial F$ in $F$;

\item{(5)} for each $t\in(1/2,1]$ the 1-form $\omega_t$ is smooth at all points
$\psi_t(x,y,0)$ outside $cl(F- B)$ or in a small neighborhood $N$ 
of $\partial[cl(F- B)]$ in $cl(F- B)$.
\trip

\noindent
In fact, slope-extendability is necessary (and sufficient) for conditions (1) and (2) only;
the other conditions can be fulfilled without referring to this property.
In particular, conditions (4) and (5) can be achieved due to the smoothness
of the foliations by the images under $\psi_t$ of the curves $\{x=\text{const},z=\text{const}\}$
in $\cal B$, at the corresponding points.

Conditions (1)-(5) above imply that putting 
$$
\xi_t(p)=\left\{\begin{array}{cc}\xi(p) & \hbox{for $p\notin\psi_t({\cal B})$}\\
 \ker\omega_t & \hbox{for $p\in\psi_t({\cal B})$} \end{array}\right.
$$

\noindent we get a semismooth $\sigma$-deformation $\xi_t$ that eliminates the bridge $B$
as required, thus finishing the proof.
\end{proof}

\begin{proof} {\it (Lemma \ref{DelicateExistLem}.)}
Choose a semi-smooth $\cal B$-chart $\psi$ consistent with $\varphi$ so that it satisfies
the following additional condition: images by $\psi$ of the lines $\{y=0, z=\text{const}\}$ in
$\cal B$ are everywhere tangent to $\xi$. It is not difficult to convince oneself that this can always be done. Note that in such a $\cal
B$-chart
$\psi$ the slope function $f^\pm$ for $\xi$ satisfies the condition: $f^\pm(x,0,z)=0$ for all $x,z\in(-1,1)$.
Likewise, due to monotonicity of $f^\pm$ with respect to $y$ outside the integral 
surface $F$, we have $f^\pm(x,1-\delta,z)>0$ and $f^\pm(x,-1+\delta,z)<0$ for any
$(x,z)$ outside $\{z=0,|x|\ge1-\epsilon\}$.
We will show the existence of a slope-extendable smoothing $X$-deformation $\psi_t$ for any
semi-smooth $\cal B$-chart $\psi$ as above.

We will concentrate on constructing the appropriate final map $\psi_1$ for
a slope-extendable smoothing $X$-deformation $\psi_t:t\in[0,1]$. We will comment on 
constructing the whole family $\psi_t:t\in[0,1]$ only after constructing $\psi_1$.

Choose rectangles $R_+:=\{y=1-\delta\}$, $R_-:=\{y=-1+\delta\}$ and $R_0:=\{y=0\}$
in $\cal B$, and view them as subsets of $M$ expressed in the chart $\psi$.
Consider characteristic foliations $\zeta_+,\zeta_-,\zeta_0$ of $\xi$ at $R_+,R_-,R_0$
respectively. Then, due to additional condition on $\psi$ as above, leaves of $\zeta_0$
all have the form $\{z=\text{const}\}$. Moreover, again due to the monotonicity with respect
to $y$, the slopes of $\zeta_+$ ($\zeta_-$)
are positive (negative respectively) everywhere except at $\{z=0, |x|\ge 1-\epsilon\}$,
where they are zero, see Figure \ref{ContactRectanglesCharac}.

%\vskip5truecm
%\centerline{Figure RectanglesCharacteristicFig.}

\begin{figure}[ht]
\centering
\scalebox{1.0}{\includegraphics{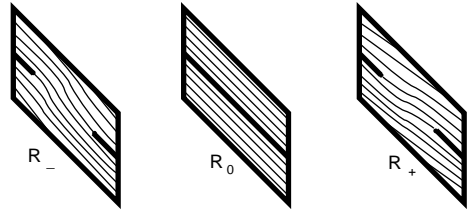}}
\caption{\small Characteristic foliation on rectangles.}
\label{ContactRectanglesCharac}
\end{figure}

Consider now smooth coordinates $x',y',z'$ in $\psi({\cal B})$ with 
$x'=x$, $y'=y$ and $z'=z'(x,y,z)$, such that $\{z=0\}=\{z'=0\}$ (i.e. $z'(x,y,0)=0$).
In these coordinates $\zeta_0$ is smooth everywhere except at $\{z'=0\}$, while
$\zeta_+$ and $\zeta_-$ are smooth everywhere except at $\{z'=0, |x'|\ge 1-\epsilon+\rho\}$.
Let $\chi_+$ and $\chi_-$ be the 1-foliations in $R_+,R_-$ respectively, with all leaves of
form $\{z=\text{const}\}$ in the initial chart $\psi$. 
Then the slopes of $\chi_+$ ($\chi_-$) are smaller (respectively bigger) than the slopes
of $\zeta_+$ ($\zeta_-$ respectively), except at $\{z'=0,|x'|\ge 1-\epsilon\}$, where
the slopes of both foliations are equal to zero (horizontal). Moreover, the foliations
$\chi_+$ and $\chi_-$ expressed in $(x',z')$-coordinates are semi-smooth
(smooth when restricted to $\{z'\ge0\}$ and to $\{z'\le0\}$, but in general not smooth
at $\{z'=0\}$).

A crucial step in constructing a chart $\psi_1$ consists of replacing the foliations
$\chi_+$ and $\chi_-$ by foliations $\chi'_+$ and $\chi'_-$ respectively,
satisfying the following conditions:

\item{(1)} $\chi'_+$ ($\chi'_-$) coincides with $\chi_+$ ($\chi_-$ respectively)
at the subset in $R_+$ ($R_-$ respectively) described as $\{ |x|\ge1-\delta
\hbox{ or } |z|\ge1-\delta \}$ in the coordinates $(x,z)$ from the chart $\psi$;

\item{(2)} the curve $\{z'=0\}$ in $R_+$ ($R_-$) is still a leaf of $\chi'_+$
($\chi'_-$ respectively);

\item{(3)} the foliations $\chi'_\pm$ are smooth everywhere except possibly at 
$\{z'=0,|x'|\ge 1-\epsilon+\rho \}$;

\item{(4)} $\chi'_+$ ($\chi'_-$) is still transverse to (the restriction to $R_+$,
or $R_-$ respectively, of) the vector field $X$, and its slopes 
are still smaller (bigger respectively) than the slopes of $\zeta_+$ ($\zeta_-$), except
at $\{z'=0, |x'|\ge 1-\epsilon \}$, where these slopes are equal.
\trip

\noindent
Without going into complete detail, we observe that foliations $\chi'_+,\chi'_-$ as described above exist by $C^0$-approximation of a
continuous foliation by a smooth one  (here we view foliations as line fields).

Choose a 2-foliation $\cal F$ in $\psi({\cal B})$ transverse to the vector field $X$
such that:

\item{(1)} leaves of $\chi'_+$ and $\chi'_-$ are contained in the leaves of $\cal F$;

\item{(2)} images by $\psi$ of the pieces of lines $\{x=\text{const},z=\text{const}\}$ contained in
${\cal B}_\delta$ are contained in the leaves of $\cal F$; 

\item{(3)} $\cal F$ is smooth except at $\{z'=0, |x'|\ge 1-\epsilon+\rho \}$.
\trip

\noindent
Such a foliation $\cal F$ can easily be constructed due to the properties of the
previously constructed foliations $\chi'_+$ and $\chi'_-$.

Choose a transverse parameter with respect to which $\cal F$ is smooth
as in (3), and take this parameter as $z$-coordinate of the chart $\psi_1$, keeping
the $x$- and $y$- coordinates as in $\psi$.
Then, due to conditions (1)-(3) above, 
the chart $\psi_1$ and its partial slope function $h_1^\pm$ fulfill the relevant
requirements from the definition of a slope-extendable smoothing $X$-deformation.
In particular, condition (1) above and condition (4) for $\chi'_+$ and $\chi'_-$ imply
that  $h^\pm(x,1-\delta,z)>0$ and $h^\pm(x,-1+\delta,z)<0$ for any
$(x,z)$ outside $\{z=0,|x|\ge1-\epsilon\}$, hence the slope-extendability of $\psi_1$.

A construction as above for $\psi_1$ can be performed parametrically, starting with
constructing a parametric family $(\chi'_\pm)_t$ of modified foliations, smoothly connecting
$(\chi'_\pm)_0=\chi_\pm$ with $(\chi'_\pm)_1=\chi'_\pm$. After taking care for 
all the necessary details one gets a slope-extendable smoothing $X$-deformation as required.
We omit further details.  
\end{proof}

We now conclude the proof of Lemma \ref{HalfSmoothPurifyLem} and Proposition
\ref{HalfSmoothSigmaDeformationProp}. To prove the first one, it is sufficient
to note that, similarly to the above, one can define and use moving charts for semi-smooth
$\cal C$-charts in order to construct semi-smooth $\sigma$-deformations that eliminate
a disc component from the maximal integral surface. Extending these methods slightly,
one can also get a parametric elimination of the maximal integral surface from any
semi-smooth $\sigma$-deformation, again using an argument similar to the one in Section \ref{Reduction}.

To prove Proposition \ref{HalfSmoothSigmaDeformationProp} observe that we can
parametrically eliminate the maximal integral surface in such a way, that whenever 
this surface is initially empty, it remains empty during the elimination process.
This clearly leads to an isotopy between any two contact structures connected
by a semi-smooth $\sigma$-deformation, which finishes the proof.

\section{Pseudo-transversal knots} \label{Pseudotransversal}

\medskip
In this section we discuss in detail pseudotransversal knots for branched surfaces
embedded in 3-manifolds. We start with the definition.

\begin{definition}\label{PseudotransversalDef}
Let $B\embed M$ be a branched surface embedded in a 3-manifold $M$, with product
complementary pieces. A {\it pseudotransversal knot} for $B$ is a null-homologous
knot $k$ in $M$ satisfying the following conditions:

\item{(1)} $k$ is disjoint with the branched locus of $B$;

\item{(2)} $k$ intersects $B$ transversely;

\item{(3)} for each connected component $U$ of the complement $M- B$ the
intersection $k\cap U$ is (isotopic rel boundary $\partial U$ to) a {\it braid} in $U$.
\trip

\noindent
Condition (3) above needs more detailed explanation. We say that an intersection
$k\cap U$ is a {\it braid} in $U$, if for some identification of $U$ with the product
$\hbox{int}F\times(0,1)$ each strand (i.e. component) of $k\cap U$ projects to the
factor $(0,1)$ in the product injectively.
\end{definition}

Recall now that for a null-homologous knot $k$ transverse to a plane field $\xi$
in a 3-manifold $M$ one defines the self-linking number $l(k)$ (relative to $\xi$)
as follows. Let $\Sigma$ be a Seifert surface for $k$ in $M$. Then we can find a
non-zero section $\sigma$ of $\xi|_\Sigma$. We then use $\sigma$ to get a parallel copy $k'$
of $k$ by pushing towards $\sigma$. The self-linking number $l(k)$ is then the
linking number of $k$ and $k'$. One easily verifies that this number does not depend
on the choices of $\Sigma$ and $\sigma$ (compare \cite{YE:LegendrianKnots} or \cite{EE:ContactLectures}).

If $k$ is a pseudotransversal knot for $B$, one can define the self-linking number $l(k)$
(relative to $B$) by means of a plane field $\xi_B$ in $M$ associated to $B$ as follows.
At any point of $B$ define $\xi_B$ to be equal to the tangent plane to $B$. Extend $\xi_B$
to the components $U$ of $M- B$ so that, identifying $U$ with 
$\hbox{int}F\times(0,1)$, restriction of $\xi_B$ to $U$ coincides with the foliation
of $U$ by the surfaces $\hbox{int}F\times\{z\}:z\in(0,1)$. This clearly defines
a continuous plane field in $M$ (we do not require smoothness here),
and this field is well-defined up to homotopy of continuous plane fields.
Observe now that, according to the definition of a pseudotransversal knot, a plane field
$\xi_B$ as above can be chosen in such a way that $k$ is everywhere transverse to $\xi_B$.
The self-linking number $l(k)$ relative to $\xi_B$ 
is then uniquely determined (i.e. does not depend on the
choice of $\xi_B$), and we take it as the {\it self-linking number for the pseudotransversal
knot $k$ relative to $B$}.

We now turn to proving the following crucial lemma which directly implies Theorem
\ref{PseudoTransversalThm} of the introduction, by referring to the Bennequin inequality,
mentioned there.

\begin{lemma}\label{TransverseFromPseudoLem}
Suppose that a branched surface $B\embed M$ carries a contact structure $\xi$, and let $k$
be a pseudotransversal knot for $B$. Then there exists a knot $\tilde k$ which is 
everywhere transverse to $\xi$, such that $l(\tilde k)=l(k)$ (where the self-linking
of $\tilde k$ is viewed relative to $\xi$).  
\end{lemma}  

\begin{proof} Let $\xi_1$ be a pure contamination supported by $V(B)$ obtained from $\xi$
by a splitting. Then $k$ determines a knot $k_1$ which intersects $V(B)$ along vertical 
fibres (which correspond to points of transverse intersection of $k$ with $B$)
and whose intersections $k_1\cap U$ with connected components of $M- V(B)$
are the same braids as the corresponding intersections of $k$ with $M- B$.

Fix a component $U$ of the complement $M- V(B)$ and suppose its boundary 
consists of the two copies of a surface $F$, so that $U$ is homeomorphic to
$\hbox{int}F\times(0,1)$. Suppose that the braid $k_1\cap U$ consists of $m$ strings.
This braid can be then viewed as an isotopy of the set of $m$ points in $\hbox{int}(F)$.
This isotopy can be extended to an isotopy of $F$ in itself (rel boundary), i.e., a 1-parameter family of diffeomorphisms of
$f_t\from F\to F$, starting at $\hbox{id}_F$ and equal the identity on $\partial F$ at all times.
We denote by $f_U$ the final diffeomorphism of such a
1-parameter family. Consider then an operation of pinching that collapses $U$, with $f_U$
taken as a glueing diffeomorphism. Observe that after doing such pinchings for all
connected components of $M- V(B)$, the pieces of $k_1\cap V(B)$ glue together to
form a knot $k_2$ transverse to the contact structure $\xi_2$ obtained by these pinchings.

It follows from the description above that the self-linking $l(k_2)$ relative to $\xi_2$
is equal to the self-linking $l(k)$ relative to $B$. One can see this by constructing
a homotopy between the plane field $\xi_B$ associated to $B$ and the contact structure
$\xi_2$, followed by a homotopy between $k$ and $k_2$ through knots transverse to the
homotoped plane fields. On the other hand, according to Theorem \ref{SplitPinchThm},
$\xi_2$ is isomorphic to the initial contact structure $\xi$, and thus the lemma
follows. 
\end{proof}

Our next observation is the following.

\medskip\noindent
\begin{lemma} \label{UniqueCarriedUpToHomotopyLem}
Suppose that $B\embed M$ is a branched surface that carries a contact structure $\xi$,
and let $\xi_B$ be a continuous plane field in $M$ associated to $B$ (as above, in the
definition of the self-linking relative to $B$). Then $\xi$ and $\xi_B$ are homotopic
as plane fields.
\end{lemma}

\begin{proof} Let $\xi_s$ be a pure contamination carried by $V(B)$ obtained from $\xi$
by a splitting. We extend $\xi_s$ to a continuous plane field $\xi_{s,B}$ just as we extended the tangent plane field of $B$ to
$\xi_B$. It is then clear that $\xi_{s,B}$ is homotopic both to $\xi$ (by collapsing the product components
of $M- V(B)$, which corresponds to pinching $\xi_s$ to return to $\xi$) and to $\xi_B$
(by collapsing $V(B)$ along vertical fibres), and thus the lemma follows.
\end{proof}

We summarize above results in the following.

\begin{proposition}\label{AtMostOneProp}
Let $B\embed M$ be a branched surface with product complementary pieces and suppose
that it admits a pseudotransversal knot $k$ which violates the Bennequin inequality. 
Then $B$ carries at most one (up to isotopy) contact structure which, if it exists,
is overtwisted.
\end{proposition}

\begin{proof} According to Theorem \ref{PseudoTransversalThm}, any contact structure
carried by $B$ is overtwisted. On the other hand, it follows from Lemma
\ref{UniqueCarriedUpToHomotopyLem} that any two contact structures carried by $B$
are homotopic as plane fields in $M$. By the classification theorem for overtwisted
contact structures (\cite{YE:OvertwistedClassification}), there is only one overtwisted contact structure
in each homotopy class of plane fields in $M$, hence the Proposition.
\end{proof}

The next result, whose proof uses pseudotransversal knots, strengthens Theorem
\ref{SplitToCarryThm} in the case of overtwisted contact structures.

\begin{proposition}\label{CarryUniqueOvertwistedProp}
For any overtwisted contact structure $\xi$ in a 3-manifold $M$ there exists a
branched surface $B\embed M$ such that $\xi$ is the only contact structure it carries.
\end{proposition}
\medskip\noindent

\begin{proof} To construct a branched surface $B$ as required we refer to the construction
from the proof of Theorem \ref{SplitToCarryThm} in Section \ref{SplitToCarry}, adding
appropriate modifications. 

Let $k$ be a knot in $M$ transverse to $\xi$ and violating the Bennequin inequality (since
$\xi$ is overtwisted, such a knot always exists). Let $\Omega$ be a 1-foliation in $M$
transverse to $\xi$, for which $k$ is a leaf. By applying to $\Omega$ the construction
from the proof of Theorem \ref{SplitToCarryThm} in Section \ref{SplitToCarry}
we get a branched surface $B$ for which $k$ is a pseudotransversal knot, with
the self-linking relative to $B$ equal to its self-linking relative to $\xi$ and thus
violating the Bennequin inequality. It then follows from Proposition \ref{AtMostOneProp} that $\xi$
is the only contact structure carried by $B$, which finishes the proof.
\end{proof}

We will now describe an easily detected property of branched surfaces which implies
the existence of a pseudotransversal knot that violates the Bennequin inequality.

\begin{figure}[ht]
\centering
\scalebox{1.0}{\includegraphics{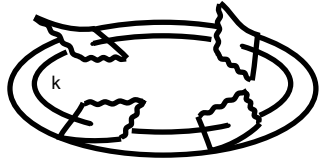}}
\caption{\small Twisted disc.}
\label{ContactTwistedDisk}
\end{figure}

\begin{definition}\label{TwistedDiskDef}
Given a branched surface $B\embed M$, an {\it embedded disc negatively twisted on one side}
in $B$ is a disc $D$ smoothly embedded in $B$, with $\partial D$ transverse to the
branch locus of $B$, with distinguished side, such that all branching of $B$ at branch
curves meeting $\partial D$ on the distinguished side of $D$ have the same sense, a sense we call
negative, as shown in Figure \ref{ContactTwistedDisk}. Observe that if all the branching
as above has the same sense, but different from the one in Figure \ref{ContactTwistedDisk},
then $D$ is not negatively twisted on the distinguished side. Observe also that
for a disc negatively twisted on one side, the branching on the other side
can be arbitrary.
\end{definition}

\begin{figure}[ht]
\centering
\scalebox{1.0}{\includegraphics{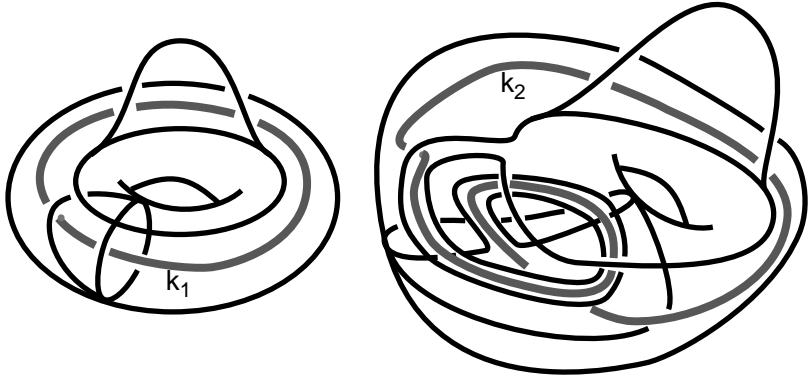}}
\caption{\small Pseudo-transversals in $S^3$.}
\label{ContactPseudoTransversalS3}
\end{figure}

\begin{lemma} \label{TwistedDiskLem}
Suppose that a branched surface $B\embed M$ with product complementary pieces contains
an embedded disc $D$ negatively twisted on one side. 
Then there exists an unknotted curve $k$
pseudotransversal to $B$ with $lk(k)=1$. Since this curve violates the Bennequin inequality,
the branched surface $B$ carries at most one contact structure which, if it exists, is
overtwisted.
\end{lemma}

\begin{proof} In view of Proposition \ref{AtMostOneProp}, it is sufficient to prove the first
assertion. A curve $k$ as asserted is shown at Figure \ref{ContactTwistedDisk}, being
a small shift of $\partial D$ towards the distinguished side.  It can be checked directly
that $lk(k)=1$. 
\end{proof}

As an application of Proposition \ref{AtMostOneProp} to concrete examples we give
the following.

\begin{proposition} \label{S3UniqueCarryProp}
The branched surfaces $B_1,B_2$ in $S^3$ that carry the overtwisted contact structures
$\xi_1,\xi_2$ respectively, as described in Section \ref{S3Examples}, both carry
no other contact structures.
\end{proposition}

\begin{proof} Figure \ref{ContactPseudoTransversalS3} shows the pseudotransversal unknotted
curves $k_1$ and $k_2$ for the branched surfaces $B_1$ and $B_2$ respectively, for which
we have $lk(k_1)=1$ and $lk(k_2)=1$. Since these curves violate the Bennequin inequality,
the proposition follows.
\end{proof}

\section{Invariance of tightness for pure contaminations}\label{ContaminationTightnessInvariance}

In the remainder of the paper we shift our attention to contaminations.
The goal of this section is to clarify the notion of tightness for pure contaminations
and to establish basic classes of deformations of pure
contaminations that preserve tightness. This is important for our long-term goal of using tight
contaminations to study 3-manifolds. 

We start this section by showing that the definition of tightness of a pure
contamination as given in the introduction (Definition 1.14) can be simplified slightly.

Given a manifold with inward cusps $V$ embedded in a 3-manifold $M$, a {\it compressing disc} for
the horizontal boundary $\partial_hV$ is a disc $D$
embedded in $M$, with $\partial D\subset\partial_hV$ and $D$ transversal to $\partial_hV$
at $\partial D$, for which there is no disc $H\subset\partial_hV$ with $\partial H=\partial D$.
We will say that
the horizontal boundary $\partial_hV$ is {\it incompressible} in $M$ if 
there is no compressing disc for it.
Note that, according to this definition, an overtwisting disc $D$ with $\partial D\subset\partial_hV$
for a pure contamiantion
$\xi$ supported by $V$  is a compressing disc for
$\partial_hV$. Thus, if $\xi$ is tight then the horizontal boundary $\partial_hV$ of the supporting
manifold $V$ is incompressible.

\begin{lemma} \label{BetterTightLemma} Suppose $\xi$ is a pure contamination supported by $V$.  
Then $\xi$ is tight if and only if
there is no overtwisting disc $D$ with $\bdry D\embed \intr (V)$ and no 
compressing disc for $\partial_hV$. 
\end{lemma}

\begin{proof} We will show that an overtwisting disc $D$ whose boundary $\partial D$
intersects $\partial_hV$, but is not contained in it, can be always replaced
by another overtwisting disc $D'$ with $\partial D'\subset int(V)$, which clearly
implies the lemma.

For a disc $D$ as above, we can subdivide the boundary $\partial D$ into consecutive
segments $I_1,I_2,\dots,I_{2n}$ such that:

\item{(1)} the intersection of $\partial D$ with $\partial_hV$ is contained
in the union of segments $I_{2k}$ with even indices only;
\item{(2)} for each $1\le k\le n$ the segment $I_{2k}$ intersects $\partial_hV$
in one component only. 
\trip

\noindent We can then push each segment $I_{2k}$ slightly, in a direction
transverse to $\xi$, to get a legendrian
segment $I_{2k}'\subset \intr(V)$. If these modifications are slight enough, we can also
perturb the remaining segments $I_{2k-1}$, which all pass through $\intr(V)$, into
new legendrian segments $I_{2k-1}'\subset \intr(V)$, so that the union of all segments
$I_j':1\le j\le 2n$ forms a smooth closed legendrian curve in $\intr (V)$. It is clear
that this curve still bounds an overtwisting disc.  
\end{proof}

We now pass to the discussion of deformations (of pure contaminations) that preserve tightness.

\begin{definition}\label{PureDeformationAwayDef} A {\it pure deformation away from
$\partial_hV$} is a pure deformation $\xi_t$ of (pure) contaminations supported by $V$
such that for some open neighbourhood $U$ of $\partial_hV$ in $V$ the deformation
$\xi_t$ is constant on $U$.
\end{definition}

The next lemma is a well-known generalization of Gray's Theorem, and we will refer to it as the relative
Gray's Theorem. 

\begin{lemma}\label{RelativeGrayLemma} Let $\xi_t:t\in[0,1]$ be a pure deformation
away from $\partial_hV$. Then $\xi_t$ is an isotopy and in particular 
the pure contaminations $\xi_0$ and $\xi_1$ are contactomorphic by a diffeomorphism
isotopic to identity.
\end{lemma}

\begin{proof} An isotopy as in the conclusion of Gray's Theorem can be constructed as a flow with respect
to some variable vector field $X_t$ determined by $\xi_t$, see \cite{ABKLR:Symplectic}. This method
can be applied in our situation without any changes, giving a variable vector field
$X_t$ in $V$ which vanishes near $\partial_hV$. The flow with respect to $X_t$ 
is then equal to identity in the neighbourhood $U$ of $\partial_hV$ in $V$, thus
giving a 1-parameter family $f_t$ of diffeomorphisms of $V$. The fact that
$f_t$ is an isotopy, as required, follows by arguments described in \cite{ABKLR:Symplectic} for the proof of the
standard Gray's Theorem.
\end{proof}

The next definition and lemma generalize the above considerations slightly.

\begin{definition}\label{SigmaDeformationAwayDef} A {\it $\sigma$-deformation away from
$\partial V$} is a $\sigma$-deformation $\xi_t$ of contaminations supported by $V$
such that for some open neighbourhood $U$ of $\partial V$ in $V$ the deformation
$\xi_t$ is constant on $U$.
\end{definition}

Observe that if a $\sigma$-deformation $\xi_t$ as above starts at a pure contamination,
then the maximal integral surface is disjoint from $\partial V$ for all $t$,
i.e. none of its components is adjacent to $\partial_hV$ at the cusp locus.
The next lemma shows that tightness is invariant under such deformations.

\begin{lemma}\label{SigmaRelativeGrayLemma} Let $\xi_t:t\in[0,1]$ be a $\sigma$-deformation
away from $\partial_hV$ connecting pure contaminations $\xi_0$ and $\xi_1$. 
Then $\xi_0$ and $\xi_1$ are contactomorphic by a diffeomorphism
isotopic to identity. In particular, $\xi_0$ is tight iff $\xi_1$ is tight.
\end{lemma}

\begin{proof}
For a $\sigma$-deformation $\xi_t$ away from $\partial V$ it is possible to construct
a parametric family $\xi_{t,s}$ of combinations of controlled reductions, which  are
all away from $\partial_hV$, such that $\xi_{t,0}=\xi_t$ and $\xi_{t,1}$ is pure
for all $t\in[0,1]$. This can be done by using the same methods as in Sections \ref{Reduction} and \ref{DeformationTheorem}.
It follows that the pure contaminations
$\xi_0=\xi_{0,0}$, $\xi_{0,1}$, $\xi_{1,1}$ and $\xi_{1,0}=\xi_1$ are connected
by a pure deformation away from $\partial V$. We finish the proof by applying 
the relative Gray's Theorem (Lemma \ref{RelativeGrayLemma}). 
\end{proof}

We now turn to deformations that change the supporting manifold $V$
of a pure confoliation. To define a class of such deformations which we
will be interested in, recall that a controlled reduction is a kind
of $\sigma$-deformation, performed by modifying slope functions in
a given special chart (a $\cal C$-chart or a $\cal B$-chart), with one
discontinuity point for the maximal integral surface $F$, where a disc
component or a bridge (1-handle) is eliminated from $F$ (compare Section 2).
Recall also that in the context of $\sigma$-contaminations, a bridge
that is eliminated by a controlled reduction may be partly adjacent to the cusp 
locus of the supporting manifold $V$.

\begin{definition}\label{SpecialControlledDef}
An {\it elementary controlled splitting} is a deformation that turns a pure contamination
$\xi$ into another pure contamination $\xi'$, and consists of the following two steps.
The first step is an operation inverse to a controlled reduction (which in our case
either produces a disc component of the maximal integral surface in the interior
of $V$, or produces a bridge adjacent to the cusp locus of $V$ along two disjoint
arcs). A second step consists of cutting $V$ along the whole maximal integral
surface (produced in the first step), to get a new $V'$ and $\xi'$.
A {\it controlled splitting} is a finite sequence of elementary controlled splittings.
A {\it controlled pinching} is an operation inverse to a controlled splitting.
Finally, a {\it controlled deformation} is a combination of controlled
splittings and pinchings (in fact, it is a combination, in sequence, of elementary controlled
splittings and their inverses). 
\end{definition}

\begin{lemma}\label{SpecialControlledLemma}
Controlled deformations of pure contaminations that preserve the property of incompressibility
for horizontal boundaries of the supporting manifolds with inward cusps preserve tightness.
\end{lemma}

\begin{proof}
According to Lemma 11.1, preserving tightness is equivalent to preservoing the following
two conditions:
\item{(1)} nonexistence of an overtwisting disc $D$ with $\partial D\subset int(V)$;
\item{(2)} incompressibility of $\partial_hV$.

\noindent
Under assumptions of this lemma, condition (2) is preserved automatically.
It is then sufficient to show that condition (1) is also preserved, and this can be
further reduced to the case of a single elementary controlled splitting.
We will deal with this case using the methods of Section \ref{Invariance} as follows.

Suppose that $\xi_t:t\in[0,1]$ is the operation inverse to the controlled reduction
occurring as a first step in a given special elementary controlled splitting. Then
$\xi_0$ is a pure contamination and $\xi_1$ has the maximal integral surface $F_1$
of one of the two forms mentioned Definition \ref{SpecialControlledDef}.
Both $\xi_0$ and $\xi_1$ are supported by the same manifold with inward cusps $V$.
An argument as in Section \ref{Invariance}
shows that the existence of an overtwisting disc $D$ for $\xi_1$ with
$\partial D\subset int(V)\setminus F_1$ is equivalent to the existence of such a $D$
for $\xi_0$ with $\partial D\subset int(V)$. 

Let $\xi_1'$ be the pure contamination (supported by $V'$) obtained from $\xi_1$
by cutting $V$ along $F_1$, i.e. the one obtained from $\xi_0$ by our fixed
special elementary controlled splitting. From what was said above, it immediately
follows that the existence of an overtwisting disc $D$ for $\xi_1'$ with
$\partial D\subset int(V')$ is equivalent to the existence of such a $D$
for $\xi_0$ with $\partial D\subset int(V)$. This finishes the proof.
\end{proof}

We finish this section by remarking that we have not proved the invariance of tightness
for the following classes of deformations of pure contaminations:
\item{(1)} splittings that are not controlled;
\item{(2)} pinchings that are not inverse to appropriate splittings;
\item{(3)} splittings which do not preserve the incompressibility property of $\partial_hV$.

\trip
\noindent
More work is required to understand invariance of tightness for pure contaminations.

\section{Purification}\label{Purification}

In many situations, given a (positive) confoliation or contamination $\xi$ defined on $M$ or $V\embed M$, which is non-integrable
on some open set in its support, it is possible to find a deformation which converts $\xi$ to a contact structure or a pure
contamination.  Such a deformation is called a {\it purification}.  The idea of purification is one of the important ideas in
the monograph of Eliashberg and Thurston, see \cite{ET:Confoliations}, though they cannot be blamed for our terminology.  It would
be useful to establish an optimal result indicating when a contamination can be purified.  We will need to purify contaminations,
but only in very special situations.  In fact, the controlled reductions of Section
\ref{Reduction} constitute a special class of purifications.  In the paper \cite{UOJS:ContaminationCarrying}, we purify using $\cal
C$-charts ($\cal B$ charts), under the assumption that the contamination is pure adjacent to the $\cal C$-chart (the ends of the $\cal
B$-chart), just as in the case of controlled reductions.  

In Section \ref{Tightness} we will need to purify in a particular situation, which we
now describe.
Let $B\embed M$ be a branched surface with boundary,
properly embedded in a 3-manifold $M$, i.e. an embedding $(B,\partial B)\embed (M,\bdry M)$.  A fibered neighbourhood $V(B)$
is associated to $B$ as before, with $V(\partial B)$
the corresponding associated fibered neighbourhood of the boundary train track.
Let $\xi$ be a confoliation carried by $B$, i.e. a confoliation in $V(B)$
that is transverse to the vertical fibers of $V(B)$ and tangent to the horizontal
boundary $\partial_h V(B)$. Let $R\subset B$ be a rectangle embedded in $B$,
disjoint from the branch locus of $B$, with two opposite sides contained in
$\partial B$. Furthermore, let $\{(x,y):|x|\le1,|y|\le1\}$ be smooth coordinates for
$R$ so that the sides $|x|=\pm1$ are contained in $\partial B$. Given $\epsilon$, 
$0<\epsilon<1$, put $R_\epsilon=\{|y|<1-\epsilon\}\subset R$ and denote by
$V(R)=\pi^{-1}(R)$, $V(R_\epsilon)=\pi^{-1}(R_\epsilon)$ the corresponding
restrictions of the fibration $\pi:V(B)\to B$. Finally, suppose that 
for some $0<\epsilon<1$ the confoliation
$\xi$ is pure in $V(R)-V(R_\epsilon)$ (except at top and bottom where it
is tangent to $\partial_h V(B)$), while it can be arbitrary inside $V(R_\epsilon)$.

\begin{lemma} \label{PurificationLemma}
Under the assumptions described above, there exists a confoliation $\xi'$ 
carried by $B$ that coincides with $\xi$ outside $V(R)$ and is pure in $V(R)$.
\end{lemma}

\begin{proof}
We construct a confoliation $\xi$ as required in two steps:
we first construct a suitable variant of a $\cal B$-chart; then we modify
$\xi$ by changing the slope function inside the $\cal B$-chart. 
The kind of $\cal B$-chart we want to use is a smooth map $\psi: [-1,1]^3\to V(B)$
such that conditions (1) and (2) of Definition \ref{BChartDef} are satisfied, as are the following:

\item{(1)} the first two coordinates $x,y$ in $[-1,1]^3$ coincide with the
coordinates $x,y$ in $R$, in the sense that $\pi\circ\psi(x,y,z)=(x,y)$
in these coordinates;
\item{(2)} the chart image $\psi([-1,1]^3)$ coincides with the preimage
$\pi^{-1}(R)=V(R)$;
\item{(3)} the images by $\psi$ of the segments $\{x=\hbox{const},z=\hbox{const}\}
\subset [-1,1]^3$ are tangent to $\xi$.

\noindent
The existence of a chart as above is as obvious as the existence of the ordinary $\cal B$-charts described in section 2.

Let $f:[-1,1]^3\to \reals$ be the slope function for the
contamination
$\xi$ in the $\cal B$-chart coordinates as above.
We then have $\partial f/\partial y\ge 0$ and $\partial f/\partial y(x,y,z)>0$
for $|y|>1-\epsilon$ and $|z|<1$. Thus if $1>y_1>1-\epsilon$ then $f(x,y_1,z)>f(x,-y_1,z)$
for all $x,z\in(-1,1)$. It follows that there exists a smooth function
$f':[-1,1]^3\to \reals$ with the following properties:
\hop
\item{(1)} $f'$ coincides with $f$ in the union of the subsets $\{|y|\ge y_1\}$ and $\{|z|=1\}$ in $\cal B$;
\item{(2)} $\partial f'/\partial y>0$ in the subset $\{-1< z<1\}$.
\trip

Since the function $f'$ coincides with $f$ near the frontier of the chart image in $V(B)$,
it can be used as a slope function for a new contamination $\xi'$ which coincides
with $\xi$ outside the chart image. By Lemma \ref{BoxChartLemma}, $\xi'$ is contact in the subset
$\{ -1< z<1\}$.
\end{proof}

\section{Internal tightness}\label{Tightness}

In this section we deal with a property of pure contaminations weaker than tightness,
namely a property we call internal tightness. We find some conditions on a branched surface embedded in 3-manifold that
imply the internal tightness of any pure contamination carried by the branched surface.  As a corollary we get Theorem
\ref{TightnessTheoremOne}. In the next section we will show that in some cases it is possible to deduce the tightness
of a pure contamination from its internal tightness. This will make it possible to describe conditions on a branched surface
sufficient to imply tightness of any pure contamination carried by the branched surface.

\begin{definition}  A pure contamination $\xi$  supported by $V$ is {\it internally tight}
if it has no overtwisting disc contained in $\intr(V)$.
It is {\it internally universally tight} if there is no overtwisting disc contained 
in the universal cover $\intr(\tilde V)$) for the induced (lifted) 
contamination $\tilde\xi$.
\end{definition}

Clearly, universal tightness implies tightness since an overtwisting disc (if exists)
can be always lifted to the universal cover. In this section we
will prove this stronger property of internal universal tightness for some
classes of pure contaminations.

To formulate our main result in this section we introduce two concepts necessary for 
describing the assumptions. A closed (smooth) curve
$\alpha$ in the branch locus of a branched surface $B$ is {\it essential} if its embedding in $B$
is a $\pi_1$-injective map. Equivalently, $\alpha$ is essential if any of its lifts
to the universal cover $\tilde B$ of $B$ is a (non-closed) line. 

Before introducing a second concept, the concept of an immersed twisted disc of contact, we give 
some preliminary definitions.

\begin{figure}[ht]
\centering
\scalebox{1.0}{\includegraphics{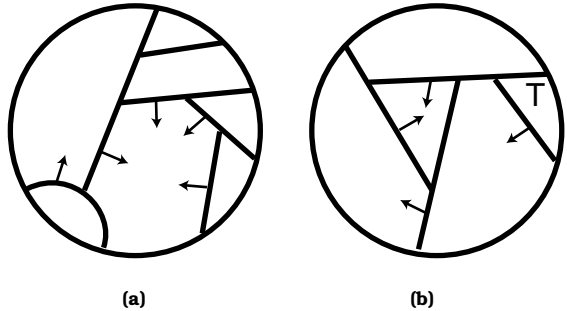}}
\caption{\small T-patterns on discs.}
\label{ContactTPattern}
\end{figure}

 Given a surface $E\embed\bdry_hV(B)$ for some branched
surface $B\embed M$, $\pi\from V(B)\to B$ defines a map of $E$ to $B$.  Pulling back the branch locus of $B$ to $E$, we
obtain a {\it T-pattern} on $E$.  It consists of a 1-complex of smooth curves possibly joined at {\it T-junctions}, see
Figure \ref{ContactTPattern}.  In addition to this 1-complex, we also consider as part of the T-pattern the data giving the sense
of branching on each smooth curve of the T-pattern.  This is indicated by an arrow transverse to the curve of the T-pattern
as shown in the figure.  The transverse arrow describes the sense of branching consistent with our definition, given in the introduction,
of inward and outward branching along the boundary of a sector.  Recall that if sectors
$W,X,Y$ are adjacent along an arc $\gamma$ of branch locus, and if
$W\cup Y$ and $X\cup Y$ are smooth, we say that branching along the arc
$\gamma\subset \bdry Y$ is {\it inward} for $Y$ and {\it outward} for $X$.  We indicate the sense of branching using a
transverse arrow pointing into
$Y$. 

T-patterns are defined in greater generality.  Let $B$ be a branched
surface and
$V(B)$ its fibred neighbourhood, and suppose that $i:E\to B$ is an immersion of a surface $E$, usually a disc in this paper.
Then $i$ determines an immersion $\hat i:E\to int(V(B))$ such that $i=\pi\circ\hat i$,
unique up to homotopy along the fibres of $V(B)$.  If $\hat i$ is an embedding in $V(B)$, then cutting $V(B)$ on $E$ yields a
neighborhood
$V(B')$ with two components of horizontal boundary $E_-$ and $E_+$ homeomorphic to $E$.  The T-patterns induced on $E_-$ and $E_+$ are
the {\it T-patterns} induced on the disc $E$ on each side of $E$ respectively.  Even when $\hat i$ is not an embedding, we can pull
back two T-patterns to $E$, for example by splitting locally.

\begin{figure}[ht]
\centering
\scalebox{1.0}{\includegraphics{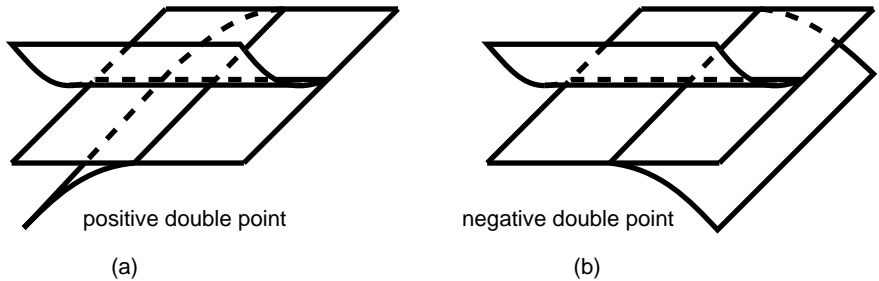}}
\caption{\small The sign of a double point of the branch locus.}
\label{ContactDoubleSign}
\end{figure}

Suppose the T-pattern on $\bdry_hV(B)$ includes a sequence of segments 
$\gamma_0,\gamma_1\ldots \gamma_n$ which can be
directed so the final end of $\gamma_i$ is attached to a point in the interior of $\gamma_{i+1}$, 
$\mod n$, forming a cycle that
bounds a disc (see Figure \ref{ContactTPattern}(b)).  We call this disc a {\it twisted disc of contact in the T-pattern}, 
which is {\it positive} ({\it
negative}) if the double points of the branch locus at the ends of the $\gamma_i$'s are all positive (negative).  Our
sign convention for double points of the branch locus is chosen for consistency with signs 
of twisted discs of contact as defined in the paper
\cite{UOJS:ContaminationCarrying}, see Figure \ref{ContactDoubleSign}.

We say that $B$ has an
{\it immersed twisted disc of contact} if there exists an immersion $i:E\to B$ of a disc $E$ such that
the induced T-pattern on one side has a twisted  disc of contact.

\begin{remark}
Deciding wheather a given branched surface has an immersed twisted disc of contact is
in general not easy. There are some obvious cases when nonexistence of these objects
is obvious. One of them is when a branched surface has no double points in the branched locus.

In our paper \cite{UOJS:ContaminationCarrying} we show how to decide whether a given branched surface
admits so-called immersed twisted
surfaces of contact. The nonexistence of these more general objects clearly implies
the nonexistence of immersed twisted discs of contact.
\end{remark}

\begin{thm} \label{TightnessTheorem}  Suppose a branched surface $B\embed M$ carries a positive pure contamination
$\xi$.  Suppose also that the branch locus of
$B$ consists of essential closed curves, $\partial_hV(B)$ contains no sphere component, and $B$ has
no immersed twisted disc of contact.  Then $\xi$ is internally universally tight.
\end{thm}

To prove above theorem we need a technical result which allows a closer study of overtwisting
discs entirely contained in $int(V(B))$. In this result we use a notion of a branched surface
with boundary. We assume that the boundary of a branched surface is a train track transverse to the branch locus
(and it avoids double points of the branch locus).  The horizontal boundary
of a neighbourhood of such a branched surface consists of components that are surfaces with corners,
where the corner points separate segments of the cusp locus and boundary segments (i.e.
segments that are projected by $\pi$ to the boundary of the branched surface), see Figure \ref{ContactBranchedModel}.

\begin{proposition}\label{ImmersedBranchedProposition}
Let $\xi$ be a pure contamination carried by a branched surface $B$, and suppose that $D\subset int(V(B))$
is an overtwisting disc for $\xi$. Then there is a compact branched surface $C$ with boundary and a map
$f:C\to B$ such that:
\item{(1)} $C$ is simply connected;
\item{(2)} $f$ is an immersion and when restricted to $int(C)$ it is a local homeomorphism;
\item{(3)} if $\bar\xi$ is a pullback of $\xi$ to $V(C)$ via $f$ then $D$ pulls back to an
overtwisting disc $\bar D\subset int(V(C))$ for $\bar\xi$.
\end{proposition}

\begin{proof}
Put the overtwisting disc $D$ in general position with respect to the projection map $\pi:V(B)\to B$.
By this we mean that the restriction of $\pi$ to $D$ is generic.  More precisely, the projection map has folds and Whitney folds, see
Figure \ref{ContactWhitney},
and the family of lines in $B$ consisting of the cusp locus of $B$ and the projection 
of the folds of $D$ and of $\partial D$ to $B$ is generic, meaning that self-intersections are at most double and transverse.
Consider then a structure of a cell complex on $B$ such that the 1-skeleton $B^{(1)}$ contains:
\item{(1)} projections of folds of $D$;

\item{(2)} the projection of $\bdry D$;

\item{(3)} the branch locus of $B$.

\trip\noindent
Pulling back, we obtain a cell structure for $D$.

\begin{figure}[ht]
\centering
\scalebox{1.0}{\includegraphics{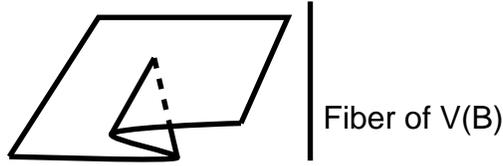}}
\caption{\small Whitney fold.}
\label{ContactWhitney}
\end{figure}

\noindent
In a cell structure as above, every cell in $D$ is mapped by $\pi$ bijectively to its image in $B$.
Modify the cell complex $D$ into a cell complex $Q$ by identifying recursively pairs of 2-cells that intersect 
(at least at a vertex) and have
the same projection to $B$. Use identifications that are compatible with projections to $B$. 
Repeat this in arbitrary order, till no pair of 2-cells as above is
left in the final quotient. Observe that the resulting complex $Q$ admits the induced from $\pi$
map $g:Q\to B$ which is an immersion (i.e. a local embedding). Observe also that the identifications we
perform produce no new generators in the fundamental groups of succesive quotients, and thus the final
complex $Q$ is simply connected, though not necessarily contractible. 

To get a branched surface $C$ as required, we want to enlarge $Q$ slightly as follows. If $g:Q\to B$
happens to be an embedding, we take $C$ to be a small regular closed neighbourhood of $g(Q)$ in $B$.
If $g$ is not an embedding, we construct a small regular neighbourhood of the image locally, so that
we get kind of a regular neighbourhood of an immersed object. In either of the cases we obtain
a branched surface $C$ with boundary for which $Q$ is a deformation retraction. We obtain also
an extension $f:C\to B$ of $g$ which, by the local nature of the construction, is an immersion
and, at each interior point of $C$, a local homeomorphisms. 

It remains to prove part (3) of the proposition. View the quotient map from $D$ to $Q$ as a map
$q:D\to C$, and observe that, viewing $D$ as a subset of $V(B)$, we have $f\circ q=\pi$.
On the other hand, we can identify the fibred neighbourhood $V(C)$ with the space
$$
\{ (x,t)\in C\times V(B):f(x)=\pi(t)  \}.
$$ 
Then the pull-back $\bar D$ of $D$ in $V(C)$ is the image of the map $\varphi:D\to V(C)$
defined by $\varphi(u)=(q(u),u)$. 
\end{proof}

One more property of a branched surface $C$ as above is described in the following.

\begin{lemma}\label{SpheresLemma}
If $C$ is an orientable simply connected compact branched surface with boundary then the
boundary $\partial V(C)$ consists of 2-spheres.
\end{lemma}

\begin{proof}
It is known that if $M$ is a compact orientable 3-manifold then the image of the boundary
map $H_2(M,\partial M)\to H_1(\partial M)$ has rank equal to one half of the rank of $H_1(\partial M)$.
Applying this to $M=V(C)$, which is simply connected, we observe that the rank of the boundary map 
is zero in this case. Therefore the rank of $H_1(\partial V(C))$ is also zero, hence the lemma.
\end{proof}

The proof of Theorem \ref{TightnessTheorem} requires further preliminaries, but we give a brief outline here. 
We will work with a branched surface $C$ as in Proposition \ref{ImmersedBranchedProposition}, containing
the overtwisting disc $\bar D$ for $\bar\xi$ in $int(V(C))$. Using a sequence of modifications
away from $\bar D$ we will convert $\bar\xi$ to a pure contamination supported by a disc
(and still containing $\bar D$ in the interior). 
This will be a contradiction, due to the following fact.

\begin{fact}\label{DiscCarriedTightFact}
A pure contamination carried by a disc is internally tight.
\end{fact}

The above fact is well known and can be easily deduced, e.g., from Proposition 3.5.6
in \cite{ET:Confoliations}.

\begin{figure}[ht]
\centering
\scalebox{1.0}{\includegraphics{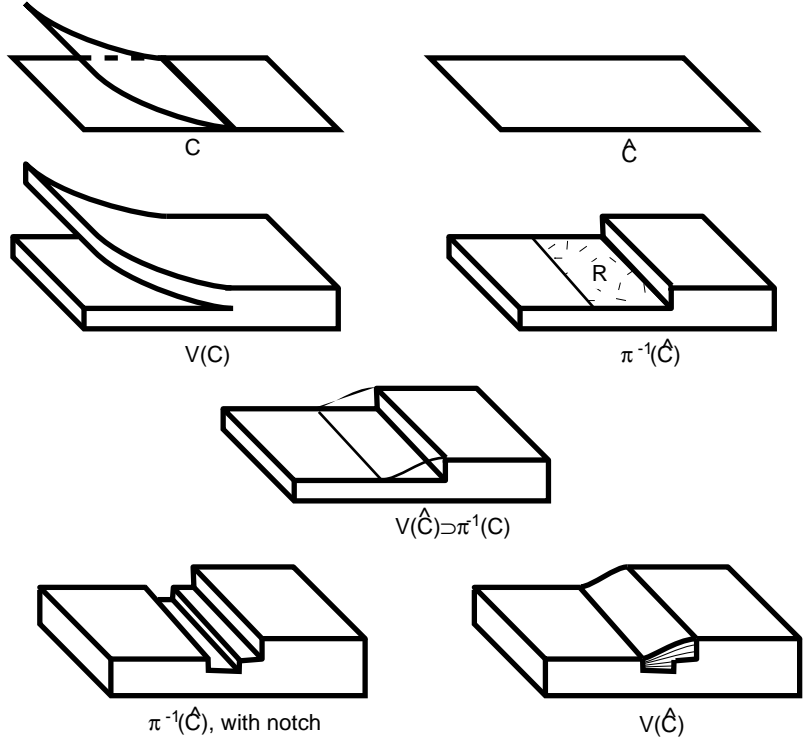}}
\caption{\small Eliminating a step.}
\label{ContactStep}
\end{figure}

We will now describe a basic modification we will be applying to some parts of $\bar\xi$.

\begin{lemma} \label{StepRemovalLemma}  Suppose $C$ is the branched surface shown in Figure 
\ref{ContactStep}, and $\hat C$ is the disc as shown, obtained from $C$ by deleting one of the sectors. 
Let $\alpha$ be the cusp arc of $V(C)$ and $R\subset\partial_hV(C)$ a rectangle
adjacent to $\alpha$, as shown in the same figure.  Suppose
$C$ carries  a pure contamination $\xi$, supported in $V(C)$.  
Identify $\pi\inverse(\hat C)$ with a subset of $ V(\hat C)$ as shown in the figure. Then there is a pure contamination $\hat\xi$
carried by $\hat C$, supported by $V(\hat C)$, such that
\item{(1)} 
$\bdry_hV(C)\cap\pi\inverse(\hat C)-R\subset \bdry_hV(\hat C)$;  
\item{(2)} 
$\hat \xi$ is the same as $\xi$ in $\pi\inverse(\hat C)$ except in an arbitrarily small neighborhood in $\pi\inverse(\hat C)$
of the rectangle $R$.
\end{lemma}

\begin{proof}  It will be convenient to choose coordinates describing $\pi\inverse(\hat C)$.  In $\reals^3$
$\pi\inverse(\hat C)$ will be the set \hfil \break
\centerline{$\{ (x,y,z): -2\le x\le 2,\ 0\le z\le 1,\ -4\le y \le 4 \}\cup$}\break
\centerline{$ \{ (x,y,z): -2\le
x\le 2,\ 1\le z\le 2,\ 1\le y \le 4 \}$.}\break
We have a confoliation $\xi$ defined on this set, which is contact except on the horizontal portions of the boundary,
where it is integrable. We choose a line segment $\beta$ at $y=b,\ z=1$, with $b<1$ arbitrarily close to the cusp arc
$\alpha$ which is at
$y=1,\ z=1$.  The arc $\beta$ is in the horizontal boundary, where $\xi$ is integrable.  We choose a Legendrian $\gamma$
whose interior is in the interior of $\pi\inverse(\hat C)$, and which satisfies $\pi(\gamma)=\pi(\beta)$.  The Legendrian
$\gamma$ can be chosen arbitrarily close to the horizontal boundary at $z=1$, but note that $\gamma$ need not lie on a
level $z=\hbox{constant}$.  Using a $\sigma$-deformation, we replace $\gamma$ by an integrable strip $G$ whose projection
in the $xy$ plane contains $b-\epsilon\le y\le b+\epsilon,\ -2\le x\le 2$, for some small $\epsilon$, $0<\epsilon<1-b$. 
We still use $\xi$ to denote the resulting confoliation.  Now we remove the portion of the slice $b-\epsilon\le y\le
b+\epsilon$ above the integral strip
$G$ in $\pi\inverse(\hat C)$, and we consider the restriction of $\xi$ to this set, which is $\pi\inverse(\hat C)$ with
a notch cut into it.  There is an induced foliation from
$\xi$ on all the vertical sides of the notched $\pi\inverse(\hat C)$.

Now, we want to insert a foliated set to construct
$V(\hat C)$, with the foliation chosen to give an extension $\xi'$ of $\xi$ as shown in Figure \ref{ContactStep}, lower right. However,
the
$\xi'$ obtained is in general not smooth (only continuous) at the vertical sides of the notched $\pi^{-1}(\hat C)$.
Thus, before inserting a foliated part, we need first to modify $\xi$ slightly,
so that it becomes a foliation near the vertical sides of the notched $\pi^{-1}(\hat C)$.
This can be done using appropriate charts and modifying slope functions. After this
preparation, and after inserting a foliated part as before,
we get a smooth $\xi'$ in $V(\hat C)$ which satisfies the assertions,
except that it is not pure. 

As a final step, we use Lemma \ref{PurificationLemma} to purify $\xi'$ and obtain a pure contamination
$\hat\xi$ supported in
$V(\hat C)$ and carried by $\hat C$.
\end{proof}

\begin{proof}[Proof of Theorem \ref{TightnessTheorem}] Suppose $B\embed M$ is as in the statement of the theorem,
with $\xi$ a pure contamination carried by $B$ and supported in $V(B)$.  Then the universal cover $\tilde B$ has a neighborhood
$\tilde V(B)$ which supports a pure contamination $\tilde \xi$.   We will work in
$\tilde V(B)$.  Assuming there is an overtwisting disc $D$ for $\tilde \xi$ in $\intr(\tilde V(B))$, we will obtain a
contradiction.  

Let $C$, $\bar\xi$ and $\bar D$ be objects as in Proposition \ref{ImmersedBranchedProposition},
for the branched surface $\tilde B$, the pure contamination $\tilde\xi$ and the overtwisting disc $D$.  
Our goal is to perform modifications on $C$ and $\bar\xi$ to  put $\bar D$ into a product
$E\times I$, as an overtwisting disc of a modified $\bar\xi$, which is a pure contamination carried by the disc $E$ and
supported by $V(E)=E\times I$.  

By the assumptions of this theorem, the branched surface $C$ has no sphere components of $\partial_hV(C)$,
no immersed twisted discs of contact, and no closed curve in the cusp locus of $V(C)$. The latter follows from the fact
that each cusp curve in the initial branched surface $B$ is essential.

We make the following claim:
\trip
\noindent{\it Claim:}  Suppose $C$ is as in the conclusion of Proposition \ref{ImmersedBranchedProposition}.  Suppose also that
$C$ has no sphere components of $\partial_hV(C)$,
no immersed twisted discs of contact, and no closed curve in the cusp locus of $V(C)$.  Then there is a branch
arc
$\alpha$ in
$C$  cutting a half-disc sector $H$ from $C$ such that the sense of
branching is outward from $H$ and $\partial H-\alpha\subset\partial C$.
\trip

\noindent We will call such an arc $\alpha$ {\it innermost}. 
Observe that the cut half-disc
$H$ can be lifted in a unique way to a half-disc component of $\partial_h V(C)$ which is
adjacent to $\alpha$ viewed as a cusp curve in $V(C)$.
Accepting the claim for the moment, 
we perform a pinching of $V(C)$ identifying
a portion of $H$ with a small neighborhood in $\bdry_hV(C)-H$ of
$\alpha$ (regarded as a cusp curve of $V(C)$).  Clearly there is such a pinching which ensures that
after the pinching,
$H$ is smaller, and $D$ is disjoint from $\pi\inverse(H)$.  We perform a $\sigma$-deformation to purify $\xi$.
Next, we let
$\hat C$ be the branched surface consisting of all sectors of $C$ except $H$.  Then
$\pi\inverse(\hat C)$ has a step corresponding to the branch curve $\alpha$. Using Lemma
\ref{StepRemovalLemma} we modify $\xi$ to obtain $\hat \xi$ carried by $\hat C$, and $\hat C$ has fewer branch curves
than $C$.   Note that all modifications could be done without disturbing the overtwisting disc
$D$.  Now we can repeat the process after first renaming as $C$ the branched surface $\hat C$, and renaming as $\xi$
the pure contamination $\hat \xi$.  The new $C$ has the properties of the old one; in particular, it has no immersed twisted
disc of contact, since the branch curves of the new $C$ are also branch curves of the old $C$.

After finitely many modifications, we must obtain a branched surface $C$ without branch arcs, which is
still simply connected. Thus $C$ is the required disc
$E$, and we have our overtwisting disc in $\intr(V(E))$, contradicting Fact \ref{DiscCarriedTightFact}.

It remains to prove the claim used in the above argument, and we devote the remainder of the section to this aim.
\end{proof}

Under the hypotheses on $C$ of the Claim we prove the following

\begin{lemma}\label{ZeroGonLemma}
For each component $S$ of $\partial V(C)$ there there are at least two disc components $\Delta$ of the
horizontal boundary $\partial_hV(C)$, contained in $S$, with $\partial\Delta$ disjoint from the cusp locus of $V(C)$.
\end{lemma}

\begin{proof}
It follows from Lemma \ref{SpheresLemma} that each component $S$ of $\partial V(C)$ is a 2-sphere.
Moreover, no such $S$ is a component of $\partial_hV(C)$. Thus on each $S$ we have a
nonempty train-track neighbourhood corresponding to $\pi^{-1}(\partial C)\subset\partial V(C)$
and a family of cusp segments of $V(C)$ (possibly empty) connecting cusp points of the train-track
neighbourhood in pairs, passing through the complement of the train-track neighbourhood.  (Recall that there are no closed cusp curves
in
$V(C)$.)

Consider the complementary regions of the train-track neighbourhood in $S$. Since $S$ is the
2-sphere, an Euler characteristic argument shows that there are at least two complementary regions $\Delta$
being discs with smooth boundary (i.e. $\partial\Delta$ contains no cusp point of the train-track
neighbourhood). Since there are thus no cusp curves inside $\Delta$, it follows that each of $\Delta$ is a component
of the horizontal boundary $\partial_hV(C)$ as required.
\end{proof}

We will call a component $\Delta$ as in the above lemma a {\it $0$-gon component} in $\partial_hV(C)$.

\hop

We say a region in the complement of a T-pattern (on a component of $\bdry_hV(B)$) is {\it minimal} 
if all branching at the boundary is outward.  

\begin{lemma}\label{MinimalRegionLemma}  If a T-pattern on a disc contains no twisted discs of
contact, then the T-pattern contains a minimal region.
\end{lemma}

\begin{proof}  Suppose that the T-pattern on a disc $D$ contains no minimal region.  
Start in any complementary region $R_1$ of the
T-pattern; there must be at least one transverse orientation on an edge of 
the T-pattern pointing into $R_1$, say from $R_2$.  Similarly,
there must be at least one transverse orientation on an edge of the T-pattern 
pointing into $R_2$, from $R_3$ say, etc.  We conclude that
there is a simple closed oriented path $\gamma$ intersecting edges of 
the T-pattern with orientation agreeing with the transverse
orientation on the edge, and intersecting the T-pattern non-trivially.  Now
$\gamma$ bounds a T-pattern on a smaller disc $H$ bounded by $\gamma$.  
Choose any endpoint of the T-pattern on $\bdry H$ and suppose
that the transverse orientation points to the right.  Follow the edge starting 
at this point until it ends at a T-junction, i.e. ends as
the upstroke of the T-junction.  Turn right at this T-junction and follow 
the new edge until it ends at another T-junction, etc.  We
must obtain a cycle of edges, each starting at a horizontal stroke of a T and 
extending to the right.  Since the T-pattern is on a disc,
the cycle represents a twisted disc of contact.  We obtain a twisted disc of contact of
the opposite sign if the transverse orientation at $\bdry H$ points to the left.
\end{proof}

To prove the existence of an innermost arc for $C$ (provided the branch locus in not empty),
consider any component $S$ of the boundary $\partial V(C)$. According to Lemma \ref{ZeroGonLemma},
there is a $0$-gon component $\Delta$ in $S$. The T-pattern on $\Delta$ is then nonempty (otherwise $C$
is already a disc) and contains no twisted disc of contact (otherwise there would be an immersed twisted
disc of contact in the initial branched surface $B$). It follows from Lemma \ref{MinimalRegionLemma}
that this T-pattern contains a minimal region $R$. If $R$ is a half-disc
(i.e. contains only one arc of the T-pattern in its boundary), then the segment of T-pattern
in its boundary is an innermost arc as required. If not (i.e. if $R$ has several arcs
of the T-pattern in its boundary), then we use the following inductive argument.   

\begin{figure}[ht]
\centering
\scalebox{1.0}{\includegraphics{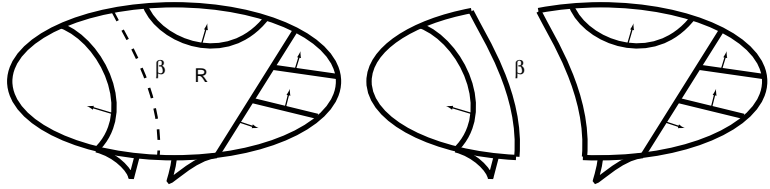}}
\caption{\small Cutting $C$ on an arc $\beta$ in a minimal region.}
\label{ContactCut}
\end{figure}

Cut $C$ along an arc $\beta$ corresponding to an arc in the minimal region $R$ connecting two distinct
segments in $\partial R$ not in the T-pattern, see Figure \ref{ContactCut}. Observe that, since $C$ is simply connected, $\beta$
disconnects $C$ yielding two simply connected branched surfaces 
(after necessary smoothing near the ends of $\beta$) that are simpler than $C$ (they contain fewer
segments in the cusp locus). Choose either of the two branched surfaces and denote it by $C'$.
Observe that after cutting along $\beta$ the part of the 0-gon component $\Delta$ contained in 
$\partial_hV(C')$ is a 0-gon component $\Delta_0$. Denote by $\Delta'$ a different 0-gon component
in $\partial_h V(C')$, which always exists due to Lemma \ref{ZeroGonLemma}. Note then that $\partial\Delta'$
does not contain a segment corresponding to the cutting arc $\beta$, because the component of $\partial_hV(C')$
other than $\Delta_0$ containing such a segment also contains a segment of the cusp locus in its boundary.

We repeat the previous steps, looking for a minimal region $R'$ of the T-pattern on $\Delta'$.
If there is such an $R'$ which is a half-disc, it yields an innermost arc for $C'$, which is also an
innermost disc for $C$ (becouse $\Delta'$ is far from $\beta$). If no such $R'$ is a half-disc,
we cut $C'$ as before, along an arc $\beta'$ in $R'$. The only difference is that now we choose for $C''$
the part of $C'$ which does not contain the previous cutting curve $\beta$.

After a finite number of steps as above, we must end with an innermost arc for the final branched surface of our sequence.  This arc is
also an innermost arc for the initial branched surface $C$, by the choices we make after each cuttting.

This completes the arguments in the proof of Theorem \ref{TightnessTheorem}.

\section{Tightness of contaminations}\label{ContaminationsTightness}

In this section we describe properties of a branched surface $B$ embedded in a 3-manifold
$M$ sufficient to guarantee that any
pure contamination carried by $B$ is tight. These conditions are based on a similar conditions sufficient for internal tightness
given in the previous section.  We view the results in this section as promising first steps towards our 
planned future development of tight pure contaminations as tools in 3-manifold topology.

First we need to recall the definition of incompressibility for a (two-sided) surfaces $S$ properly embedded in a 3-manifold $M$.
A {\it} compressing disc for a surface $S$ in a 3-manifold $M$ is a disc
$D$ embedded in $M$, with $D\cap S=\partial D$ and such that $\partial D$
does not bound a disc embedded in $S$. A surface $S$ properly embedded in $M$ is {\it incompressible}, 
if there is no compressing disc for $S$. We will apply this notion either
to closed surfaces $S$ in closed manifolds $M$, or to the boundaries $S=\partial M$
of compact manifolds $M$ with boundary.  

We start with a lemma showing that in some situations it is possible
to deduce the tightness of a pure contamination from its internal tightness.  

\begin{lemma}\label{InternalImpliesTightLemma}
Let $\xi$ be a pure contamination carried by a branched surface $B\embed M$,
and suppose that $\xi$ is internally tight. Suppose also that $\partial V(B)$ is an incompressible surface in $M$ and that
$\bdry_h V(B)$ has no disc components.  Then $\xi$ is tight.
\end{lemma}

\begin{proof}
According to lemma \ref{BetterTightLemma}, we have to exclude existence of
two types of overtwisting discs $D$: the ones with $\partial D\subset\intr(V(B))$
and the ones with $\partial D\subset\partial_h V(B)$. By the incompressibility of
$\partial V(B)$ in $M$, if there is an overtwisting disc $D$ for $\xi$ with
$\partial D\subset\intr(V(B))$, it can be replaced by a disc $D'\subset\intr(V(B))$ 
which coincides with $D$ on a neighbourhood of $\partial D=\partial D'$.
Since $D'$ is still an overtwisting disc, this contradicts internal tightness
of $\xi$. Suppose now that $D$ is an overtwisting disc for $\xi$ with
$\partial D\subset\partial_h V(B)$. By transversality of $D$ to $\partial_h V(B)$
at $\partial D$, a neighbourhood of $\partial D$ in $D$ intersects
$\partial V(B)$ only at $\partial D$. Again, by incompressibility of
$\partial V(B)$ in $M$, we can replace $D$ by a disc $D'$ which coincides with $D$
on a neighbourhood of $\partial D=\partial D'$ and which intersects
$\partial V(B)$ at $\partial D'$ only.  This is a standard argument using the Loop Theorem.  Since $D$ is an overtwisting disc,
the boundary $\partial D'=\partial D$ does not bound a disc in $\partial_h V(B)$,
and because $\bdry_hV(B)$ has not disc components, is a compressing disc for $\partial V(B)$ as well,
a contradiction. 
\end{proof}

Combining the above lemma with Theorem \ref{TightnessTheorem} we get Corollary \ref{ContaminationTightnessCorollary},
which can be used to construct examples of tight contaminations.

\begin{example} \label{TightExample} The branched surface $B$ shown in Figure 
\ref{ContactFavoriteBranched} can be embedded in $\reals^3$, though we
have drawn it as an immersed branched surface.  By the carrying criterion given in \cite{UOJS:ContaminationCarrying}, $B$
carries a positive pure contamination $\xi$.

\begin{figure}[ht]
\centering
\scalebox{1.0}{\includegraphics{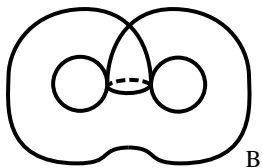}}
\caption{\small The branched surface $B$.}
\label{ContactFavoriteBranched}
\end{figure}

Regarding $B$ as embedded in $\reals^3$, we have a well-defined neighborhood $V(B)$, and clearly $\bdry_hV(B)$ consists
of two punctured tori which we label $F_+$ and $F_-$.  We can glue a 3-manifold with incompressible boundary to $\bdry
V(B)$ to obtain a closed 3-manifold $M$. 

Observe now that due to Theorem \ref{TightnessTheoremOne}
the pure contamination $\xi$ is internally tight. To get tightness by applying Lemma
\ref{InternalImpliesTightLemma} we need to show that the boundary $\partial_h V(B)$
is incompressible in $M$.

By construction of $M$, $\bdry V(B)$ is incompressible on the outside of $V(B)$.  
To prove incompressibility on the inside of $V(B)$, we argue as follows.  If
$\alpha$ is the branch locus in $B$, $\alpha$ a closed curve, then $\pi\inverse(\alpha)$ is an annulus $A$ whose core is the cusp
curve ${\cal C}(B)\subset V(B)$, cutting $A$ into two annuli $A_1$ and $A_2$. Cutting $V(B)$ on $A$ yields a product.   If $E$ were a
compressing disc for
$\bdry V(B)$ in
$V(B)$,  then making
$E$ transverse to $A$ and eliminating (trivial) closed curves of intersection, the intersection pattern would consist of arcs. 
Furthermore, we can assume these arcs are essential in $A_1$ and $A_2$, otherwise we could surger to obtain a compressing disc
intersecting $A$ in fewer arcs.   An innermost arc would cut a half-disc
$H$ from
$E$, with one arc of
$\bdry H$ in the vertical boundary of the product and one arc in the horizontal boundary of the product.  Since this is impossible, we
conclude that $\bdry V(B)$ is incompressible in $M$.
Thus, by applying Lemma \ref{InternalImpliesTightLemma}, we conclude that $\xi$ is tight.
\end{example}

Using Theorem  \ref{TightnessTheorem} it is easy to construct more tight pure contaminations in 3-manifolds, starting with branched
surfaces $B$ satisfying the assumptions of this theorem, and then glueing a manifold with incompressible boundary
to $\partial V(B)$.  The
unsatisfactory aspect of these examples, however, is that they can be constructed only in Haken manifolds. 
\bibliographystyle{amsplain}
\bibliography{ReferencesUO2}

\end{document}